\documentclass[hidelinks, 12pt]{amsart}               

\usepackage{amssymb, amsmath, amsthm,stmaryrd}
\usepackage[dvipsnames]{xcolor}
\usepackage[backref]{hyperref}
\usepackage[backrefs,lite]{amsrefs}
\usepackage{amscd}   % for commutative diagrams
\usepackage{fullpage}
\usepackage[all]{xy} % for complicated commutative diagrams
\usepackage{MnSymbol,graphicx}
\usepackage{mathrsfs}
\usepackage{makecell}
\usepackage{array}
\usepackage{tikz}
\usepackage{verbatim}
\usepackage[toc]{appendix}
\usepackage{amssymb, amsmath, amsthm,stmaryrd}
\usepackage[dvipsnames]{xcolor}
\usepackage[backref]{hyperref}

\usepackage{amscd}   % for commutative diagrams
\usepackage{fullpage}
\usepackage[all]{xy} % for complicated commutative diagrams
\usepackage{MnSymbol,graphicx}
\usepackage{mathrsfs}
\usepackage{makecell}
\usepackage{verbatim}
\usepackage{tikz-cd}
\usepackage{lipsum}
\usepackage{adjustbox}
  \tikzstyle{block} = [rectangle, draw,
      text width=8em, text centered, minimum height=4em]
  \tikzstyle{line} = [draw, -latex']
\usepackage{amssymb, amsmath, amsthm,stmaryrd}
\usepackage[backref]{hyperref}
\usepackage{amscd}   % for commutative diagrams
\usepackage{fullpage}
\usepackage[all]{xy} % for complicated commutative diagrams
\usepackage{MnSymbol,graphicx}
\usepackage{mathrsfs,enumitem}
\usepackage{multirow,qtree}% http://ctan.org/pkg/multirow
\DeclareFontEncoding{OT2}{}{} % to enable usage of cyrillic fonts

\usepackage{tikz}
\usepackage[OT2,T1]{fontenc}
\usepackage{mathtools}
\usepackage[normalem]{ulem}
\newtheorem{lemma}{Lemma}[section]
\newtheorem{theorem}[lemma]{Theorem}
\newtheorem{proposition}[lemma]{Proposition}
\newtheorem{prop}[lemma]{Proposition}
\newtheorem{cor}[lemma]{Corollary}

\newtheorem{claim*}{Claim}
\newtheorem{thm}[lemma]{Theorem}
\newtheorem{defn}[lemma]{Definition}
\theoremstyle{definition}
\newtheorem{remark}[lemma]{Remark}

\theoremstyle{definition}

\newtheorem{example}[lemma]{Example}

\newtheorem{notation}[lemma]{Notation}
\newtheorem*{theorem*}{Theorem}

\def\O{\mathcal{O}}

\newcommand{\F}{{\mathbb F}}

\newcommand{\Q}{{\mathbb Q}}

\newcommand{\Z}{{\mathbb Z}}

\newcommand{\Kbar}{{\overline{K}}}

\newcommand{\cO}{{\mathcal O}}

\newcommand{\mathfrakq}{{\mathfrak q}}
\newcommand{\mathfrakr}{{\mathfrak r}}

\newcommand{\mathfrakt}{{\mathfrak t}}

\newcommand{\s}{\mathfrak{s}}

\newcommand{\Qplus}{\mathcal{K}}
\newcommand{\sfL}{L}
\newcommand{\Qrplus}{\Q_{r}(\omega)}

\newcommand{\ie}{\textit{i.e.}}
\newcommand{\Crrp}{C_{r}}
\newcommand{\frrp}{f_{r}}
\newcommand{\gminus}{g_{r}^{-}}
\newcommand{\gplus}{g_{r}^{+}}
\newcommand{\Cminus}{C_{r}^{-}}
\newcommand{\Cplus}{C_{r}^{+}}
\newcommand{\Qzrplus}{\Qplus}
\newcommand{\Rroots}{\mathcal{R}}
\newcommand{\cR}{\mathcal{R}}

\newcommand{\qQzrplus}{\mathfrak{q}}

\newcommand{\rQzrplus}{\mathfrakr}

\DeclareMathOperator{\HH}{H}

\DeclareMathOperator{\lcm}{lcm}

\DeclareMathOperator{\Frob}{Frob}

\DeclareMathOperator{\End}{End}

\DeclareMathOperator{\Aut}{Aut}
\DeclareMathOperator{\Gal}{Gal}

\DeclareMathOperator{\Norm}{Norm}

\DeclareMathOperator{\Trace}{Tr}

\DeclareMathOperator{\Jac}{Jac}

\DeclareMathOperator{\GL}{GL}

\DeclareMathOperator{\disc}{disc}

\newcommand{\conductor}{N}

\numberwithin{equation}{section}
\numberwithin{table}{section}

\newcommand{\zr}{\zeta_{r}}

\newcommand{\Zl}{\Z_{\ell}}
\newcommand{\Ql}{\Q_{\ell}}
\newcommand{\Fl}{F_{\ell}}
\newcommand{\Flam}{F_{\lambda}}
\newcommand{\Slam}{S_{\lambda}}
\newcommand{\codim}{\operatorname{codim}}

\newcommand{\Tl}{\operatorname{T}_{\ell}}
\newcommand{\TlA}{\operatorname{T}_{\ell}(A)}
\newcommand{\Vl}{\operatorname{V}_{\ell}}
\newcommand{\Vlam}{\operatorname{V}_{\lambda}}
\newcommand{\rhol}{\rho_{\ell}}

\newcommand{\rholam}{\rho_{\lambda}}
\newcommand{\rholamtilde}{\widetilde{\rho}_{\lambda}}

\newcommand{\Gltwo}{\GL_{2}}
\newcommand{\Pplaces}{\mathcal{P}}
\newcommand{\WD}{\mathrm{WD}}
\newcommand{\condv}[1]{\mathfrak{n}_{v} (#1)}

\newcommand{\condexp}[1]{\mathfrak{n} \left(#1\right)}
\newcommand{\condtame}[1]{\mathfrak{n}_{\operatorname{tame}}(#1)}
\newcommand{\condwild}[1]{\mathfrak{n}_{\operatorname{wild}}(#1)}

\newcommand{\unr}[1]{#1^{\operatorname{unr}}}

\newcommand{\lastroot}{\gamma}

\newcolumntype{C}[1]{>{\centering\let\newline\\\arraybackslash\hspace{0pt}}m{#1}}

\usepackage{pbox}
\usepackage{tkz-graph}
\usetikzlibrary{snakes,fit,positioning,calc,shapes}
\usepackage{graphicx}
\usepackage{cleveref,tikz-cd,pbox,wrapfig}

\definecolor{amethyst}{rgb}{0.6, 0.4, 0.8}
\definecolor{atomictangerine}{rgb}{1.0, 0.6, 0.4}
\definecolor{deeppeach}{rgb}{1.0, 0.8, 0.64}
\definecolor{eggshell}{rgb}{0.94, 0.92, 0.84}
\definecolor{lightapricot}{rgb}{0.99, 0.84, 0.69}
\definecolor{lemonchiffon}{rgb}{1.0, 0.98, 0.8}
\definecolor{roundabout}{rgb}{1.0, 0.91, 0.75}
\definecolor{atomictangerine}{rgb}{1.0, 0.6, 0.4}
\definecolor{ruby}{rgb}{0.88, 0.07, 0.37}
\definecolor{sapphire}{rgb}{0.03, 0.15, 0.4}

\def\rootsep{0.05}               % horizontal space between roots in a cluster picture
\def\clustersep{0.06}            % horizontal and vertical space between parent & child cluster
\def\cnamescale{0.8}             % cluster names font size 
\def\cdepthscale{0.7}            % cluster depths and signs font size 
\def\cltopskip{1pt}              % space above a cluster picture
\def\clbottomskip{1pt}           % space below a cluster picture

\def\rootscale{0.5}   \def\rootcolor{gray}
\def\rootscaleA{0.7}  \def\rootcolorA{yellow}
\def\rootscaleB{0.5}  \def\rootcolorB{green}
\def\rootscaleC{0.4}  \def\rootcolorC{sapphire}
\def\rootscaleD{0.45}  \def\rootcolorD{ruby}
\tikzset{
  clA/.style = {very thick,black},
  clB/.style = {thick,purple}
  clC/.style = {thick, red}
}

\tikzset{
  root/.style = {circle,scale=\rootscale,fill=\rootcolor},
    rc/.style 2 args = {right=#1*1.5*\clustersep of {#2.east|-first},root}, rr/.style = {right=\rootsep of {#1.east|-first},root},
  roott/.style = {circle,inner sep=-2pt,minimum size=5pt,black,font=\ttfamily\footnotesize},
    rct/.style 2 args = {right=#1*1.5*\clustersep of {#2.east|-first},roott}, rrt/.style = {right=\rootsep of {#1.east|-first},roott},
  rootA/.style = {circle,scale=\rootscaleA,ball color=\rootcolorA},
    rcA/.style 2 args = {right=#1*1.5*\clustersep of {#2.east|-first},rootA}, rrA/.style = {right=\rootsep of {#1.east|-first},rootA},
  rootB/.style = {circle,scale=\rootscaleB,ball color=\rootcolorB},
    rcB/.style 2 args = {right=#1*1.5*\clustersep of {#2.east|-first},rootB}, rrB/.style = {right=\rootsep of {#1.east|-first},rootB},
  rootC/.style = {diamond,scale=\rootscaleC,ball color=\rootcolorC},
    rcC/.style 2 args = {right=#1*1.5*\clustersep of {#2.east|-first},rootC}, rrC/.style = {right=\rootsep of {#1.east|-first},rootC},
  rootD/.style = {circle,scale=\rootscaleD,ball color=\rootcolorD},
    rcD/.style 2 args = {right=#1*1.5*\clustersep of {#2.east|-first},rootD}, rrD/.style = {right=\rootsep of {#1.east|-first},rootD},
   rootE/.style = {circle,scale=0.1,fill=black},
    rcE/.style 2 args = {right=#1*1.5*\clustersep of {#2.east|-first},rootE}, rr/.style = {right=\rootsep of {#1.east|-first},rootE},
  cluster/.style = {draw=black!90,thick,rounded corners,inner sep=22*\clustersep,outer xsep=22*\clustersep,fit=#1},
   clusterc/.style = {draw=white,thick,rounded corners,inner sep=22*\clustersep,outer xsep=22*\clustersep,fit=#1},
  clabel/.style  = {anchor=west,scale=\cdepthscale,black,inner sep=0,outer xsep=1,outer ysep=0},
  clabelL/.style = {above right=-\clustersep of #1t.north east,clabel},
  clabelD/.style = {below right=-\clustersep of #1t.south east,clabel},
  clouter/.style = {inner sep=0,outer sep=0,fit=#1}
}

\def\Cluster #1 = #2;{\node[cluster=#2] (#1) {};}
\def\ClusterL #1[#2] = #3;{
  \node[cluster=#3] (#1t) {}; \node[clabelL=#1] (#1l) {$#2$}; \node[clouter=(#1t)(#1l)] (#1) {};}
\def\ClusterD #1[#2] = #3;{
  \node[cluster=#3] (#1t) {}; \node[clabelD=#1] (#1d) {$#2$}; \node[clouter=(#1t)(#1d)] (#1) {};}
\def\ClusterLD #1[#2][#3] = #4;{
  \node[cluster=#4] (#1t) {}; \node[clabelL=#1] (#1l) {$#2$}; 
  \node[clabelD=#1] (#1d) {$#3$}; \node[clouter=(#1t)(#1l)(#1d)] (#1) {};}
\def\ClusterLDName #1[#2][#3][#4] = #5;{
  \node[cluster=#5] (#1t) {}; \node[clabelL=#1] (#1l) {$#2$}; 
  \node[clabelD=#1] (#1d) {$#3$}; 
  \node[scale=\cnamescale,above=\clustersep/3 of #1t,inner sep=0, outer sep=0] (#1n) {$#4$}; 
  \node[clouter=(#1l)(#1d)(#1t)] (#1) {};}
\def\ClustercLDName #1[#2][#3][#4] = #5;{
  \node[clusterc=#5] (#1t) {}; \node[clabelL=#1] (#1l) {$#2$}; 
  \node[clabelD=#1] (#1d) {$#3$}; 
  \node[scale=\cnamescale,above=\clustersep/3 of #1t,inner sep=0, outer sep=0] (#1n) {$#4$}; 
  \node[clouter=(#1l)(#1d)(#1t)] (#1) {};}

\newcommand{\Root}[4][]{
  \ifx\relax#2\relax\node[rr#1=#3] (#4) {};\else\node[rc#1={#2}{#3}] (#4) {};\fi}
\newcommand{\RootT}[5][]{
  \ifx\relax#2\relax\node[rrt#1=#3] (#4) {#5};\else\node[rct#1={#2}{#3}] (#4) {#5};\fi}

\def\frob(#1)(#2){\path[draw,thick,shorten <=-22*\clustersep,shorten >=-22*\clustersep](#1.east)--(#2.west|-#1){};}

\usepackage{pbox}
\def\pb#1{\pbox[c]{\textwidth}{\hfil #1\hfil}}

\long\def\clusterpicture#1\endclusterpicture{\pb{\vbox to \cltopskip{\vfill}\\%
  \begin{tikzpicture}\node[coordinate] (first) {};#1\end{tikzpicture}\\[-11pt]\vbox to \clbottomskip{\vfill}}}   
\long\def\clusterpictureopt#1#2\endclusterpicture{\pb{\vbox to \cltopskip{\vfill}\\%
  \begin{tikzpicture}[#1]\node[coordinate] (first) {};#2\end{tikzpicture}\\[-11pt]\vbox to \clbottomskip{\vfill}}}

\def\pb#1{\pbox[c]{\textwidth}{\hfil #1\hfil}}

\begin{document}

\title[]{Conductor exponents for families of hyperelliptic curves}

\begin{abstract}
We compute the conductor exponents at odd places using the machinery of cluster pictures (developed in \cite{M2D2}) for three infinite families of hyperelliptic curves. These are families of Frey hyperelliptic curves constructed by Kraus (\cite[Section 3]{BCDF}) and Darmon (\cite{DarmonDuke}) in the study of the generalised Fermat equations of signatures $(r,r,p)$ and $(p,p,r)$ respectively. Here, $r$ is a fixed prime number and $p$ is a prime that is allowed to vary. In the context of Darmon's program, Billerey-Chen-Dieulefait-Freitas (\cite{BCDF}) computed all conductor exponents for the signature $(r,r,p)$. We recover their computations at odd places, providing an alternative approach. In a similar setup, Chen-Koutsianas (\cite{chen2022modular}) computed all conductor exponents for the signature $(p,p,5)$. We extend their work to the general case of signature $(p,p,r)$ at odd places.  
Our work can also be used to compute local arithmetic data for the curves in these families à la \cite{M2D2}. 

\end{abstract}

\author{Martin Azon}

\address{INRIA Saclay / École Polytechnique}
\email{martin.azon-y-trell@inria.fr}

\author{Mar Curc\'{o}-Iranzo}
\address{Utrecht Universiteit}
\email{mcurcoiranzo@gmail.com}

\author{Maleeha Khawaja}
\address{University of Warwick}
\email{maleeha.khawaja@warwick.ac.uk}

\author{C\'{e}line Maistret}

\address{University of Bristol}
\email{celine.maistret@bristol.ac.uk}

\author{Diana Mocanu}
\address{Max Planck Institute for Mathematics, Bonn}
\email{diana.mocanu97@outlook.com}

\makeatletter
\let\@wraptoccontribs\wraptoccontribs
\makeatother

\contrib[with an appendix by]{Martin Azon}

\date{}
\thanks{}

\makeatletter
\@namedef{subjclassname@2020}{%
  \textup{2020} Mathematics Subject Classification}
\makeatother


\maketitle

\vspace{-2.1em}

\tableofcontents
%%%%%%%%%%%%%%%%%%%%%%%%%%%%%%%%%%%%%%%%%%%%%%%%%%%%%%%%%%%%%%%%%%%%%%%%%%%%%
%%%%%%%%%%%%%%%%%%%%%%%%%%%%%%%%%%%%%%%%%%%%%%%%%%%%%%%%%%%%%%%%%%%%%%%%%%%%%
\section{\large Introduction}
%%%%%%%%%%%%%%%%%%%%%%%%%%%%%%%%%%%%%%%%%%%%%%%%%%%%%%%%%%%%%%%%%%%%%%%%%%%%%%%%%%%%%%%%%%%%%%%%%%%%%%%%%%%%%%%%%%%%%%%%%%%%%%%%%%%%%%%%%%%%%%%%%%%%%%%%%%
Let $C: y^2=f(x)$ be a hyperelliptic curve defined over a number field $K$. 
The conductor $\conductor(C)$ of $C/K$ is the ideal of the ring of integers $\mathcal{O}_{K}$ of $K$ given by
\begin{equation*}
        \conductor(C) \coloneqq \prod_{\mathfrakq} \mathfrakq^{\condexp{C / K_{\mathfrakq}}},
\end{equation*}
where the product runs over all finite  places $\mathfrakq$ of $K$ and $K_{\mathfrakq}$ denotes the completion of $K$ at $\mathfrakq$. We refer to the integer $\condexp{C / K_{\mathfrakq}}$ as the \textit{conductor exponent} of $C$ at the place $\mathfrakq$ (Definition \ref{def:conductor}). The conductor is one of the key arithmetic invariants associated to a curve and its Jacobian as it is divisible only by the prime ideals at which the Jacobian has bad reduction.  
If $\deg(f)=3$ or $4$, i.e. when $C$ is an elliptic curve, the conductor of $C$ can be computed using Tate's algorithm \cite[pp. 364-368]{SilvermanAdvanced}. 

Tate's algorithm does not have analogues for higher genus curves at present. However, in the case of hyperelliptic curves defined over local fields of odd residue characteristic, the machinery of cluster pictures introduced by Dokchitser, Dokchitser, Maistret and Morgan in \cite{M2D2} provides a formula to compute their conductor exponents in terms of arithmetic data of the roots of their defining polynomials. The method goes roughly as follows. For simplicity, consider a hyperelliptic curve $C/\Q$, described by the Weierstrass equation
$$C: y^2 = f(x).$$ 
Compute the discriminant of $C/\Q$ (see \cite[\S 1]{Lockhart}) and choose an odd prime of bad reduction, say $q$. Compute the roots of $f(x)$ over $\overline{\Q_q}$,  the algebraic closure of the $q$-adic field $\Q_q$, and order them in $q$-adic discs by computing the $q$-adic valuations of their pairwise differences. A \textit{cluster} is a non-empty subset of roots cut out by such a $q$-adic disc (Definition \ref{def:cluster}). The conductor exponent can then be computed from elementary data associated to each cluster and an explicit description of the action of the absolute Galois group of $\Q_q$ on the roots. We refer the reader to Section~\ref{subsec:clustersbackground} for a precise definition of clusters and cluster pictures, and their link to conductor exponents. A more comprehensive overview with examples of such computations and other applications can be found in \cite{hyperusersguide}.

The purpose of this paper is to use this method to compute conductor exponents at odd places in the context of three infinite families of Frey hyperelliptic curves. Our goal is twofold. On one hand, our interest is to provide new results in conductor exponent computations with application towards the modular method.  
On the other hand, our methodology implies direct applications towards the computation of several other local arithmetic data of these families of curves. Indeed, the first step of our approach is to compute the cluster pictures of the curves at each place of bad reduction. In particular, as detailed in \cite{hyperusersguide}, one can obtain local data such as an explicit description of the minimal regular model and the computation of Tamagawa numbers (in the semistable case), an explicit description of the special fibre of the minimal SNC model and the computation of the root numbers (in the tame case). 

As such, this article provides the first instance of using cluster pictures to compute local data in infinite families of hyperelliptic curves as well as using cluster pictures in the context of the modular method.  

\medskip
%%%%%%%%%%%%%%%%%%%%%%%%%%%%%%%%%%%%%%%%%%%%%%%%%%%%%%%%%%%%%%%%%%%%%%%%%%%%%
\subsection{Setup}
%%%%%%%%%%%%%%%%%%%%%%%%%%%%%%%%%%%%%%%%%%%%%%%%%%%%%%%%%%%%%%%%%%%%%%%%%%%%%
Let $p, r \ge 5$ be 
%distinct 
prime numbers. We consider below infinite families constructed by Kraus and Darmon in the study of the generalised Fermat equations of signatures $(r,r,p)$ and $(p,p,r)$ 
\begin{equation}\label{GeneralisedFermatEq}
    x^r + y^r = z^p \quad \text{ and } \quad x^p + y^p = z^r,
\end{equation}
(see \cite[Section 3]{BCDF}, \cite{DarmonDuke} respectively). We note that the roles played by $r$ and $p$ are not symmetric. In the context of the modular method, $r$ is viewed as a fixed prime, whereas $p$ is considered a variable. Our interest for these curves was prompted by the recent work of Billerey-Chen-Dieulefait-Freitas and Chen-Koutsianas, who have used them in this context (see  \cite{BCDF} and \cite{chen2022modular} respectively). Indeed, one key step in the implementation of this method is the computation of the conductor of the curve (see Appendix \ref{sec:compatibility} for the compatibility between their conductor computations and that of the curves). We refer the reader to \cite{DarmonDuke} or \cite{MaleehaSamir} for an overview of the modular method and its relation to conductor computations. 

\begin{remark}
    The case $r = 3$ has been addressed in the literature (see for example \cite{Kraus33p}, \cite{Freitas33p} for the signature $(3, 3, p)$ or \cite{DarmonMerel},\cite{P} for the signature $(p, p, 3)$). Moreover, when $r = 3$, the curves considered below are elliptic curves, and the authors of \cite{M2D2} assume that a hyperelliptic curve has genus at least $2$. However, the conductor of an elliptic curve can be computed using Tate's algorithm. This is why we focus on the case $r \geq 5$. 
    We also note that we do not require $r$ and $p$ to be distinct primes in our computations, although the equation \eqref{GeneralisedFermatEq} has of course been resolved when $r=p$ \cite{Wiles}. 
\end{remark}

\begin{notation}\label{no:global} 
% Throughout the paper, $p, r \geq 5$ will be distinct odd primes. 
We fix $\overline{\Q}$ an algebraic closure of $\Q$.
We consider a primitive $r$-th root of unity $\zr \in \overline{\Q}$ and let $\Q(\zr)$ be the $r$-th cyclotomic field. For any $1 \leq j \leq \frac{r-1}{2}$, we write $\omega_j \coloneqq \zr^{j} + \zr^{-j}$, and set $\omega\coloneqq\omega_1$. We will denote by $\Qzrplus= \Q(\omega)$ the maximal totally real subfield of $\Q(\zr)$, which is generated by 
\begin{equation*}
    h_r(x) \coloneqq \prod_{j = 1}^{\frac{r-1}{2}} (x - \omega_j) \in \Z[x].
\end{equation*}

\end{notation}
\begin{defn}[$C_r$ family]\label{def:cr}
Let $(a, b, c) \in \Z^3$ satisfying $abc\ne 0$, $\gcd(a,b,c)=1$ and $a^r+b^r = c^p$. Following Kraus in \cite[Section 3]{BCDF}, we define the hyperelliptic curve $\Crrp/\Q : y^2 = \frrp(x)$, where
$$
\frrp(x)\coloneqq (ab)^{(r-1)/2} \, x \, h_r\left(\frac{x^{2}}{ab}+2\right)+b^{r}-a^{r} \in \Z[x].
$$
\end{defn}

The polynomial $\frrp$ has degree $r$, so $\Crrp$ has genus $\frac{r-1}{2}$. We note that $\frrp$ and $\Crrp$ only depend on $(a, b)$ and not on the triple $(a,b,c)$. To ease notation however, we write them without any reference to $(a, b)$. We refer the reader to \cite{BCDF} for more details about $C_r$ in the context of the modular method.

\begin{example} 
For small values of $r$, the curve $\Crrp$ is given by
 \begin{align*}
    C_5 : y^2 & = x^5 + 5abx^3 + 5a^2b^2x + b^5 -a^5, \\
    C_7 : y^2 & = x^7 + 7abx^5 + 14a^2b^2x^3 + 7a^3b^3x + b^7 - a^7, \\
    C_{11} : y^2 & =x^{11}+11 a b x^9+44 a^2 b^2 x^7+77 a^3 b^3 x^5+55 a^4 b^4 x^3+11 a^5 b^5 x+b^{11}-a^{11}.
    \end{align*}
\end{example}

\begin{defn}[$\Cminus, \Cplus$ families]\label{def: Cminus & Cplus}
Let $(a, b, c) \in \Z^3$ satisfying $abc\ne 0$, $\gcd(a,b,c)=1$ and $a^p+b^p = c^r$. Following Darmon in \cite{DarmonDuke}, we define the hyperelliptic curves $\Cminus/\Q : y^2 = \gminus(x)$ and $\Cplus/\Q : y^2 = \gplus(x)$, where
  \begin{align*}
        \gminus(x) & \coloneqq (-1)^{\frac{r-1}{2}} c^{r-1} x \, h_r \left( 2 - \frac{x^2}{c^2}\right) -2 \left(a^p - b^p \right) \in \Z[x], \\
        \gplus(x) & \coloneqq \gminus(x) (x+ 2c) \in \Z[x]. 
    \end{align*}
     
\end{defn}

The polynomials $\gminus$ and $\gplus$ have degree $r$ and $r+1$ respectively, so $\Cminus$ and $\Cplus$ both have genus $\frac{r-1}{2}$. As above, to lighten notation, we write $\gminus$, $\gplus$, $\Cminus$ and $\Cplus$ without mentioning $(a, b, c)$. We refer the reader to \cite{chen2022modular} and \cite{DarmonDuke} for more details about $\Cminus$ and $\Cplus$ in the context of the modular method.

\begin{example}
For small values of $r$, the curves $\Cminus$ %and $\Cplus$ 
are given by
\begin{align*}
C_{5}^{-} : y^2 & = x^5-5 c^2 x^3+5 c^4 x-2\left(a^p-b^p\right), \\
C_{7}^{-} : y^2 & = x^7-7 c^2 x^5+14 c^4 x^3-7 c^6 x-2\left(a^p-b^p\right), \\
C_{11}^{-} : y^2 & = x^{11}-11 c^2 x^9+44 c^4 x^7-77 c^6 x^5+55 c^8 x^3-11 c^{10} x-2\left(a^p-b^p\right).\\
%, \\
%\text{and} \qquad \qquad & \\
%C_{5}^{+} : y^2 & = (x+2 c)\left(x^5-5 c^2 x^3+5 c^4 x-2\left(a^p-b^p\right)\right), \\
%C_{7}^{+} : y^2 & = (x+2 c)\left(x^7-7 c^2 x^5+14 c^4 x^3-7 c^6 x-2\left(a^p-b^p\right)\right), \\
%C_{11}^{+} : y^2 & = (x+2 c)\left(x^{11}-11 c^2 x^9+44 c^4 x^7-77 c^6 x^5+55 c^8 x^3-11 c^{10} x-2\left(a^p-b^p\right)\right).
\end{align*}

\end{example}

%\begin{remark}
 %   The construction of the $\Cminus, \Cplus$ families relies on the condition $a^p + b^p = c^r$, while that of the $\Crrp$ family only depends on the coprimality of $a, b$ and $a^r + b^r$.
%\end{remark}

%%%%%%%%%%%%%%%%%%%%%%%%%%%%%%%%%%%%%%%%%%%%%%%%%%%%%%%%%%%%%%%%%%%%%%%%%%%%%
\medskip
\subsection{Main results}
%%%%%%%%%%%%%%%%%%%%%%%%%%%%%%%%%%%%%%%%%%%%%%%%%%%%%%%%%%%%%%%%%%%%%%%%%%%%%

To allow for applications both in the context of the modular method and that of computing local arithmetic data, we consider the curves $\Crrp$, $\Cminus$ and $\Cplus$ as defined over $\Q$ and $\Qplus$. 
The conductor exponent computations are summarized in the two following theorems. 
\begin{theorem}[Theorems \ref{cor:conductorQ}, \ref{thm: conductorCminusQ}, \ref{thm: conductorCplusQ}]  
Consider $\Crrp/\Q$, $\Cminus/\Q$ and $\Cplus/\Q$ as in Definitions \ref{def:cr} and \ref{def: Cminus & Cplus}. 
  \begin{enumerate}
  \item The odd primes of bad reduction for $\Crrp / \Q$ are those dividing $r(a^r + b^r)$. At such a prime $q$, the conductor exponent of $\Crrp / \Q_q$ is 
    $$
    \condexp{\Crrp / \Q_q}=
     \begin{cases} \frac{r-1}{2}  &\text{ if $q\nmid r$ and $q \mid a^r+b^r$}, \\
      r-1 &\text{ if $q=r$}. \end{cases}
     $$
     
     \item The odd primes of bad reduction for $\Cminus / \Q$ are those dividing $rab$.  At such a prime $q$, the conductor exponent of $\Cminus / \Q_q$ is %of $\Cminus /\Q$ at $q$ is 
    $$
    \condexp{\Cminus / \Q_q}=
     \begin{cases}
        \frac{r-1}{2} & \text{ if } q \neq r, \text{ and } r \mid ab, \\ 
        r-1 & \text{ if } q = r,\text{ and } \gminus \text{ is reducible over } \Q_r \, ,\\
        r & \text{ if } q = r, \text{ and } \gminus \text{ is irreducible over } \Q_r \, .    
    \end{cases}
     $$

\item The odd primes of bad reduction for $\Cplus / \Q$ are those dividing $rab$. At such a prime $q$, the conductor exponent of $\Cplus / \Q_q$ is %of $\Cplus /\Q$ at $q$ is 
    $$
    \condexp{\Cplus / \Q_q}=
    \begin{cases}
        \frac{r-1}{2} & \text{if } q \neq r, \text{ and } r \mid ab, \\
        r-1 & \text{if } q=r, \ r \nmid ab  \text{ and } \gminus \text{ is reducible over } \Q_r \, , \\
        r & \text{if } q=r, \ r\nmid ab \text{ and } \gminus \text{ is irreducible over } \Q_r \, , \\
        r-2 & \text{if } q=r, \text{ and } r \mid a,\\
        r-1 & \text{if } q=r, \text{ and } r \mid b.\\
 \end{cases}
     $$
     \end{enumerate}
\end{theorem}

\newpage
\begin{theorem} [Theorems \ref{thm:conductorQplus}, \ref{thm: conductorCminusQpl}, \ref{thm: conductorCplusQpl}]  \label{thm:main}
Consider $\Crrp/\Qplus$, $\Cminus/\Qplus$ and $\Cplus/\Qplus$ as in Definitions \ref{def:cr} and \ref{def: Cminus & Cplus} respectively. Let $\mathfrakq$ be an odd place of $\Qplus$. Write $\Qplus_{\mathfrakq}$ to denote the completion of $\Qplus$ at $\mathfrakq$. 
  \begin{enumerate}
  \item The conductor exponent of $\Crrp / \Qplus_{\mathfrakq}$ is%of $\Crrp /\Qplus$ at $\mathfrakq$ is 
    $$
    \condexp{\Crrp / \Qplus_{\mathfrakq}}=
     \begin{cases} \frac{r-1}{2}  &\text{ if $\mathfrakq\nmid r$ and $\mathfrakq \mid a^r+b^r$}, \\
      r-1 &\text{ if $\mathfrakq\mid r$}. \end{cases}
     $$
     \item The conductor exponent of $\Cminus / \Qplus_{\mathfrakq}$ is%of $\Cminus /\Qplus$ at $\mathfrakq$ is 
    $$
    \condexp{\Cminus /\Qplus_{\mathfrakq}}=
     \begin{cases}
        \frac{r-1}{2} & \text{ if } \qQzrplus \nmid r, \text{ and } \qQzrplus \mid ab,\\ 
        r-1 & \text{ if } \qQzrplus \mid r, \text{ and } \gminus \text{ is reducible over } \Q_r \, ,\\
        \frac{3(r-1)}{2} & \text{ if } \qQzrplus\mid r, \text{ and } \gminus \text{ is irreducible over } \Q_r \, .
    \end{cases}
     $$
      
\item The conductor exponent of $\Cplus / \Qplus_{\mathfrakq}$ is% of $\Cplus /\Qplus_{\mathfrakq}$ at $\mathfrakq$ is 
    $$
    \condexp{\Cplus / \Qplus_{\mathfrakq}}=
     \begin{cases}
        \frac{r-1}{2} & \text{if } \mathfrakq \nmid r, \text{ and } \qQzrplus \mid ab, \\
        r-1 & \text{if } \mathfrakq \mid r, \ r \nmid ab  \text{ and } \gminus \text{ is reducible over $\Q_r$} \, , \\
        \frac{3(r-1)}{2} & \text{if } \mathfrakq\mid r, \ r\nmid ab \text{ and } \gminus \text{ is irreducible over $\Q_r$} \, , \\
        \frac{r-1}{2} & \text{if } \mathfrakq\mid r, \text{ and } r \mid a,\\
        r-1 & \text{if } \mathfrakq\mid r, \text{ and } r \mid b.\\
 \end{cases}
     $$

     \end{enumerate}
\end{theorem}

\begin{remark}
    Note that when $r\mid ab$ the polynomial $g_{r}^{-}$ is always reducible over $\Q_{r}$, see Proposition~\ref{prop: elephant i0}.
\end{remark}

The following corollary shows that it is possible to twist $\Crrp/\Qzrplus_{\mathfrakq}$ (resp. $\Cminus/\Qzrplus_{\mathfrakq}$ and $\Cplus/\Qzrplus_{\mathfrakq}$) to find curves with smaller conductor. This is useful when working within the modular method. 

\begin{cor}[Corollaries~\ref{cor:twistCr}, \ref{cor:twistCmius} and \ref{cor:twistCplus}]\label{cor: Twist Intro}
When $\qQzrplus \mid r$, we have 
\begin{enumerate}
\item If $r \mid (a^r + b^r)$, the twist of $\Crrp$ by a uniformizer of $\Qplus_{\qQzrplus}$ has conductor exponent $\frac{r-1}{2}$.

\item If $r \mid ab$, the twist of $\Cminus$ by a uniformizer of $\Qplus_{\qQzrplus}$ has conductor exponent $\frac{r-1}{2}$.

\item If $r \mid b$, the twist of $\Cplus$ by a uniformizer of $\Qplus_{\qQzrplus}$ has conductor exponent $\frac{r-1}{2}$.
\end{enumerate}
\end{cor}

\begin{remark} In the context of solving
generalised Fermat equations (\ref{GeneralisedFermatEq})
one assumes that the triple $(a,b,c)\in \Z^3$ is a non-trivial primitive solution (\ie \ $abc \neq 0$ and $\gcd(a, b, c) = 1$) and seeks a contradiction via the modular method (see \cite{DarmonDuke} and \cite{MaleehaSamir} for more details). As such, this triple should not exist. In particular, a priori, one cannot perform numerical verification of the above results for the families $\Cminus$ and $\Cplus$. However, one can work around this issue by setting some values for say, $a, c, p, r$, and input the integer $b^p$ or $b^r$ into the defining polynomial of the curves. We did so and ran extensive numerical verification using Magma \cite{Magma} functions for conductor computations and Tim Dokchitser's implementation of cluster pictures. All codes and computations are available in \cite{Git}.
\end{remark}
%%%%%%%%%%%%%%%%%%%%%%%%%%%%%%%%%%%%%%%%%%%%%%%%%%%%%%%%%%%%%%%%%%%%%%%%%%%%%
\subsection{Related work}
%%%%%%%%%%%%%%%%%%%%%%%%%%%%%%%%%%%%%%%%%%%%%%%%%%%%%%%%%%%%%%%%%%%%%%%%%%%%%
Theorem \ref{thm:main} and Appendix \ref{sec:compatibility} show that, at odd places $\mathfrakq$ of $\Qplus$, we recover the computations performed in \cite{BCDF} for the $C_r$ family as well as those in \cite{chen2022modular} for the $C_5^-$ and $C_5^+$ families. Our results also agree with the upper bounds on conductor exponents given by Darmon in \cite{DarmonDuke}. However, we found Proposition $1.16 \, (1)$ in \textit{loc. cit.} inaccurate as currently stated. The curves $\Cminus, \Cplus$ need to be twisted by a uniformizer of the local field as in Corollary~\ref{cor: Twist Intro} to achieve the claimed conductor exponents.

We note that we did not contribute results at places above 2, since the theory of clusters is limited to odd residue characteristics. Also, in the context of the modular method, in addition to computing the conductor exponents at bad places for the variety considered, it is useful to compute the so-called inertial type at these places (see \cite{DembeleFreitasVoight} for precise role and definition). For these issues we refer to \cite{AzonGFE}, where the first author generalizes the results presented below to the case of even residue characteristic, and describes the inertial types of the Galois representations attached to the curves considered here. 

%We expect to be able to extend the present work to include the description of the inertial types in each case considered above. This is work in progress.

%%%%%%%%%%%%%%%%%%%%%%%%%%%%%%%%%%%%%%%%%%%%%%%%%%%%%%%%%%%%%%%%%%%%%%%%%%%%%
\subsection{Overview}
%%%%%%%%%%%%%%%%%%%%%%%%%%%%%%%%%%%%%%%%%%%%%%%%%%%%%%%%%%%%%%%%%%%%%%%%%%%%%
In Section \ref{sec:background} we review background material on curves, Jacobians and conductors (\ref{subsec:CJC}), cluster pictures (\ref{subsec:clustersbackground}) and how to compute conductor exponents at odd places from them (Theorems \ref{thm:condss} and \ref{thm: allcond}). Section \ref{part:rrp} focuses on the $C_r$ family associated to the signature $(r, r, p)$. We first exhibit algebraic expressions for the roots of the defining polynomials (Proposition~\ref{prop:rrproots}), construct the cluster pictures at places of bad reduction for $C_r/\Q$ (Theorem~\ref{thm: Clusters rrp Q}) and $C_r/\Qplus$ (Corollary~\ref{cor: Clusters Crrp Qplus}), and proceed to compute the conductor exponents in each case (Theorems~\ref{cor:conductorQ} and \ref{thm:conductorQplus}). Section \ref{part:ppr} is concerned with the pair of families $\Cminus$ and $\Cplus$ associated to the signature $(p, p, r)$. We start by presenting an algebraic expression of the roots of the defining polynomials (Proposition~\ref{prop: Rootsppr}) and construct cluster pictures at places of bad reduction for $\Cminus/\Q$ (Theorem~\ref{thm: ClusterCminusQ}), $\Cplus/\Q$ (Theorem~\ref{thm: ClusterCplusQ}), $\Cminus/\Qplus$ (Corollary~\ref{cor: ClusterCminusQpl}) and $\Cplus/\Qplus$ (Corollary~\ref{cor: ClusterCplusQpl}). Before proceeding to compute the conductor exponents in Section \ref{sec: Conductors ppr}, we investigate the ramification index of the splitting field of the defining polynomials and compute some necessary discriminants. This is done in Section \ref{sec: Ramification indices}. Appendix \ref{sec:compatibility} presents the compatibility between our computations and those of Billerey-Chen-Dieulefait-Freitas \cite{BCDF} and Chen-Koutsianas \cite{chen2022modular} in the context of the modular method.

%We note that throughout the paper, when referring to statements, we may mean the statement itself or an argument developed in its proof. \\

\subsection{Acknowledgements.} We would like to warmly thank the referees for their careful reading of this manuscript and their meaningful suggestions and comments. We would also like to thank Nicolas Billerey for suggesting this project, and Julian Lyczak and Ross Paterson for organizing the Workshop on Arithmetic and Algebra of Rational Points, at which the project started. We are grateful to Nicolas Billerey, Matthew Bisatt, Imin Chen, Angelos Koutsianas and Samir Siksek for helpful discussions. 
The third and fourth authors were supported by the Royal Society Dorothy Hodgkin Fellowship of the fourth author. 
%%%%%%%%%%%%%%%%%%%%%%%%%%%%%%%%%%%%%%%%%%%%%%%%%%%%%%%%%%%%%%%%%%%%%%%%%%%%%
%%%%%%%%%%%%%%%%%%%%%%%%%%%%%%%%%%%%%%%%%%%%%%%%%%%%%%%%%%%%%%%%%%%%%%%%%%%%%
\section{\large Notation and Background} \label{sec:background}
%%%%%%%%%%%%%%%%%%%%%%%%%%%%%%%%%%%%%%%%%%%%%%%%%%%%%%%%%%%%%%%%%%%%%%%%%%%%%
%%%%%%%%%%%%%%%%%%%%%%%%%%%%%%%%%%%%%%%%%%%%%%%%%%%%%%%%%%%%%%%%%%%%%%%%%%%%%

%%%%%%%%%%%%%%%%%%%%%%%%%%%%%%%%%%%%%%%%%%%%%%%%%%%%%%%%%%%%%%%%%%%%%%%%%%%%%
\subsection{Notation}\label{sec:notation}
%%%%%%%%%%%%%%%%%%%%%%%%%%%%%%%%%%%%%%%%%%%%%%%%%%%%%%%%%%%%%%%%%%%%%%%%%%%%%
As in Notation \ref{no:global}, we let $p, r \geq 5$ denote prime numbers, $\zr \in \overline{\Q}$ a primitive $r$-th root of unity and $\Q(\zr)$ the $r$-th cyclotomic field. For any $1 \leq j \leq \frac{r-1}{2}$, we write $\omega_j \coloneqq \zr^{j} + \zr^{-j}$, and set $\omega\coloneqq\omega_1$. We denote by $\Qzrplus= \Q(\omega)$ the maximal totally real subfield of $\Q(\zr)$.

For a fixed rational prime $q$, we let $\Q_q$ be the field of $q$-adic numbers and $v_q$ be the valuation on $\overline{\Q_q}$ normalized with respect to $q$, \ie , satisfying $v_q(q)=1$. For a prime ideal $\mathfrakq$ in the ring of integers $\cO_{\Qplus}$ of $\Qplus$, we denote by $\Qplus_{\mathfrakq}$ the completion of $\Qplus$ at $\mathfrakq$. We let $v_{\mathfrakq}$ be the valuation on $\overline{\Qplus_{\mathfrakq}}$ normalized with respect to $\mathfrakq$ \ie , satisfying $v_{\mathfrakq}(\pi)=1$ for any uniformizer $\pi$ of $\Qplus_{\mathfrakq}$. We point out that any such valuation is Galois-invariant (see \cite[Chapter II \S 8]{Neukirch}).

If we suppose that $\mathfrakq$ is a prime above $q$ in $\cO_{\Qplus}$ standard local field theory implies the isomorphism
\begin{equation}\label{eq: Isom compl+ext}
    \Qzrplus_{\qQzrplus} \simeq \Q_q(\omega). 
\end{equation}

Recall that the extension $\Q(\zr)/\Q$
is unramified at all primes $q\ne r$, and totally ramified at $r$. In particular, the same holds for $\Qplus/\Q$.
Denote by $\rQzrplus$ the unique prime ideal above $r$ in $\Qplus$. For $a \in \Qplus$, we have the following link between valuations and ramification indices.  
\begin{equation}\label{eq:valsrelation}
v_{\mathfrak{q}}(a)= 
\begin{cases}
v_{q}(a) & \text{if } \mathfrak{q} \neq \rQzrplus, \\
\frac{r-1}{2}v_{r}(a) & \text{ if } \mathfrak{q}=\rQzrplus.
\end{cases}
\end{equation}

Given any variety $V$ defined over a field $K$ and a field extension $L / K$, we denote by $V / L\coloneq V\times_K L$ the base change of $V$ to $L$.

Given a hyperelliptic curve $C$ over a number field $\mathcal{F}$ and an odd place $\mathfrak p$ of bad reduction, we refer to the cluster picture of $C$ at $\mathfrak p$ for the cluster picture of $C$ over the local field $\mathcal{F}_{\mathfrak p}$ as introduced in Definition \ref{def:cluster}. 
 We will use the notation $\condexp{C/\mathcal{F}_{\mathfrak p}}$ for the conductor exponent of $C$ at $\mathfrak p$ (see Definition \ref{def:conductor}).
%%%%%%%%%%%%%%%%%%%%%%%%%%%%%%%%%%%%%%%%%%%%%%%%%%%%%%%%%%%%%%%%%%%%%%%%%%%%%
\subsection{Curves, Jacobians and conductors}\label{subsec:CJC}
%%%%%%%%%%%%%%%%%%%%%%%%%%%%%%%%%%%%%%%%%%%%%%%%%%%%%%%%%%%%%%%%%%%%%%%%%%%%%
Consider a smooth projective curve $C$ of genus $g$ defined over a local field $K$. The Jacobian variety $\Jac(C)$ of $C$ is an abelian variety of dimension $g$ whose elements are linear equivalence classes of degree $0$ line bundles on $C$. We refer the reader to \cite[Section A.8]{Milne86} for details on the construction of the Jacobian. Let $\ell$ be a prime distinct from the residue characteristic of $K$ and consider the $\ell$-adic Tate module of $\Jac(C)$
\begin{equation*}
    \Tl \left( \Jac(C)\right) \coloneqq \varprojlim_{n \geq 1} \Jac(C) \left[\ell^n \right].
\end{equation*}
As explained in \cite{SerreTate}, $\Tl \left( \Jac(C)\right)$ is a free $\Zl$-module of rank $2g$ on which $G_K$ acts continuously. Then $\Vl \coloneqq \Tl \left( \Jac(C)\right) \otimes_{\Zl} \Ql$ gives a continuous representation of the absolute Galois group $G_K$ of $K$
\begin{equation*}
    \rhol : G_{K} \rightarrow \Aut (\Vl) \simeq \GL_{2g}(\Ql).
\end{equation*}

For any finite Galois group $G$ and any real number $u \geq -1$, we denote by $G^u$ the $u$-th ramification group of $G$ in upper numbering (see \cite[Chapter II, §10]{Neukirch} or \cite[§3]{Ulmer}). The absolute $u$-th ramification group of $K$ in upper numbering is defined by
\begin{equation}\label{eq: Ramif groups}
G_K^u \coloneqq \varprojlim_{\overline{K}/L/K}\Gal(L/K)^u,
\end{equation}
where the inverse limit is taken over all intermediate finite Galois extensions ${\overline{K}/L/K}$. The absolute inertia group $I_K$ of $K$ matches $G^0$.
For any group $\mathcal{G}$ acting on $\Vl$, we denote by $\Vl^{\mathcal{G}}$ the subspace of $V_\ell$ fixed by $\mathcal{G}$. 
We now introduce the conductor exponent of the curve $C$ defined over the local field $K$. 
\begin{defn}\label{def:conductor}
The tame and wild parts of the conductor exponent of $C / K$ are defined respectively by
\begin{equation*}
    \condtame{C/K} \coloneqq  \int_{-1}^{0} \codim_{\Q_\ell} \Vl^{G_K^u}\, du \quad \text{ and } \quad \condwild{C/K}\coloneqq\int_{0}^{\infty} \codim_{\Q_\ell} \Vl^{G_K^u}\, du.
\end{equation*}
We define the conductor exponent of $C / K$ by $\condexp{C/K} \coloneqq \condtame{C/K} + \condwild{C/K}$.
\begin{comment}
\[
\condexp{C/K}\coloneqq\int_{-1}^{\infty} \codim_{\Q_\ell} V^{G_K^u}\, du.
\]
Moreover, we define the tame conductor exponent of $C/K$ by 
\[
 \condtame{C/K} \coloneqq  \int_{-1}^{0} \codim_{\Q_\ell} V^{G_K^u}\, du = \codim_{\Q_\ell} V^{(I_K)},
\]
where $I_K$ is the inertia subgroup of $G_K$ and the wild conductor exponent (or the Swan conductor) of $C/K$ as
\[
\condwild{C/K}\coloneqq\int_{0}^{\infty} \codim_{\Q_\ell} V^{\rho_{\ell}(G_K^u)}\, du.
\]
\end{comment}
\end{defn}

We can describe the behaviour of the wild conductor with respect to tame base change.

\begin{lemma}\label{lem:wildthoughtame}
Let $K$ be a local field of odd residue characteristic. Suppose that $F/K$ is a finite tamely ramified extension with ramification index $e_{F/K}$.
Then 
\[
\condwild{C/F}=e_{F/K} \, \condwild{C/K}.
\]
\end{lemma}
\begin{proof}
Since $F/K$ is a tame extension, then one can deduce from \cite[IV. Proposition 15]{SerreLF} that $G_{F}^{u}= G_K^{u/e}$, for all $u>0$, where $e\coloneqq e_{F/K}$. Hence, we compute 
\begin{align*}
\condwild{C/F}= \int_{0}^{\infty} \codim_{\Q_\ell} V^{G_{F}^u}\, du & = \int_{0}^{\infty} \codim_{\Q_\ell} V^{G_{K}^{u/e}}\, du\\
&= e\int_{0}^{\infty} \codim_{\Q_\ell} V^{G_{K}^{t}}\, dt = e\,\condwild{C/K}.
\end{align*}
Here, the third equality follows from a change of variables given by $t\coloneqq u/e$.
\end{proof}

\begin{remark}
    The $\ell$-adic Tate module of an abelian variety $A / K$ is isomorphic to the dual of the first $\ell$-adic cohomology group $\HH^{1}(A / \Kbar, \Zl)$ (see \cite[Theorem 15.1]{Milne84}). In the specific case of Jacobians, there is an isomorphism $\HH^{1}(C / \Kbar, \Zl) \simeq \HH^{1}(\Jac(C) / \Kbar, \Zl)$ (\cite[Corollary 9.6]{Milne86}), so we obtain $\HH^{1}(C / \Kbar, \Ql) \simeq \Vl(\Jac(C))^{\vee}$. Thus, one can identify the representation $\rhol$ and the representation arising from the action of $G_K$ on $\HH^{1}(C / \Kbar, \Ql).$ Therefore, it makes sense to talk about the conductor exponent of the curve itself without referring to its Jacobian.
\end{remark}
%%%%%%%%%%%%%%%%%%%%%%%%%%%%%%%%%%%%%%%%%%%%%%%%%%%%%%%%%%%%%%%%%%%%%%%%%%%%%
\subsection{Cluster pictures}\label{subsec:clustersbackground}
%%%%%%%%%%%%%%%%%%%%%%%%%%%%%%%%%%%%%%%%%%%%%%%%%%%%%%%%%%%%%%%%%%%%%%%%%%%%%

Let $K/\Q_q$ be a finite extension with $q$ an odd prime number. Let $C/K$ be a hyperelliptic curve of genus $g$ given by the affine equation
$$
y^2 = f(x) = c \prod_{\gamma \in \mathcal{R}}(x-\gamma),
$$
where $c \in K^{\times}$, $f \in K[x]$ is separable, and $\mathcal{R}$ denotes the set of roots of $f$ in $\Kbar$. We let $v$ be a valuation on $\overline{K}$, normalized with respect to a uniformizer of $K$.

\begin{defn}[Clusters and cluster pictures]
\label{def:cluster}
A \textit{cluster} is a non-empty subset $\s \subseteq \mathcal{R}$ of the form $\s = D \cap \mathcal{R}$ for some disc $D=\{x\!\in\! \Kbar \mid v(x-z)\!\geq\! d\}$, some $z\in \Kbar$ and some $d\in \Q$.

A cluster $\s$ with $|\s|>1$ is called a \textit{proper} cluster. For a proper cluster $\s$, its \textit{depth} $d_\s$ is the maximal $d$ for which $\s$ is cut out by such a disc, \ie ,  $d_\s\! =\! \min_{\gamma,\gamma' \in \mathfrak{s}} v(\gamma\!-\!\gamma')$. If moreover $\s\neq \mathcal{R}$, we define the \textit{parent} cluster $P(\s)$ of $\s$ to be the smallest cluster with $\s\subsetneq P(\s)$. Then the \textit{relative depth} of $\s$ is $\delta_\s\! =\! d_\s\! -\!d_{P(\s)}$.

We refer to this data as the \textit{cluster picture} of $C$.
\end{defn}

\begin{remark}\label{CLrem1}
The absolute Galois group $G_K$ of $K$ acts on clusters via its action on the roots. This action preserves the depths and containments of clusters. 
\end{remark}

We now recall some terminology and notations concerning clusters.

\begin{defn}\label{defclusterinv} Consider the cluster picture of $C/K$.
\begin{itemize}
\item If $\s'\subsetneq \s$ is a maximal subcluster, we refer to $\s'$ as a \textit{child} of $\s$.
For two clusters $\s_1$, $\s_2$ write $\s_1\wedge \s_2$ for the smallest cluster containing both of them.

\item A cluster $\s$ such that $|\s|$ is odd (resp. even) is called an \textit{odd} (resp. \textit{even}) cluster. We denote by $\tilde{\s}$ the set of odd children of $\s$. 
\item A cluster consisting of precisely two roots is called a \textit{twin}. 
\item  An even cluster whose children are all even is called \textit{\"ubereven}. 
\item A non-\"ubereven cluster with a child of size $2g$ is called a \textit{cotwin}.
\item A cluster $\s$ is called principal if one of the following holds
\[
\begin{cases}
|\s|\neq 2g+2 \text{ and } \s \text{ is proper, not a twin or a cotwin;}\\
|\s|=2g+2 \text{ and } \s \text{ has more than }3 \text{ children.}
\end{cases}
\]

\item \begin{comment} We let $L\coloneqq K(\Rroots)$ be the splitting field of $f(x)$ and $I_{L}$ the inertia subgroup of $\Gal(L/K)$.
\end{comment}
For a cluster $\s$, we denote by $I_{\s}$ the stabiliser of $\s$ under $I_K$. 
%and $I_{L, \s}$  the stabiliser of $\s$ under $I_{L}$. Note that $I_{L,\s} = I_\s \cap \Gal(L/K)$.

\item Fix a cluster $\s$. If there exists a unique fixed child $\s'$ under the action of $I_{\s}$, then we call $\s'$ an orphan of $\s$. 

\item For a proper cluster $\s$, we define the following quantities 
\begin{equation*}
    \tilde\lambda_{\s}=\frac{1}{2}\left( v(c)+|\tilde{\s}|d_{\s}+ \sum_{\gamma \not\in \s}d_{ \{\gamma\} \wedge \, \s } \right), \qquad 
    \nu_{\s}= v(c)+|{\s}|d_{\s}+ \sum_{\gamma \not\in \s}d_{ \{\gamma\} \wedge \, \s}.
\end{equation*}
\end{itemize}

\end{defn}

In a given cluster picture, each root $r \in \mathcal{R}$ is represented by \smash{\raise4pt\hbox{\clusterpicture\Root[D]{1}{first}{r1};\endclusterpicture}} and each proper cluster is represented by an oval. 
The subscript on the largest cluster $\mathcal{R}$ is the depth of $\cR$, whilst the subscripts on the other clusters are their respective relative depths. For example, consider the following cluster picture.
$$
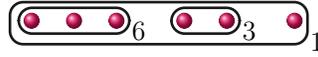
\scalebox{1.25}{\clusterpicture            
  \Root[D] {1} {first} {r1};
  \Root[D] {3} {r1} {r2};
  \Root[D] {3} {r2} {r3};
  \Root[D] {6} {r3} {r4};
  \Root[D] {3} {r4} {r5};
  \Root[D] {6} {r5} {r6};
  \ClusterLD c1[][{6}] = (r1)(r2)(r3);
  \ClusterLD c2[][{3}] = (r4)(r5);
  \ClusterLD c3[][{1}] = (c1)(c2)(r6);
\endclusterpicture}
$$
In this example, the largest cluster $\cR$ has depth $1$, the cluster containing $3$ roots has depth $7$, and the twin has depth $4$. For further details, we refer the reader to \cite{hyperusersguide} wherein the authors compute numerous examples of cluster pictures as well as explain how to extract various arithmetic data for the curve in question.

Cluster pictures provide an explicit criterion to check if the hyperelliptic curve $C/K$ is semistable.

\begin{defn}\label{def:ss}
We say that $C / K$ satisfies \textit{the semistability criterion} if the following conditions hold:
\begin{enumerate}
\item The field extension $K(\mathcal R)/K$ given by adjoining the roots of $f(x)$ has ramification degree at most $2$.
\item Every proper cluster is invariant under the action of the absolute inertia group $I_K$.
\item Every principal cluster $\s$ has $d_{\s} \in \Z$ and $\nu_{\s} \in 2\Z$.
\end{enumerate}
\end{defn}

\begin{theorem}[Theorem 10.3 in \cite{M2D2}]\label{thm:ss}
The hyperelliptic curve $C/K$, equivalently $\Jac(C)$, is semistable if and only if it satisfies the semistability criterion.
\end{theorem}

%Following \cite{M2D2}, we refer to Theorem \ref{thm:ss} as the \textit{semistability criterion}. 
The conductor exponent of a semistable hyperelliptic curve can easily be computed from its cluster picture via the following result.

\begin{theorem}[Corollary 9.4 in \cite{M2D2}]\label{thm:condss}
Assume that $C/K$ is semistable. Its conductor exponent $\condexp{C / K}$ satisfies
\begin{equation*}
\condexp{C / K} =
\begin{cases}
\; |A|-1 & \text{ if }\mathcal{R}\text{ is \"ubereven},\\
\; |A| & \text{ else},
\end{cases}
\end{equation*}
where $A \coloneqq \{ \text{even clusters }\mathfrak{s}\neq \mathcal{R} \mid \s \text{ not \"ubereven} \}$.
\end{theorem}

When $C/K$ is not semistable, one can still compute its conductor exponent from its cluster picture but more notation is needed. Given a proper cluster $\s$ and $a\in\Q$, we define 
$$\xi_{\s}(a) = \max \{ -v_2([I_K:I_{\s}]a) ,\; 0\}.$$ 

To the cluster picture of a hyperelliptic curve $C/K$, we associate the sets
\begin{eqnarray*}
	U &=& \{ \text{odd clusters } \s \neq \mathcal{R} \mid \xi_{P(\s)}(\tilde\lambda_{P(\s)}) \leq \xi_{P(\s)}(d_{P(\s)}) \}; \\
	V &=& \{ \text{proper non-\"{u}bereven clusters } \s \mid \xi_{\s}(\tilde\lambda_{\s})=0 \}.
\end{eqnarray*}

%We now compute the conductor exponent in terms of the data described above.

\begin{theorem}[Theorem 12.3 in \cite{hyperusersguide}]
\label{thm: allcond}
%Let $C/K$ be a hyperelliptic curve. We have
The tame and wild part of the conductor exponent of $C/K$ are given by
\begin{enumerate}
\item $\condtame{C/K} = 2g - \#(U/I_K) + \#(V/I_K) + \begin{cases} 1  &\text{ if $|\mathcal{R}|$ and $v(c)$ are even,} \\ 0 &\text{ else;} \end{cases}$
\medskip
\item $\condwild{C/K} = \sum\limits_{\gamma \in \, \mathcal{R}/G_K} \left( v(\Delta_{K(\gamma)/K}) - [K(\gamma):K] + f_{K(\gamma)/K} \right),$\\ where $\Delta_{K(\gamma)/K}$ and $f_{K(\gamma)/K}$ are the discriminant and residue degree of $K(\gamma)/K$ respectively.
\end{enumerate}
\end{theorem} 

\begin{remark}\label{rem:wildplus}
Suppose we are given two curves $C/K$ and $C'/K$ with sets of roots $\cR$ and $\cR'=\cR \cup \{\gamma\}$ respectively, where $\gamma \in K$. Then an immediate consequence of Theorem~\ref{thm: allcond} is that $\condwild{C/K}=\condwild{C'/K}$.
\end{remark}
Next, we prove that roots inside tame extensions do not contribute to the wild part of the conductor. 
\begin{proposition}\label{cor:nwild0}
If $K(\mathcal{R})/K$ is a tame extension then $\condwild{C/K}=0$.
\end{proposition}
\begin{proof}
Suppose that $L/K$ is a tamely ramified extension of $K$ with normalized valuation $v_L$ extending $v_K$.
We claim that
\begin{equation}\label{lem:nwild0}
v_K(\Delta_{L/K}) - [L:K] + f_{L/K}=0
\end{equation}
where $\Delta_{L/K}$ is the discriminant of $L/K$, and $f_{L/K}$ is the residue degree of $L/K$. Since $L/K$ is tame \cite[III. Proposition 13]{SerreLF} yields $v_L(\mathcal D_{L/K})=e_{L/K}-1$, where $\mathcal D_{L/K}$ is the different and $e_{L/K}$ is the ramification index of $L/K$. Hence $v_K(\Delta_{L/K})=v_K(\Norm(\mathcal D_{L/K}))= f_{L/K}(e_{L/K}-1)$.
This together with the fact that $[L:K]=e_{L/K}f_{L/K}$ can be rearranged to give \eqref{lem:nwild0}. The statement then follows from the expression of $\condwild{C/K}$ given in Theorem~\ref{thm: allcond}.
\end{proof}
When computing the tame part of the conductor, %$\condtame{C/K}$
the following result will be useful. 
\begin{theorem}[Theorem 1.3(iv) in \cite{Bisatt}] \label{thm:tameindex}
Assume that $K(\cR)/K$ is a tame extension. Let $\s$ be a cluster. Then
\[[I_{K} : I_{\s}] = \lcm_{\s \subsetneq \s'} \mathrm{denom} \ d^{\ast}_{\s'},\]
where for a cluster $\s' \supsetneq \s$,
\[ d^{\ast}_{\s'} = \begin{cases}
        1 & \text{if the child of } \s' \text{ containing } \s \text{ is an orphan,} \\
        d_{\s'} &  \text{else.} \end{cases}\]
\end{theorem}
~
%%%%%%%%%%%%%%%%%%%%%%%%%%%%%%%%%%%%%%%%%%%%%%%%%%%%%%%%%%%%%%%%%%%%%%%%%%%%%
\subsection{On the factorization of \texorpdfstring{$X^r + Y^r$}{}} 
%%%%%%%%%%%%%%%%%%%%%%%%%%%%%%%%%%%%%%%%%%%%%%%%%%%%%%%%%%%%%%%%%%%%%%%%%%%%%
Keeping notation as in Section \ref{subsec:clustersbackground}, we discuss some results about the factorization of the two-variable polynomial $X^r + Y^r$, which we view as an element of $\Z[X, Y]$. Recall that $\zr \in \overline{\Q_r}$ is a primitive $r$-th root of unity and  $\omega_j = \zr^j + \zr^{-j}$.
%.\in K[X,Y]$. 
%This will inform the link between the defining polynomial of $\Crrp$ (resp. $\Cminus$ and $\Cplus$) and the Diophantine equations of signature $(r, r, p)$ (resp. $(p, p, r)$). 

\begin{defn}\label{def: phir}
    We define the two-variable polynomial $\phi_r(X, Y) \in \Z[X, Y]$ by the formula
    \begin{equation*}
        \phi_r(X, Y) \coloneqq \frac{X^r + Y^r}{X+Y} = \sum_{j = 0}^{r-1} (-1)^j X^{r-1-j} Y^j.
    \end{equation*}
    One can factorize $\phi_r$ over the cyclotomic field $\Q(\zr)$ as follows
    \begin{equation}\label{eq: Factorisation phir}
        \phi_r(X, Y) = \prod_{j = 1}^{r-1} \left(X+ \zr^j Y \right) = \prod_{j = 1}^{\frac{r-1}{2}} \left( X^2 + \omega_j X Y + Y^2 \right).
    \end{equation}
\end{defn}

The following lemma gives an alternative way to write the polynomial $X^r + Y^r$ using $h_r$, the defining polynomial of $\Qplus = \Q(\omega)$. This will help us compute the roots of the defining polynomials $\frrp$ and $\gminus, \gplus$. 

\begin{lemma} \label{lem:phiandh}
For any odd prime $r$, we have the following equality of polynomials
\begin{equation*}
(-X Y)^{\frac{r-1}{2}} (X + Y) h_{r}\left(2 - \frac{(X + Y)^2}{X Y}\right) = X^r + Y^r. 
\end{equation*}
\end{lemma}

\begin{proof}
    Using \eqref{eq: Factorisation phir}, we compute
    \begin{multline*}
        (-X Y)^{\frac{r-1}{2}}  h_{r}\left(2 - \frac{(X + Y)^2}{X Y} \right) = (-X Y)^{\frac{r-1}{2}}  h_{r}\left(2 - \frac{X^2 + Y^2}{X Y} - 2 \right) \\
         = (-X Y)^{\frac{r-1}{2}} \prod_{j = 1}^{\frac{r-1}{2}} \left(- \frac{X^2 + Y^2}{X Y} - \omega_j \right) = \prod_{j = 1}^{\frac{r-1}{2}} \left( X^2 + \omega_j X Y + Y^2 \right) =  \phi_{r}(X, Y),
    \end{multline*}
    and the equality $X^r + Y^r = (X + Y) \, \phi_r(X, Y)$ allows us to conclude.
\end{proof}

As explained in Definition~\ref{def: phir}, we have $X^r + Y^r = (X+ Y) \phi_r(X, Y) = \prod_{j = 0}^{r-1} (X + \zr^j Y)$ (see \eqref{eq: Factorisation phir}). We now describe the valuation of each of the factors in the product in terms of the valuation of the left-hand side. Recall that we denote by $K$ a local field of odd residue characteristic and by $v$ a valuation on $\Kbar$, normalized with respect to a uniformizer of $K$.

\begin{lemma} \label{lem:generalproduct}
Let $\alpha, \beta \in \mathcal{O}_K$ be such that $v(\alpha^r + \beta^r)>0$ and $v(\alpha) v(\beta) = 0$.
\begin{enumerate}
    \item If $q \neq r$ there exists some $0 \leq j_{0} \leq r-1$ such that $v(\alpha +\zr^{j_{0}}\beta) = v(\alpha^r + \beta^r ) >0$ and $v(\alpha +\zr^{k}\beta)=0$ for all $0 \leq k \leq r-1$, $k \neq j_{0}$. 

    \item If $q = r$ then $v(\alpha + \zr^j \beta) >0$ for any $0 \leq j \leq r-1$. Moreover, if $v(\alpha^r+\beta^r) > r \, v(1-\zr)$ then there exists some $0 \leq j_0 \leq r-1$ such that $v\left( \alpha + \zr^{k} \beta \right) = v\left( 1 - \zr \right)$ for any $0 \leq k \leq r-1$, $\ k \neq j_0$ and $v(\alpha + \zr^{j_0} \beta) > v(1 - \zr)$.
\end{enumerate}
Furthermore, if $\zr \notin K$, then we have $j_0 = 0$.
%in cases (1) and (2) respectively. 
\end{lemma}

\begin{proof} 
The assumptions $v(\alpha^r + \beta^r) > 0$ and $v(\alpha)v(\beta)=0$ combined imply $v(\beta)=v(\alpha) = 0$. 
Suppose $q \neq r$. Since $\alpha^r + \beta^r = \prod_{j=0}^{r-1} \left( \alpha + \zr^j \beta \right)$ there is some $0 \leq j_0 \leq r-1$ such that $v(\alpha + \zr^{j_0} \beta)>0$. For any $0 \leq k \leq r-1$, we have 
\begin{equation} \label{eq:decomp}
\alpha + \zr^k \beta = \left(\alpha + \zr^{j_0} \beta \right) + \zr^{k} \left( 1 - \zr^{j_0 -k} \right) \beta,
\end{equation}
which implies that, when $k \neq j_0$, we have
%Hence, when $q \neq r$, for all $0 \leq k \leq r-1$, $k \neq j_{0}$ we have
\begin{equation} \label{eq:valcoprime}
v(\alpha + \zr^k \beta) = \text{min}\{v\left(\alpha + \zr^{j_0} \beta \right), v(\zr^{k} \left( 1 - \zr^{j_0-k} \right) \beta)\} =0,
\end{equation}
because $v(\zr^{k} \left( 1 - \zr^{j_0-k} \right) \beta) =0$. Therefore we deduce $v(\alpha^r + \beta^r) = v(\alpha + \zr^{j_0} \beta)$.

Suppose now $q=r$. The algebraic expression~\eqref{eq:decomp} still holds, so we have $v(\alpha + \zr^k \beta) > 0$ for any $0 \leq k \leq r-1$. If we assume $v(\alpha^r + \beta^r) > r \, v(1 - \zr)$, then there exists some $0 \leq j_0 \leq r-1$ such that $v(\alpha + \zr^{j_0} \beta) > v(1 - \zr)$. Indeed, if we had $v(\alpha + \zr^{k} \beta) \leq v(1 - \zr)$ for all $0 \leq k \leq r-1$, then $v(\alpha^r + \beta^r) = \sum_{k=0}^{r-1} v(\alpha + \zr^k \beta) \leq \sum_{k = 0}^{r-1} v(1 - \zr) = r \, v(1-\zr)$, contradicting the assumption.
%Let $\nu$ be the valuation on $K(\zr)$ extending $v$. The assumption $v(\alpha^r+\beta^r) > r v(1-\zr)$ translates into $\nu(\alpha^r+\beta^r) > r \nu(1-\zr)$.  By the pigeon-hole principle and \eqref{eq: Factorisation phir}, there exists some $j_0$ such that $\nu(\alpha+\zr^{j_0}\beta) > \nu(1-\zr)$, which means that $v(\alpha+\zr^{j_0}\beta) > v(1-\zr)$.
We deduce the existence of such $j_0$. Then for $k \neq j_0$,
we have $v(\zr^{j_0} (1 - \zr^{k - j_0}) \beta ) = v(1-\zr^{k - j_0})=v(1-\zr)$, so using~\eqref{eq:decomp} we get
\begin{align*}
v(\alpha + \zr^{k} \beta) & = v( \zr^{j_0}(1-\zr^{k - j_0}) \beta +  (\alpha + \zr^{j_{0}} \beta)) \\
& = \text{min}\{v(1-\zr), v(\alpha + \zr^{j_{0}} \beta)\} = v(1-\zr).
\end{align*}
To conclude, assume $\zr \notin K$, and let us prove that $j_0 = 0$. We focus on the case $q \neq r$, the case $q = r$ is treated similarly. Assume by contradiction that $j_0 \neq 0$. There exists a map $\sigma \in \text{Gal}(K(\zr)/K)$ not fixing $\zr^{j_0}$, \ie \ $\sigma(\zr^{j_{0}})=\zr^{k_{0}}$ for some $1 \leq k_{0} \leq r-1$, $k_{0} \neq j_{0}$. Hence, using the Galois-invariance property of the valuation $v$, we get 
$$ v( \alpha + \zr^{j_{0}} \beta)=v( \sigma(\alpha + \zr^{j_{0}} \beta))= v(\alpha + \zr^{k_{0}} \beta)>0,$$
which contradicts \eqref{eq:valcoprime}, and we conclude that $j_{0}=0$. 
\end{proof}

%\begin{remark}
    %Note that the valuation $v$ in Lemma~\ref{lem:generalproduct} is not required to be normalized with respect to $K$.
%\end{remark}

%%%%%%%%%%%%%%%%%%%%%%%%%%%%%%%%%%%%%%%%%%%%%%%%%%%%%%%%%%%%%%%%%%%%%%%%%%%%%
%%%%%%%%%%%%%%%%%%%%%%%%%%%%%%%%%%%%%%%%%%%%%%%%%%%%%%%%%%%%%%%%%%%%%%%%%%%%%
\section{\large \texorpdfstring{$\Crrp$}{Crrp} family associated to the signature \texorpdfstring{$(r, r, p)$}{rrp}}\label{part:rrp}
%%%%%%%%%%%%%%%%%%%%%%%%%%%%%%%%%%%%%%%%%%%%%%%%%%%%%%%%%%%%%%%%%%%%%%%%%%%%%
%%%%%%%%%%%%%%%%%%%%%%%%%%%%%%%%%%%%%%%%%%%%%%%%%%%%%%%%%%%%%%%%%%%%%%%%%%%%%
In this section we consider the family of Frey hyperelliptic curves associated to the generalised Fermat equation of signature $(r, r, p)$. Let $(a, b, c) \in \Z^3$ satisfying $abc\ne 0$, $\gcd(a,b,c)=1$ and $a^r+b^r = c^p$. Recall from Definition (\ref{def:cr}) that the curves in the $\Crrp$-family are given by 
$\Crrp/\Q:y^2=f_r(x)$, where 
$$
\frrp(x)= (ab)^{(r-1)/2} \, x \, h_r\left(\frac{x^{2}}{ab}+2\right)+b^{r}-a^{r} \text{ and } h_r(x) = \prod_{j = 1}^{\frac{r-1}{2}} (x - \omega_j).
$$
We will compute the cluster pictures and conductor exponents of the curve $C_r$ over $\Q_q$ and $\Qplus_{\mathfrakq}$, where $q$ is odd and  $\mathfrakq$ is a prime ideal above $q$, both of bad reduction for $\Crrp$. 
%%%%%%%%%%%%%%%%%%%%%%%%%%%%%%%%%%%%%%%%%%%%%%%%%%%%%%%%%%%%%%%%%%%%%%%%%%%%%
\subsection{Roots of the defining polynomial \texorpdfstring{$\frrp$}{}}\label{sec: Roots rrp}
%%%%%%%%%%%%%%%%%%%%%%%%%%%%%%%%%%%%%%%%%%%%%%%%%%%%%%%%%%%%%%%%%%%%%%%%%%%%%
Let $q$ be an odd prime. We first exhibit an algebraic expression for the roots of \texorpdfstring{$f_r$}{fr} in \texorpdfstring{$\overline{\Q_q}$}{Qqbar}. Denote by $\zr \in \overline{\Q_q}$ the choice of a primitive $r$-th root of unity. 

\begin{defn}
    For any $j \in \Z$, we define $\gamma_j \coloneqq \zr^{j} \, a - \zr^{-j}  b \in \overline{\Q_q}$.
\end{defn}

\begin{lemma}\label{cor: Dif Roots rrp}
For any $j,k \in \Z$, %with $j \neq k$, 
we have
\begin{equation*}
\gamma_k-\gamma_j=\zeta_r^k(1-\zeta_r^{j-k})(a+\zeta_r^{-j-k}b). 
\end{equation*}
For any $0 \leq k \leq r-1$, we have
    \begin{equation*}
        \prod_{\substack{0 \leq j \leq r-1 \\ j \neq k}} \left( \gamma_k - \gamma_j \right) = \frac{r \left(a^r + b^r \right)}{\zr^k \left(a + \zr^{-2k} b \right)}.
    \end{equation*}
\end{lemma}

\begin{proof} The first equality is obtained by developing the right-hand side. Now for any $0\!\leq\!k\!\leq\!r\!-\!1$, we have
    \begin{equation*}
        \prod_{j \neq k} (\gamma_k - \gamma_j) = \prod_{j \neq k} \zr^{k} (1 - \zr^{j-k}) (a + \zr^{-j-k} b) = \zr^{-k} r \prod_{j \neq k} (a + \zr^{-j-k} b) = \frac{r \left(a^r + b^r \right)}{\zr^k \left(a + \zr^{-2k} b \right)}.
    \end{equation*}
\end{proof}

\begin{proposition}\label{prop:rrproots}
The set of roots of the polynomial $f_r$ is $\Rroots \coloneqq \lbrace \gamma_0, \ldots, \gamma_{r-1} \rbrace$.
\end{proposition} 

\begin{proof} 
Fix $0\leq j \leq r-1$. 
Using Lemma~\ref{lem:phiandh}, we compute
\begin{align*}
\frrp(\gamma_j)&=(ab)^{\frac{r-1}{2}}(\zr^{j}a-\zr^{-j}b) \cdot h_r\left(\frac{(\zr^{j}a-\zr^{-j}b)^{2}}{ab}+2\right) + b^{r}-a^{r} \\
&=(\zr^{j}a\zr^{-j}b)^{\frac{r-1}{2}}(\zr^{j}a-\zr^{-j}b) \cdot h_r\left(2 - \frac{(\zr^{j}a-\zr^{-j}b)^{2}}{-\zr^{j}a\zr^{-j}b}\right) + b^{r}-a^{r}\\
&=(\zr^ja)^r+(-\zr^{-j}b)^r + b^{r}-a^{r}=0. 
\end{align*}
Therefore, $\gamma_0, \ldots, \gamma_{r-1}$ are roots of $f_r$. By assumption $a^r+b^r=c^p \neq 0$, so Lemma~\ref{cor: Dif Roots rrp} yields $\gamma_{k}-\gamma_{j}\neq 0$ for all $0 \leq j, k \leq r-1$, with $j \neq k$. Since $f_r$ has degree $r$, and the $\gamma_{j}$'s are all distinct, we conclude that these are all the roots of $f_r$. 
\end{proof}

%To ease notation, for any $1 \leq j \leq r-1$, we let $\gamma_{-j} \coloneqq \gamma_{r-j}$. 
%Next, we present a useful result concerning the differences 
%of the $\gamma_{j}$'s. 
%of the roots of $\frrp$.

We can easily describe the splitting field of the polynomial $\frrp$.

\begin{lemma}\label{lem: Splitting field frrp}
The splitting field $\Q_q(\Rroots)$ of the polynomial $\frrp$ is $\Q_q(\Rroots) = \Q_q(\zr)$.
\end{lemma}

\begin{proof}
The formula defining each $\gamma_j$ already shows that $\Q_q(\Rroots) \subseteq \Q_q(\zr)$. The reverse inclusion follows from the identity $\zr= \left(a\gamma_1+b\gamma_{-1} \right) / (a^2-b^2) \in \Q_q(\Rroots)$. Note that $a^2 - b^2 \neq 0$, because if that was the case then $\lbrace a, b \rbrace = \lbrace 1, -1 \rbrace$, implying that $c = 0$, a contradiction.
\end{proof}

\subsection{Cluster pictures of \texorpdfstring{$\Crrp$}{}}\label{sec: Clusters rrp}
%%%%%%%%%%%%%%%%%%%%%%%%%%%%%%%%%%%%%%%%%%%%%%%%%%%%%%%%%%%%%%%%%%%%%%%%%%%%%
In this section, we construct the cluster pictures of $\Crrp$ over $\Q_q$ and $\mathcal{K}_{\mathfrak{q}}$ at odd places of bad reduction.

\begin{proposition}\label{eq: DiscriminantCrrp}
The discriminant of the defining equation of $\Crrp$ is given by
$$\Delta(\Crrp)=(-1)^\frac{r-1}{2} \, 2^{2(r-1)} \, r^r \left(a^r+b^r \right)^{r-1}.$$
In particular, the odd places of bad reduction for $C_r$ divide $r(a^r + b^r)$.
\end{proposition}

\begin{proof} Since $r$ is odd, recall that $\Delta(\Crrp) =2^{4g}\Delta(f_{r})$, with $g = \frac{r-1}{2}$ (see \cite[Definition 1.6]{Lockhart}). 
From Lemma~\ref{cor: Dif Roots rrp} we have 
$$ \Delta(f_{r})= 
(-1)^{r\frac{r-1}{2}}\prod_{0 \leq k \leq r-1}\prod_{\substack{0 \leq j \leq r-1 \\ j \neq k}}(\gamma_{k}-\gamma_{j}) = (-1)^{\frac{r-1}{2}}\prod_{0 \leq k \leq r-1}\frac{r \left(a^r + b^r \right)}{\zr^k \left(a + \zr^{-2k} b \right)}.$$ 
Since $\prod_{0 \leq k \leq r-1} \zr^{k} = 1$ and $\prod_{0 \leq k \leq r-1} \left(a + \zr^{-2k} b \right)=a^r+b^r$, the result follows. 
\end{proof}

We note that the proof of Proposition~\ref{eq: DiscriminantCrrp} provides an alternative proof to \cite[(3.6), pg.15]{BCDF}.

\begin{theorem}\label{thm: Clusters rrp Q}
 Let $q$ be an odd prime number of bad reduction for $\Crrp/\Q$. 
\begin{enumerate}
    \item If $q \neq r$ and $q \mid (a^r + b^r)$, then the cluster picture of $\Crrp / \Q_q$ is
    \begin{center}
   \clusterpicture           
  \Root[D] {2} {first} {r2};
  \Root[D] {4} {r2} {r3};
  \ClusterLDName c1[][n][] = (r2)(r3);
  \Root[D] {2} {c1} {r4};
  \Root[D] {4} {r4} {r5};
  \ClusterLDName c2[][n][] = (r4)(r5);
  \Root[E] {6} {c2} {r6};
  \Root[E] {2} {r6} {r7};
  \Root[E] {2} {r7} {r8};
  \Root[D] {7} {r8} {r9};
  \Root[D] {4} {r9} {r10};
  \ClusterLDName c3[][n][] = (r9)(r10);
  \Root[D] {4} {c3} {r11};
   \ClustercLDName c5[][][] = (r11);
  \ClusterLD c4[][\, 0] = (c1)(c1n)(c2)(c2n)(r6)(r7)(r8)(c3)(c3n)(c5)(c5n);
\endclusterpicture
\end{center}
 where $n = v_{q}(a^r+b^r)$. If we denote the isolated root by $\gamma_{i_0}$ then each twin consists of the pair of roots $\{\gamma_{k}, \gamma_{2i_0-k}\}$, for $0\leq k \leq \frac{r-1}{2}$, $k \not \equiv i_0 \mod r$.
 \item If $q = r$ and $r\nmid (a^r+b^r)$ then the cluster picture of $\Crrp / \Q_r$ is
 \begin{center}
    \clusterpicture            
        \Root[D] {2} {first} {r2};
        \Root[D] {4} {r2} {r3};
        \Root[D] {4} {r3} {r4};
        \Root[E] {4} {r4} {r5};
        \Root[E] {1} {r5} {r6};
        \Root[E] {1} {r6} {r7};
        \Root[D] {4} {r7} {r8};
        \Root[D] {4} {r8} {r9};
        \Root[D] {4} {r9} {r10};
        \ClusterLD c1[][\, \dfrac{1}{r-1}] = (r2)(r3)(r4)(r5)(r6)(r7)(r8)(r9)(r10);
        \endclusterpicture.
\end{center}

    \item If $q = r$ and $r\mid (a^r+b^r)$ then the cluster picture of $\Crrp / \Q_r$ is
    \begin{center}
     \clusterpicture         
  \Root[D] {2} {first} {r2};
  \Root[D] {4} {r2} {r3};
  \ClusterLDName c1[][m][\gamma_1, \gamma_{-1}] = (r2)(r3);
  \Root[D] {2} {c1} {r4};
  \Root[D] {4} {r4} {r5};
  \ClusterLDName c2[][m][\gamma_2, \gamma_{-2}] = (r4)(r5);
  \Root[E] {6} {c2} {r6};
  \Root[E] {2} {r6} {r7};
  \Root[E] {2} {r7} {r8};
  \Root[D] {7} {r8} {r9};
  \Root[D] {4} {r9} {r10};
  \ClusterLDName c3[][m][\gamma_{\frac{r-1}{2}}, \gamma_{\frac{r+1}{2}}] = (r9)(r10);
  \Root[D] {4} {c3} {r11};
   \ClustercLDName c5[][][\gamma_0] = (r11);
  \ClusterLD c4[][\, \dfrac{2}{r-1}] = (c1)(c1n)(c2)(c2n)(r6)(r7)(r8)(c3)(c3n)(c5)(c5n);
\endclusterpicture
\end{center}
where $m = v_r(a+b)-\frac{1}{r-1}$.
\end{enumerate}
\end{theorem}

\begin{proof}
    Recall that $a$ and  $b$ are coprime, so that $v_q(a) v_q(b) = 0$ for any odd prime $q$.
    
    $(1)$
    %\mar{Suppose $q \neq r$ and $q\mid (a^r + b^r)$. Since $a$ and $b$ are coprime, Lemma~\ref{lem:generalproduct} implies that $v_q(a^r + b^r) = v_q(a+b)$ and $v_q(a + \zr^j b) = 0$ for any $1 \leq j \leq r-1$. Combining this with Corollary \ref{cor: Dif Roots rrp} yields $ v_q(\gamma_0 - \gamma_j) = 0$ for any $1 \leq j \leq r-1$.} 
    Suppose $q \neq r$ and $q\mid (a^r + b^r)$. Lemma~\ref{lem:generalproduct}, applied with $K = \Q_q$, $\alpha = a$, $\beta = b$, states that there exists some $0 \leq j_0 \leq r-1$ such that $v_q(a + \zr^{j_0} b) = v_q(a^r + b^r)$ and $v_q(a + \zr^k b) = 0$ for $k \neq j_0$. Since $q \neq r$, Lemma~\ref{cor: Dif Roots rrp} implies that, for any $0 \leq j \leq r-1$ such that $j \not \equiv k \mod r$, we have $v_q(\gamma_k - \gamma_j) > 0$ if and only if $j \equiv -k - j_0 \mod r$. It follows that
    \begin{equation}\label{eq:twins}
        v_q (\gamma_{k} - \gamma_{-k-j_0}) =v_q(a + \zr^{j_0} b) = v_q(a^r + b^r).
    \end{equation}
    In particular, letting $i_0$ be such that $-2 i_0 \equiv j_0 \mod r$, any choice of $k \neq i_0$ satisfies $- k - i_0 \not \equiv j_0 \mod r$, implying that $v_q(\gamma_k - \gamma_{i_0}) = 0$. We deduce that $\gamma_{i_0}$ is an isolated root in the cluster picture. Moreover,  
     %Moreover, for any $1 \leq k \leq \frac{r-1}{2}$, we have $k - (2i_0 - k) \equiv j_0 \mod r$.
     %so $v_q(\gamma_{k} - \gamma_{2i_0 - k}) > 0$. 
     %We can explicitly compute this valuation:
    %The conventions made in Remark~\ref{rmk:conventionrrp} yield $ v_q(\gamma_0 - \gamma_j) = 0$ for any $1 \leq j \leq r-1$. Thus $\gamma_0$ is an isolated root and does not belong to any cluster other than $\Rroots$. 
    %Similarly, for any $1 \leq k \leq r-1$, Lemma \ref{lem:generalproduct} combined with the convention made in Remark~\ref{rmk:conventionrrp} asserts that $v_q( a + \zr^{-k-j} b) > 0$ if and only if $- k  -j \equiv 0 \bmod \, r$. 
    % \begin{equation*}
    %     v_q (\gamma_{k} - \gamma_{2i_0 - k}) = v_q \left( \prod_{\substack{0 \leq j \leq r-1 \\ j \not \equiv i_0 + k \mod r}} (\gamma_{i_0 + k} - \gamma_j) \right) = v_q \left( \frac{r (a^r + b^r)}{\zr^k(a + \zr^{-2k} b)} \right) = v_q(a^r + b^r).
    % \end{equation*}
  there are $\frac{r-1}{2}$ twins of depth $v_q(a^r + b^r)$, each of them consisting of the roots $\lbrace \gamma_{k}, \gamma_{2i_0 - k} \rbrace$ where $k\neq i_0$. This is because the map $f: k \mapsto 2i_0 -k$ is an involution on $\F_r \setminus \{ i_0\}$ without fixed points, and thus it induces $\frac{r-1}{2}$ orbits $\{k, f(k)\}=\{k, 2i_0-k\}$ which satisfy $v_q(\gamma_k-\gamma_{2i_0-k})=n$ by \eqref{eq:twins}.

    % We therefore obtain $\frac{r-1}{2}$ twins , giving the desired cluster picture.

    $(2)$ Suppose $q= r$ and $r \nmid (a^r + b^r)$. Then, for any $0 \leq j, k \leq r-1$ with $j \neq k$, we have $v_r (a + \zr^{-j -k} b) = 0$. This implies that $v_r(\gamma_k - \gamma_j) = v_r \left( 1 - \zr^{j-k} \right) = \frac{1}{r-1}$. We deduce that all the roots lie in a common cluster whose depth is $\frac{1}{r-1}$.

    $(3)$ Suppose $q = r$ and $r \mid (a^r + b^r)$. Then $v_r(a^r + b^r) \geq 2$, because $a^r + b^r = c^p$, and $v_r (1 - \zr) = \frac{1}{r-1}$, therefore it holds that $v_r(a^r + b^r) > r v_r(1 - \zr)$. We apply Lemma \ref{lem:generalproduct} to $K = \Q_r$, $\alpha=a$ and  $\beta=b$. Since $\zr \notin \Q_r$, we have that $j_0 = 0$ and therefore $v_r \left( a+ \zr^j b \right) = v_r (1 - \zr) = \frac{1}{r-1}$ for any $1 \leq j \leq r-1$. Moreover
    $$
    v_r( \gamma_0 - \gamma_j) = v_r((1 - \zr^{j}) (a + \zr^{-j} b)) = 2 v_r(1 - \zr) = \tfrac{2}{r-1}.
    $$
    Suppose now that $k\neq 0$.
    Then, for any $j \neq k$, we have
    \begin{equation*}
		v_r(\gamma_k - \gamma_j) = v_r((1 - \zr^{j-k}) (a + \zr^{-j-k} b)) = \begin{cases}
			\frac{2}{r-1} & \text{ if } j \not \equiv -k \bmod\:r, \\
			\frac{1}{r-1} + v_r(a + b)  & \text{ if } j \equiv -k \bmod\:r.
		\end{cases}
	\end{equation*}
	Therefore all the roots belong to a cluster having depth $\frac{2}{r-1}$ and the root $\gamma_0$ is isolated. Further, there are $\frac{r-1}{2}$ twins, and each twin has relative depth $v_r(a + b) - \frac{1}{r-1}$.
\end{proof}

\begin{cor}\label{cor: Clusters Crrp Qplus}
Let $\mathfrakq$ be an odd place of $\Qzrplus$ lying above the rational prime $q$. The cluster picture of $C_r / \Qzrplus_{\mathfrakq}$ is obtained from the cluster picture of $C_r / \Q_q$ as in Theorem~\ref{thm: Clusters rrp Q}, by multiplying all depths by the index of ramification $e_{\Qzrplus_{\mathfrakq} / \Q_q} = \begin{cases}
1 & \text{ if } q \neq r, \\
\frac{r-1}{2} & \text{ if } q = r.
\end{cases}$ 

\end{cor}

\begin{proof}
This is an immediate consequence of Theorem~\ref{thm: Clusters rrp Q} and \eqref{eq:valsrelation}. 
\end{proof}

%%%%%%%%%%%%%%%%%%%%%%%%%%%%%%%%%%%%%%%%%%%%%%%%%%%%%%%%%%%%%%%%%%%%
\subsection{Conductor exponents for \texorpdfstring{$\Crrp$}{} at odd places}\label{sec: Conductors rrp}
%%%%%%%%%%%%%%%%%%%%%%%%%%%%%%%%%%%%%%%%%%%%%%%%%%%%%%%%%%%%%%%%%%%%%%%%%%%%%
In this section, we compute the conductor exponent of $\Crrp$ over $\Q_{q}$ and $\mathcal{K}_\qQzrplus$ at odd places of bad reduction from their cluster pictures. 

\begin{theorem} \label{cor:conductorQ} 
    The conductor exponent of \texorpdfstring{$\Crrp / \Q$}{Crrp/Q} at an odd prime of bad reduction $q$ is 
    $$
      \condexp{\Crrp / \Q_q} = 
      \begin{cases}\frac{r-1}{2} &\text{ if } q \ne r \text{ and } q \mid a^r + b^r,\\
								r-1 &\text{ if } q = r.
\end{cases}
$$
    
\end{theorem}

\begin{proof}   
    \textbf{Case 1:} Assume first $q \neq r$ and $q\mid (a^r + b^r)$. We check that the cluster picture of $\Crrp / \Q_q$, depicted in Theorem~\ref{thm: Clusters rrp Q} $(1)$, satisfies the semistability criterion (Definition \ref{def:ss}) so that $\Crrp / \Q_q$ is semistable by Theorem~\ref{thm:ss}. Since $q \neq r$, Lemma~\ref{lem: Splitting field frrp} implies that the extension $\Q_q (\Rroots) / \Q_q$ is unramified and the inertia group $I_{\Q_{q}}$ acts trivially on $\mathcal{R}$ and on every twin. % $\mathfrak{t}_{k}$ for  
    %$ \coloneqq \lbrace \gamma_{j}, \gamma_{-j} \rbrace$, 
   % $1 \leq k \leq \frac{r-1}{2}$.
     Moreover, the only principal cluster is $\mathcal{R}$, which has depth $d_{\Rroots} = 0$, and using the definition of $\nu_{\Rroots}$ in Definition \ref{defclusterinv}, one checks that $\nu_{\Rroots} = 0$. Hence, $\Crrp / \Q_q$ is semistable and Theorem~\ref{thm:condss} implies that $\condexp{\Crrp / \Q_q}$ equals the number of twins appearing in its cluster picture, \ie, $\frac{r-1}{2}$.

    For the next two cases, we will use Theorem \ref{thm: allcond} to compute the conductor exponent.
    
    \textbf{Case 2:} Assume now $q=r$ and $r \nmid (a^r+b^r)$. Note first that $\Q_r(\Rroots) / \Q_r$ is a tame extension, so Proposition~\ref{cor:nwild0} implies that $\condwild{\Crrp / \Q_r} = 0$. It remains to compute $\condtame{\Crrp / \Q_r}$. The cluster picture depicted in Theorem~\ref{thm: Clusters rrp Q} $(2)$ shows that the candidates for $U$ are odd clusters $\neq \Rroots$, which are the singletons $\lbrace \gamma_j \rbrace$. Fix $0 \leq j \leq r-1$, then $P(\lbrace \gamma_{i} \rbrace)=\mathcal{R}$, whose depth is $d_{\Rroots}=\frac{1}{r-1}$, and since the inertia group stabilizes $\Rroots$ we have $[I_{\Q_{r}}:I_{\Rroots}]=1$. One can compute $\tilde{\lambda}_{\Rroots}=\frac{r}{2(r-1)}$, and then we have 
    \begin{equation*}
        \xi_{\Rroots}(\tilde{\lambda}_{\Rroots}) = v_2(2(r-1)) > v_2(r-1) = \xi_{\Rroots}(d_{\Rroots}).
    \end{equation*}
    Thus $\lbrace \gamma_j \rbrace \notin U$ for any $0 \leq j \leq r-1$, so $U = \emptyset$. To compute $V$, note that the only proper non-\"ubereven cluster is $\mathcal{R}$. But as we just saw $\xi_{\Rroots}(\tilde{\lambda}_{\Rroots}) = v_2(2(r-1)) \neq 0$, so $\Rroots \notin V$ and therefore $V = \emptyset$. Since $\Crrp$ has genus $\frac{r-1}{2}$ and $\frrp$ has odd degree, we conclude that $\condtame{\Crrp / \Q_r} = r-1$.

    \textbf{Case 3:} Assume now $q = r$ and $r\mid (a^r+b^r)$. Again $\Q_r(\Rroots) / \Q_r$ is a tame extension, so $\condwild{\Crrp / \Q_r} = 0$.  We focus on computing $\condtame{\Crrp / \Q_r}$. For $ 1 \leq j \leq \frac{r-1}{2}$, we denote by $\mathfrakt_j$ the twin $\mathfrakt_j \coloneqq \left \lbrace \gamma_j, \gamma_{-j} \right \rbrace$. 
    Just as in the previous case, the odd clusters $\neq \Rroots$ (candidates for $U$) are the singletons $\lbrace \gamma_j \rbrace$. We first consider $\lbrace \gamma_0 \rbrace$, whose parent is $\Rroots$, which has depth $d_{\Rroots}=\frac{2}{r-1}$, and again $[I_{\Q_{r}}:I_{\Rroots}]=1$. This time one computes $\tilde{\lambda}_{\Rroots}=\frac{1}{r-1}$, so $\xi_{\Rroots}(\tilde{\lambda}_{\Rroots}) > \xi_{\Rroots}(d_{\Rroots})$, and therefore we have $\lbrace \gamma_0 \rbrace \notin U$. Fix now $1 \leq j \leq r-1$. The parent of $\lbrace \gamma_j \rbrace$ is $\mathfrakt_{j}$ if $1 \leq j \leq \frac{r-1}{2}$ and $\mathfrakt_{r-j}$ if not, which has depth $d_{\mathfrakt_j}=\frac{1}{r-1} + v_r(a + b)$. For simplicity assume that $1 \leq j \leq \frac{r-1}{2}$. We compute $ \tilde{\lambda}_{\mathfrakt_j} = v_{r}(a+b) + 1$.
  
    The explicit expression of $\gamma_j$ and $\gamma_{-j}$ shows that the stabilizer of $\mathfrakt_j$ under inertia is the set $I_{\mathfrakt_j}=\{\sigma \in I_{\Q_r} : \sigma (\zr + \zr^{-1}) = \zr + \zr^{-1}\}$, whose index in $I_{\Q_r}$ is 
    \begin{equation*}
        \left[ I_{\mathbb{Q}_r} : I_{\mathfrakt_j} \right] = \left \mid \Gal(\unr{\Q_r}(\zr+\zr^{-1})/\unr{\Q_r}) \right \mid =\frac{r-1}{2}.
    \end{equation*}
   One computes 
    \begin{equation*}
        v_2 \left( \left[ I_{\mathbb{Q}_r} : I_{\mathfrakt_j} \right] d_{\mathfrakt_j} \right) = v_2 \left( \frac{r-1}{2} \left( \frac{1}{r-1} + v_r(a+b) \right) \right) = -1,
    \end{equation*}
    so $\xi_{\mathfrakt_j} \left(d_{\mathfrakt_j} \right) = 1$. On the other hand we have $v_2 \left( \tilde{\lambda}_{\mathfrakt_j}\right) \geq 0$ so $\xi_{\mathfrakt_j} \left( \tilde{\lambda}_{\mathfrakt_j}\right) = 0$. Therefore we have $U = \left \lbrace \lbrace \gamma_j \rbrace ; 1 \leq j \leq r-1 \right \rbrace$. The candidates for $V$ are $\Rroots$ and $\mathfrakt_{1}, \ldots, \mathfrakt_{\frac{r-1}{2}}$. From the discussion above, we have $\xi_{\Rroots}(\tilde{\lambda}_{\Rroots}) > 0$ and $\xi_{\mathfrakt_j}(\tilde{\lambda}_{\mathfrakt_j}) = 0$ for any $1 \leq j \leq \frac{r-1}{2}$, hence, $V = \left \lbrace \mathfrakt_j ; 1 \leq j \leq \frac{r-1}{2} \right \rbrace$. Since $I_{\Q_r}$ contains the maps $\zr \mapsto \zr^{j}$ for $1 \leq j \leq r-1$, $I_{\Q_r}$ acts transitively on $U$ and $V$ giving $\# (U/I_{\Q_r})= \# (V/I_{\Q_r})=1$. We conclude that $\condtame{\Crrp / \Q_r} = r-1$.
\end{proof}

%\newpage
\begin{theorem} \label{thm:conductorQplus}
The conductor exponent $\Crrp / \Qplus$ at an odd place of bad reduction $\mathfrakq$ is 
    $$
      \condexp{\Crrp / \Qplus_{\mathfrakq}} = 
      \begin{cases}\frac{r-1}{2} & \text{ if } \mathfrakq \nmid r \text{ and } \mathfrak{q} \mid a^r + b^r ,\\
	   r-1 & \text{ if } \mathfrakq \mid r.
\end{cases}
$$
\end{theorem}

\begin{proof}
We know from \eqref{eq: Isom compl+ext} that $\Qplus_{\qQzrplus} \simeq \Q_q(\omega)$. Moreover, by Lemma~\ref{lem: Splitting field frrp}, $\omega$ belongs to the splitting field $\Q_q(\Rroots)$, so we can identify $\Qplus_{\qQzrplus}(\Rroots)$ as an extension of $\Q_q(\omega)$. Recall that $\mathfrak r$ is the unique prime above $r$ in $\Qplus$.

\textbf{Case 1:} Assume first $\mathfrakq \nmid r$, the curve $\Crrp/\Q_q$ is semistable by the proof of Theorem~\ref{cor:conductorQ}, Case 1.
Hence $\Crrp/\Qplus_\mathfrakq$ is also semistable by \cite[Proposition 3.15]{Liu}.  Theorem~\ref{thm:condss} implies that $\condexp{\Crrp / \Qplus_{\mathfrakq}}$ equals the number of twins appearing in its cluster picture, \ie , $\frac{r-1}{2}$.

In the next two cases, we use Theorem \ref{thm: allcond} to compute conductor exponents. As $\Qplus_{\rQzrplus}/\Q_r$ is tamely ramified, the computation of $\condwild{\Crrp / \Qplus_{\rQzrplus}}$ follows immediately from the proof of  Theorem~\ref{cor:conductorQ}, cases 2 and 3, and Lemma~\ref{lem:wildthoughtame}, therefore we focus on the computation of $\condtame{\Crrp / \Qplus_{\rQzrplus}}$.

\textbf{Case 2:} Assume now $\mathfrakq \mid r$ (so that $\mathfrakq = \mathfrak r$) and $r \nmid (a^r+b^r)$.  
%By Lemma~\ref{lem: Splitting field frrp} we have $\Qplus_{\rQzrplus}(\mathcal{R})=\Q_r(\zr)$, which is a tamely ramified extension of $\Q_r(\omega)$, so $\condwild{\Crrp / \Qplus_{\rQzrplus}}=0$ by Proposition~\ref{cor:nwild0}. It remains to compute $\condtame{\Crrp / \Qplus_{\rQzrplus}}$.
Corollary~\ref{cor: Clusters Crrp Qplus} shows that the potential elements of $U$ are the singletons $\lbrace \gamma_j \rbrace$ for $0\leq j\leq r-1$. Fix $0\leq j\leq r-1$. 
Note that $P(\lbrace \gamma_j \rbrace)=\Rroots$ whose depth is $d_{\Rroots}=\frac{1}{2}$ and, since inertia stabilizes $\Rroots$, we have $[I_{\Qrplus}:I_\Rroots]=1$. One can compute $\tilde{\lambda}_{\Rroots}=\frac{r}{4}$, so $\xi_{\Rroots}(\tilde{\lambda}_{\Rroots}) = 2$ and $\xi_{\Rroots}(d_{\Rroots}) = 1$. Thus none of the singletons $\lbrace \gamma_j \rbrace$ belongs to $U$ and hence $U=\emptyset$. The only potential element of $V$ is $\mathcal{R}$, but as we just saw $\xi_{\Rroots}(\tilde{\lambda}_{\Rroots}) \neq 0$, so $\Rroots \notin V$, and thus $V=\emptyset$. We conclude that $\condtame{\Crrp / \Qplus_{\rQzrplus}} = r-1$. 

\textbf{Case 3:} Assume now $\mathfrakq \mid r$ (so that $\mathfrakq = \mathfrak r$) and $r \mid (a^r+b^r)$. 
%Again $\Qplus_{\rQzrplus}(\Rroots) / \Qplus_{\rQzrplus}$ is a tame extension, so $\condwild{\Crrp / \Qplus_{\rQzrplus}} = 0$. We focus on computing $\condtame{\Crrp / \Qplus_{\rQzrplus}}$.
For any $ 1 \leq j \leq \frac{r-1}{2}$, we denote by $\mathfrakt_j$ the twin $\mathfrakt_j \coloneqq \left \lbrace \gamma_j, \gamma_{-j} \right \rbrace$. Again by Corollary~\ref{cor: Clusters Crrp Qplus}, the potential elements of $U$ are the singletons $\lbrace \gamma_j \rbrace$ for $0\leq j \leq r-1$. First consider $\lbrace \gamma_{0} \rbrace$, whose parent is $\Rroots$, which has depth $d_{\Rroots}=1$ and again $[I_{\Qrplus}:I_\Rroots]=1$. One can compute $\tilde{\lambda}_{\cR}=\frac{1}{2}$, so $\xi_{\Rroots}(\tilde{\lambda}_{\Rroots}) = 1$ and $\xi_{\Rroots}(d_{\Rroots}) = 0$, so we deduce that $\gamma_0 \notin U$. Fix $1 \leq j \leq \frac{r-1}{2}$ (this is without loss of generality because $P(\lbrace \gamma_j \rbrace )= P(\lbrace \gamma_{-j} \rbrace)$. Then $P(\lbrace \gamma_j \rbrace )=\mathfrakt_j$, which has depth $d_{\mathfrakt_j}=\frac{1}{2}((r-1)v_{r}(a+b)+1)$. One can check that $\tilde{\lambda}_{\mathfrakt_j} = \frac{r-1}{2}(v_r(a+b)+1)$. Note that $I_{\Qrplus}$ fixes $\zr+\zr^{-1}$, so it fixes $\mathfrakt_j$ and hence $[I_{\Qrplus} : I_{\mathfrakt_j}] = 1$. Therefore $\xi_{\mathfrakt_j}(d_{\mathfrakt_j})=1$ and $\xi_{\mathfrakt_j}(\tilde{\lambda}_{\mathfrakt_j})=0$. We conclude that $U = \left \lbrace \lbrace \gamma_j \rbrace ; 1 \leq j \leq r-1 \right \rbrace$. The potential elements of $V$ are $\Rroots$ and the $\mathfrakt_j$'s. It follows from the above discussion that $V=\{\mathfrakt_j: 1\leq j\leq \frac{r-1}{2}\}$. The orbits of $\Rroots$ under the action of inertia $I_{\Qrplus}$ are the pairs $\{\gamma_j, \gamma_{-j}\}$, so $\#(U/I_{\Qrplus})=\#(V/I_{\Qrplus})=\frac{r-1}{2}$. We conclude that $\condtame{\Crrp / \Qplus_{\rQzrplus}} =r-1$.
\end{proof}

\begin{remark}
One could describe the reduction of the Néron model of the Jacobian of $C_r / \Qplus_{\mathfrakq}$ \textit{à la} \cite{M2D2}, and then deduce from this the computation of $\condtame{C_r / \Qplus_{\mathfrakq}}$ using \cite[Exposé IX, \S 4]{SGA}. When $\mathfrakq \nmid 2r$, it has totally toric reduction, and when $\mathfrakq \mid r$, this is totally unipotent. Since the unipotent rank can only decrease when base-changing to field extensions, one could also recover from this the computation $\condtame{C_r / \Q_r} = r-1$.
\end{remark}

When using Frey hyperelliptic curves in the modular method, one wishes to minimize the conductor exponent. Using the cluster picture of $\Crrp$, one can find a curve with lower conductor exponent.

\begin{cor}\label{cor:twistCr}
    When $r \mid (a^r + b^r)$, the twist of $\Crrp$ by a uniformizer of $\Qplus_{\rQzrplus}$ is semistable and has conductor exponent $\frac{r-1}{2}$.
\end{cor}

\begin{proof} Recall that $\mathfrak r$ denotes the unique prime ideal above $r$ in $\Qplus$, and let $\pi$ be a uniformizer of $\Qplus_{\rQzrplus}$. We describe the quadratic twist of $\Crrp$ by $\pi$ by the model $y^2 = \pi f_r(x)$. The set of roots of the RHS polynomial is still $\Rroots$, so the cluster picture of the twist is the same as that of $C_r$, but the leading coefficient of the defining polynomial has now valuation $1$. In particular, the semistability criterion (Definition \ref{def:ss}) is now satisfied. Indeed, recall that $\Qplus_{\rQzrplus}(\Rroots) \simeq \Q_r(\zr)$, hence $e_{\Qplus_{\rQzrplus}(\Rroots) / \Qplus_{\rQzrplus}} = 2$. We also know from the proof of Theorem~\ref{thm: Clusters rrp Q} that $I_{\Qplus_{\rQzrplus}}$ does not permute the twins. Lastly, since the valuation of the leading term is one, we check that $\nu_{\Rroots}$ for the twist is even. The computation of the conductor exponent follows from Theorem~\ref{thm:condss}.
\end{proof}

%%%%%%%%%%%%%%%%%%%%%%%%%%%%%%%%%%%%%%%%%%%%%%%%%%%%%%%%%%%%%%%%%%%%%%%%%%%%%%%%%%%%%%%%%%%%%%%%%%%%%%%%%%%%%%%%%%%%%%%%%%%%%%%%%%%%%%%%%%%%%%%%%%%%%%%%%%
\section{\large \texorpdfstring{$\Cminus$}{} and \texorpdfstring{$\Cplus$}{} families associated to the signature \texorpdfstring{$(p, p, r)$}{}}\label{part:ppr}
%%%%%%%%%%%%%%%%%%%%%%%%%%%%%%%%%%%%%%%%%%%%%%%%%%%%%%%%%%%%%%%%%%%%%%%%%%%%%
%%%%%%%%%%%%%%%%%%%%%%%%%%%%%%%%%%%%%%%%%%%%%%%%%%%%%%%%%%%%%%%%%%%%%%%%%%%%%
In this section, we consider the families of Frey hyperelliptic curves associated to the equation of signature $(p, p, r)$. Let $(a, b, c) \in \Z^3$ satisfying $abc\ne 0$, $\gcd(a,b,c)=1$ and $a^p+b^p = c^r$. 
Recall from Definition \ref{def: Cminus & Cplus} that the curves in the families $\Cminus$ and $\Cplus$ are  given by $\Cminus / \Q : y^2 = g_r^-(x)$ and $\Cplus / \Q : y^2 =g_r^+(x)$, where
\begin{align*}
     \gminus(x) & =(-1)^{\frac{r-1}{2}} c^{r-1} x \, h_r \left( 2 - \frac{x^2}{c^2}\right) -2 \left(a^p - b^p \right), \\
    \gplus(x) & = \left( (-1)^{\frac{r-1}{2}} c^{r-1} x \, h_r \left( 2 - \frac{x^2}{c^2}\right) -2 \left(a^p - b^p \right) \right) (x+ 2c),
\end{align*} 
and $h_r(x) = \prod_{j = 1}^{\frac{r-1}{2}} (x - \omega_j)$.
%\diana{where $h_r=..$. This part is asymmetric to rrp.} 
We will compute the cluster pictures and conductor exponents of the curves $\Cminus$ and $\Cplus$ over $\Q_q$ and $\Qplus_{\mathfrakq}$, where $q$ is odd and  $\mathfrakq$ is a prime ideal above $q$, both of bad reduction for $\Cminus$ or $\Cplus$.

%%%%%%%%%%%%%%%%%%%%%%%%%%%%%%%%%%%%%%%%%%%%%%%%%%%%%%%%%%%%%%%%%%%%%%%%%%%%%
\subsection{Roots of the defining polynomials \texorpdfstring{$\gminus$ and $\gplus$}{gminus and gplus}}\label{sec: Roots ppr}
%%%%%%%%%%%%%%%%%%%%%%%%%%%%%%%%%%%%%%%%%%%%%%%%%%%%%%%%%%%%%%%%%%%%%%%%%%%%%
Let $q$ be an odd prime. In this section, we first exhibit an expression for the roots of $\gminus$ and $\gplus$ in $\overline{\Q_q}$. We will compute the (products of) differences of the roots of $\gminus$ and $\gplus$, and the discriminants of $\Cminus$ and $\Cplus$.
%Finally, we will discuss the reducibility of $\gminus$ over $\Q_q$. 

%For any $\alpha \in \overline{\Q_q}$, we denote by $\sqrt[r]{\alpha}$ a choice of an $r$-th root of $\alpha$.
%For any $a \in \overline{\Q_q}$, we will denote by $\sqrt[r]{a}$ a fixed choice of an $r$-th root of $a$, i.e. a choice of a root of $x^r-a \in \overline{\Q_q}[x]$. Whenever $x^r-a$ is reducible, this choice will be specified. 

\begin{defn}\label{defn: Rootsppr}
Let $i \in \overline{\Q_q}$ be a square root of $-1$, $\zr \in \overline{\Q_q}$ be a primitive $r$-th root of unity, and $\sqrt{a}$, $\sqrt{b}$ be square roots of $a$ and $b$ respectively. In order to ease notation, denote by $\sqrt{-ab} \coloneqq i \sqrt{a} \sqrt{b}$. Define $\alpha_0$ to be an $r$-th root of $\left( \sqrt{a}^p + i \sqrt{b}^p \right)^2$, and let $\beta_0 \coloneqq  c^2 / \alpha_0$. For any $0 \leq j \leq r-1$, define $\alpha_j \coloneqq \zr^j \alpha_0$ and $\beta_j \coloneqq \zr^{-j} \beta_0$. 
%For simplicity, we let $\alpha_j \coloneqq \zr^j \alpha_0$ and $\beta_j \coloneqq \zr^{-j} \beta_0$ for $0 \leq j \leq r-1$.

%Moreover, we let $\alpha_j \coloneqq \zr^j \, \alpha_0$ and $\beta_j \coloneqq \zr^{-j} \beta_0$.
\end{defn}

\begin{lemma}\label{lem: useful eqs}
For any $0 \leq j \leq r-1$, we have the following properties :
\begin{enumerate}
\item The element $\beta_0$ is an $r$-th root of $\left( \sqrt{a}^p - i \sqrt{b}^p \right)^2$,
%\item $\alpha_0\beta_0= c^2$,
\item $\alpha_j^r+\beta_j^r= 2(a^p-b^p) $,
\item $\alpha_j \beta_j = c^2$, 
\item $\alpha_j^r-\beta_j^r= 4i\sqrt{a}^p\sqrt{b}^p $,
\item $\alpha_j^r+c^r=2\sqrt{a}^p(\sqrt{a}^p+i\sqrt{b}^p)$,
%\item $(\alpha_0+\zr^kc)(\alpha_0-\zr^kc)=\alpha_0(\alpha_0-\zr^{2k}\beta_0)$ for any $k\in \Z$.
\item 
%$\alpha_j\beta_j=c^2$, and 
for any odd prime number $q$ and any $0 \le j \le r-1$, we have $v_q(\alpha_j)v_q(\beta_j) = 0$.

\end{enumerate}
\end{lemma}
\begin{proof}
Statements \textit{(1) - (5)} are algebraic manipulations that follow from Definition \ref{defn: Rootsppr}. For statement \textit{(6)}, assume by contradiction that $v_q(\alpha_j), v_q(\beta_j)$ are both positive. Then $v_q(\alpha_j^r - \beta_j^r) = p v_q(ab)/2 > 0$, and $v_q(\alpha_j^r + \beta_j^r) = v_q(a^p - b^p) > 0$, which contradicts $a, b$ being coprime.
%follows from $v_q(\alpha_0^r) = 2 v_q(\sqrt{a}^p + i \sqrt{b}^p) = $, as $a, b$ are coprime. The same reasoning applies to $\beta_0$.
\end{proof}

\begin{remark}\label{rmk: alpha0 in Qq}
    By \cite[Section VI, Theorem 9.1]{Lang}, the polynomial $x^{r} - (\sqrt{a}^p+i\sqrt{b}^p)^{2} \in \Q_{q}(\sqrt{-ab})[x]$ is reducible if and only if it has a root in $\Q_{q}(\sqrt{-ab})$. Whenever this is the case, we adopt the convention that $\alpha_0 \in \Q_{q}(\sqrt{-ab})$ is such a root. Therefore, the degree of $\Q_q(\alpha_0) / \Q_q(\sqrt{-ab})$ is either equal to $1$ or $r$.
\end{remark}

\begin{defn}
For any $j \in \Z$, we define $\gamma_j \coloneqq \alpha_{j} + \beta_{j} \in \overline{\Q_q}$. We also let $\lastroot \coloneqq -2c$.
\end{defn}

We note that, for any $j \in \Z$, we have $\alpha_j = \alpha_{j + r}$, and similarly $\beta_j = \beta_{j + r}$, hence $\gamma_j = \gamma_{j + r}$. Thus, the definition of $\gamma_j$ depends only on the congruence class of $j \mod r$. 

\begin{lemma}\label{cor: Dif 2 Roots ppr}
For any $j,k \in \Z$, 
%with $j \neq k$, 
we have
$$\gamma_{k}-\gamma_{j} = \zeta_{r}^{k}(1-\zeta_{r}^{j-k})(\alpha_{0} - \zeta_{r}^{-k-j}\beta_{0}).$$
For any $0 \leq k \leq r-1$, we have
    \begin{equation*}
        \prod_{\substack{0 \leq j \leq r-1 \\ j \neq k}} (\gamma_k - \gamma_j) = \frac{4 i r \sqrt{a}^p \sqrt{b}^p}{\zr^k (\alpha_0 - \zr^{-2k} \beta_0)}.
    \end{equation*}
\end{lemma}

\begin{proof}
The first equality is obtained by developing the right-hand side.
For the second part, let
$0 \leq k \leq r-1$ and observe that
    \begin{equation*}
        \prod_{\substack{0 \leq j \leq r-1 \\ j \neq k}} (\gamma_k - \gamma_j) = \prod_{\substack{0 \leq j \leq r-1 \\ j \neq k}} \zr^{k} (1 - \zr^{j-k}) (\alpha_0 - \zr^{-j-k} \beta_0) = \zr^{-k} r \prod_{\substack{0 \leq j \leq r-1 \\ j \neq k}} (\alpha_0 - \zr^{-j-k} \beta_0).
    \end{equation*}
    The product in the last term equals  $\frac{\alpha_0^r - \beta_0^r}{\alpha_0 - \zr^{-2k} \beta_0} = \frac{4 i \sqrt{a}^p \sqrt{b}^p}{\alpha_0 - \zr^{-2k} \beta_0}$, by Lemma~\ref{lem: useful eqs} \textit{(4)}.
\end{proof}

\begin{proposition}\label{prop: Rootsppr}
The set of roots of the polynomial $\gminus$ is $\Rroots \coloneqq \lbrace \gamma_0, \ldots, \gamma_{r-1} \rbrace$.
\end{proposition}

\begin{comment}

\begin{proposition}
    \label{prop: Rootsppr}
The set of roots $\cR$ of the polynomial $g_r^-$ is $\cR: = \{\gamma_0, \ldots, \gamma_{r-1}\} \subset \overline{\Q_q}$, where
        $$
        \gamma_{j} = \zr^{j} \, \alpha_{0} + \zr^{-j}  \beta_{0}, \quad \text{ with }$$ 
        $$\alpha_{0}\coloneqq \sqrt[r]{(\sqrt{a^{p}} + i\sqrt{b^{p}})^{2}} \in \overline{\Q_q} \quad  \text{ and } \quad
    \beta_{0} \coloneqq \frac{c^{2}}{\alpha_{0}}=\sqrt[r]{(\sqrt{a^{p}} - i\sqrt{b^{p}})^{2}}\in \overline{\Q_q}.$$
The set of roots $\Rroots^{+} \subset \overline{\Q_q}$ of the polynomial $\gplus$ is $\Rroots^{+} = \Rroots \cup \lbrace -2c \rbrace$.
\end{proposition}
\end{comment}

\begin{proof}
Fix $0 \leq j \leq r-1$.
%For simplicity, write $\alpha_j \coloneqq \zr^j \alpha_0$ and $\beta_j \coloneqq \zr^{-j} \beta_0$.
Using the equality $\alpha_j \beta_j = c^2$, Lemma~\ref{lem:phiandh} and Lemma~\ref{lem: useful eqs} (2), we compute \vspace{-0.5em}
\begin{align*}
    \gminus(\gamma_j) & = (-1)^{\frac{r-1}{2}} c^{r-1} (\alpha_j + \beta_j) \, h_r \left( 2 - \frac{(\alpha_j + \beta_j)^2}{c^2}\right) -2 \left(a^p - b^p \right) \\
    & = (- \alpha_j \beta_j)^{\frac{r-1}{2}} (\alpha_j + \beta_j) \, h_r \left( 2 - \frac{(\alpha_j + \beta_j)^2}{\alpha_j \beta_j}\right) -2 \left(a^p - b^p \right) \\
    & = \alpha_j^r + \beta_j^r - 2(a^p - b^p) \\
    & = 2(a^p - b^p) - 2(a^p - b^p) = 0
\end{align*}
from Lemma~\ref{lem: useful eqs} $2)$. Therefore, $\gamma_0, \ldots, \gamma_{r-1}$ are roots of $\gminus$. We assume that $ab \neq 0$, so $\alpha_0^r - \beta_0^r = 4i\sqrt{a}^p \sqrt{b}^p \neq 0$, and thus Lemma~\ref{cor: Dif 2 Roots ppr} implies that $\gamma_k - \gamma_j \neq 0$ for all $0 \leq j, k \leq r-1$, with $j \neq k$. Since $\gminus$ has degree $r$ and the $\gamma_j$'s are all distinct, we conclude that these are all the roots of $\gminus$. 
\end{proof}

%Throughout the rest of the paper, to ease notation, we will denote by $\gamma_{-j} \coloneqq \gamma_{r-j}$ for any $1 \leq j \leq r-1$. Moreover, we let $\lastroot \coloneqq -2c$, which clearly is a root of $\gplus$. 
We now describe the differences between $\lastroot = -2c$ and the roots of $\gminus$.

\begin{proposition}\label{prop extra root}
For any $0 \leq j \leq r-1$, we have
\begin{equation*}
\gamma_j - \lastroot = \frac{\zr^j}{\alpha_0} \left( \alpha_0 + \zr^{-j} c\right)^2 \qquad \text{ and } \qquad \prod_{j= 0}^{r-1} (\gamma_j - \lastroot) = 4 a^p.
\end{equation*}
In particular, for any $0 \leq j \leq r-1, \, \gamma_j \neq \lastroot$, and so all the roots of $\gplus$ are distinct.
\end{proposition}

\begin{proof}
By definition of $\beta_0$, we have $\alpha_0 \beta_0 = c^2$, so we compute 
\begin{equation*}
\frac{\zr^j}{\alpha_0} \left( \alpha_0 + \zr^{-j} c\right)^2 \ = \ \zr^j \alpha_0 + 2c + \zr^{-j} \frac{c^2}{\alpha_0} \  = \ \zr^j \alpha_0 + 2c + \zr^{-j} \beta_0 \ = \ \gamma_j - \lastroot.
\end{equation*}
The formula for the product of the differences follows from Lemma~\ref{lem: useful eqs} \textit{(5)}. By assumption we have $a \neq 0$, so $\prod_{j= 0}^{r-1} (\gamma_j - \lastroot) \neq 0$, hence the last claim.
\end{proof}

We can easily deduce the discriminants of the curves $\Cminus, \Cplus$.

\begin{proposition}\label{prop:discppr}
The discriminants of the defining equations of $\Cminus$ and $\Cplus$ are given by
%$\Cminus(a,b,c)$ and $\Cplus(a,b,c)$ are given by 
\begin{align*}
\Delta(\Cminus) & = 2^{4(r-1)}r^{r}a^{p(r-1)/2}b^{p(r-1)/2}, \\
\Delta(\Cplus) & = 2^{4r}r^{r}a^{p(r+3)/2}b^{p(r-1)/2}.   
    %\Delta(\Cminus) & = (-1)^{(r-1)/2}2^{2(r-1)}r^{r}a^{p(r-1)/2}b^{p(r-1)/2} \\
    %\Delta(\Cplus) & = (-1)^{(r-1)/2}2^{2(r+1)}r^{r}a^{p(r+3)/2}b^{p(r-1)/2}.
\end{align*}
In particular, the odd places of bad reduction for $\Cminus$ and $\Cplus$ divide $rab$.
\end{proposition}

\begin{proof}
We recall that $\Delta(C_r^{\pm}) =2^{4g}\Delta(g_{r}^{\pm})$, 
%and $\Delta(\Cplus) =2^{4(g+1)}\Delta(\gplus)$,
with $g = \frac{r-1}{2}$. 
From Lemma~\ref{cor: Dif 2 Roots ppr} we have 
$$ \Delta(\gminus)= 
(-1)^{r\frac{r-1}{2}}\prod_{\substack{0 \leq j, k \leq r-1 \\ j \neq k}}(\gamma_{k}-\gamma_{j}) = (-1)^{\frac{r-1}{2}}\prod_{0 \leq k \leq r-1}\frac{4 i r \sqrt{a}^p \sqrt{b}^p}{\zr^k (\alpha_0 - \zr^{-2k} \beta_0)}.$$ 
But we know that 
$$\prod_{0 \leq k \leq r-1} \zr^{k} = 1 \quad \text{ and } \quad \prod_{0 \leq k \leq r-1} \left(\alpha_0 - \zr^{-2k} \beta_0 \right)=\alpha_0^r-\beta_0^r=4i\sqrt{a}^p \sqrt{b}^p,$$ 
by Lemma~\ref{lem: useful eqs} \textit{(4)}, so it follows that 
$$\Delta(\gminus)= 2^{2(r-1)}r^{r}a^{p(r-1)/2}b^{p(r-1)/2},$$
and the result follows. Moreover, we note that $\Delta(\gplus)=\Delta(\gminus)\prod_{j= 0}^{r-1} (\gamma_j - \lastroot)^2 $, and so we deduce the discriminant of $\Cplus$ using Proposition~\ref{prop extra root}.
\end{proof}

\begin{remark}
Note that the discriminant of $\Cplus$ is not symmetric in $a$ and $b$. As we will see in Sections  \ref{sec: Clusters ppr} and \ref{sec: Conductors ppr}, the local behaviour of $\Cplus$ is not the same at places dividing $a$ and places dividing $b$.
\end{remark}

\subsection{Cluster pictures of \texorpdfstring{$\Cminus$}{C-} and \texorpdfstring{$\Cplus$}{C+}}\label{sec: Clusters ppr}
%%%%%%%%%%%%%%%%%%%%%%%%%%%%%%%%%%%%%%%%%%%%%%%%%%%%%%%%%%%%%%%%%%%%%%%%%%%%%

In this section, we construct the cluster pictures of $\Cminus$ and $\Cplus$ over $\Q_q$ and $\Qplus_{\mathfrakq}$ at odd places of bad reduction. 
\medskip 

%Whenever the prime $q$ divides $ab$, Lemma~\ref{lem: useful eqs} \textit{(3)} implies $v_q(\alpha_0^r - \beta_0^r) > 0$. Lemma \ref{lem:generalproduct} implies the existence of some $0 \leq j_0 \leq r-1$ such that $v_q(\alpha_0 - \zr^{j_0} \beta_0) = v_q(\alpha_0^r - \beta_0^r)$ 

\begin{defn}\label{def: i0 j0 ppr}
    If $q \neq r$ and $q \mid ab$, we denote by $0 \leq j_0 \leq r-1$ the unique integer satisfying $v_q(\alpha_0 - \zr^{j_0} \beta_0) = v_q(\alpha_0^r - \beta_0^r)$ and $v_q(\alpha_0 - \zr^j \beta_0) = 0$ for $j \neq j_0$. We let $0 \leq i_0 \leq r-1$ be the unique integer satisfying $-2 i_0 \equiv j_0 \mod r$. If $q = r$, we define $i_0 = j_0 = 0$.
\end{defn}

\begin{remark}\label{rmk: i0 j0 ppr}
When $q \neq r$, Lemma~\ref{lem: useful eqs} \textit{(4)} and \textit{(6)} imply $v_q(\alpha_0^r - \beta_0^r) > 0$ and $v_q(\alpha_0) v_q(\beta_0) = 0$. The existence of $j_0$ is then guaranteed by Lemma~\ref{lem:generalproduct} applied with $K = \Q_q$, $\alpha = \alpha_0$, $\beta = - \beta_0$. We note that, whenever $\zr \notin \Q_q$, then we have $j_0 = i_0 = 0$.

When $q = r$, we have $\zr \notin \Q_r$ and $v_r( \alpha_0^r - \beta_0^r) = \frac{p}{2} \, v_r(ab) > \frac{3}{2} \geq r v_r(1 - \zr) = \frac{r}{r-1}$. We can also apply Lemma~\ref{lem:generalproduct}, which states that $v_r(\alpha_0 - \beta_0) > \frac{1}{r-1}$ and $v_r(\alpha_0 - \zr^j \beta_0) = \frac{1}{r-1}$ for $1 \leq j \leq r-1$. 

Whether $q = r$ or not, we have $v_q(\alpha_0 + \zr^{j_0} \beta_0) > v_q(\alpha_0 + \zr^{j} \beta_0)$ for all $j \neq j_0$.
\end{remark}

\begin{thm}\label{thm: ClusterCminusQ}
    Let $q$ be an odd prime of bad reduction for $\Cminus/\Q$.
    \begin{enumerate}
        \item If $q \neq r$ and $q \mid ab$, then the cluster picture of $\Cminus / \Q_q$ is
        \begin{center}
            \clusterpicture            
          \Root[D] {2} {first} {r2};
          \Root[D] {4} {r2} {r3};
          \ClusterLDName c1[][n][] = (r2)(r3);
          \Root[D] {2} {c1} {r4};
          \Root[D] {4} {r4} {r5};
          \ClusterLDName c2[][n][] = (r4)(r5);
          \Root[E] {6} {c2} {r6};
          \Root[E] {2} {r6} {r7};
          \Root[E] {2} {r7} {r8};
          \Root[D] {7} {r8} {r9};
          \Root[D] {4} {r9} {r10};
          \ClusterLDName c3[][n][] = (r9)(r10);
          \Root[D] {4} {c3} {r11};
           \ClustercLDName c5[][][] = (r11);
          \ClusterLD c4[][\, 0] = (c1)(c1n)(c2)(c2n)(r6)(r7)(r8)(c3)(c3n)(c5)(c5n);
        \endclusterpicture 
        \end{center}
        where $n = \frac{p}{2} \, v_q(ab)$. Using the notation $i_0$ introduced in Definition~\ref{def: i0 j0 ppr}, each twin consists of the pair of roots $\{\gamma_{k}, \gamma_{2i_0-k}\}$, for $0\leq k \leq \frac{r-1}{2}$, $k \not \equiv i_0 \mod r$.
        \vspace{1em}
        
        \item If $q=r$ and $r \nmid ab$ then the cluster picture of $\Cminus / \Q_r$ is 
        \begin{center}
        \clusterpicture            
            \Root[D] {2} {first} {r2};
            \Root[D] {4} {r2} {r3};
            \Root[D] {4} {r3} {r4};
            \Root[E] {4} {r4} {r5};
            \Root[E] {1} {r5} {r6};
            \Root[E] {1} {r6} {r7};
            \Root[D] {4} {r7} {r8};
            \Root[D] {4} {r8} {r9};
            \Root[D] {4} {r9} {r10};
            \ClusterLD c1[][\, \dfrac{1}{r-1}] = (r2)(r3)(r4)(r5)(r6)(r7)(r8)(r9)(r10);
            \endclusterpicture.
        \end{center}
        \vspace{1em}
        
        \item  If $q=r$ and $r \mid ab$ then the cluster picture of $\Cminus / \Q_r$ is 
        \begin{center}
        \clusterpicture          
          \Root[D] {2} {first} {r2};
          \Root[D] {4} {r2} {r3};
          \ClusterLDName c1[][m][\gamma_{1}, \gamma_{-1}] = (r2)(r3);
          \Root[D] {2} {c1} {r4};
          \Root[D] {4} {r4} {r5};
          \ClusterLDName c2[][m][\gamma_2, \gamma_{-2}] = (r4)(r5);
          \Root[E] {6} {c2} {r6};
          \Root[E] {2} {r6} {r7};
          \Root[E] {2} {r7} {r8};
          \Root[D] {7} {r8} {r9};
          \Root[D] {4} {r9} {r10};
          \ClusterLDName c3[][m][\gamma_{\frac{r-1}{2}}, \gamma_{\frac{r+1}{2}}] = (r9)(r10);
          \Root[D] {4} {c3} {r11};
           \ClustercLDName c5[][][\gamma_0] = (r11);
          \ClusterLD c4[][\, \dfrac{2}{r-1}] = (c1)(c1n)(c2)(c2n)(r6)(r7)(r8)(c3)(c3n)(c5)(c5n);
        \endclusterpicture 
        \end{center}
        where $m = \frac{p}{2} \, v_r(ab) - \frac{r}{r-1}$.
    \end{enumerate}
\end{thm}

\begin{proof}    
% $(1)$ Suppose that $q \neq r$ and $q \mid ab$.
% \mar{Suppose that $q \neq r$ and $q \mid ab$. Note that $v_{q}(\alpha_0^r - \beta_0^r) = v_{q}(4\sqrt{-a^p b^p})>0$. Furthermore $v_{q}(\alpha_0)v_{q}(\beta_0)=0$ by Lemma~\ref{lem: alpha0 beta0 coprime}. Since $q \neq r$, Lemma \ref{lem:generalproduct} implies that $v_{q}(\phi_r(\alpha_0, - \beta_0)) = 0$. Therefore the product $\prod_{j \neq 0} (\gamma_0 - \gamma_j)$ has $q$-adic valuation $0$, and moreover $\gamma_0$ is an isolated root.} 
% The conventions made in Remark~\ref{rmk:conventionppr} yield $ v_q(\gamma_0 - \gamma_j) = 0$ for any $1 \leq j \leq r-1$. Thus $\gamma_0$ is an isolated root and does not belong to any cluster other than $\Rroots$. 
% For any $1 \leq k \leq r-1$, Lemma \ref{lem:generalproduct} combined with the convention made in Remark~\ref{rmk:conventionppr}implies that $v_{q} \left(\alpha_0 - \zr^{-k-j} \beta \right) > 0$ if and only if $-k-j \equiv 0 \bmod\:r$, and then by Corollary~\ref{cor: Dif 2 Roots ppr}, we have
% 	\begin{equation*}
% 		v_q (\gamma_k - \gamma_{-k}) = v_q \left(\prod_{j \neq k} (\gamma_k - \gamma_{j} ) \right) = v_q \left( \frac{4r\, \sqrt{-a^p b^p}}{\zr^k (\alpha_0 - \zr^{-2k} \beta_0)} \right) = \tfrac{p}{2} \, v_q(ab).
% 	\end{equation*}
% 	Thus we obtain $\frac{r-1}{2}$ twins of depth $\frac{p}{2} \, v_q(ab)$, and the description of the cluster picture follows.

%     \bigskip
    
$(1)$ Suppose that $q \neq r$, $q \mid ab$ and let $i_0, j_0$ as in Definition~\ref{def: i0 j0 ppr}.
%From Lemma~\ref{lem: useful eqs}, we know that $v_q(\alpha_0) v_q(\beta_0) = 0$. Lemma~\ref{lem:generalproduct}, applied with $K = \Q_q$, $\alpha = \alpha_0$, $\beta = -\beta_0$, states that there exists some $0 \leq j_0 \leq r-1$ such that $v_q(\alpha_0 - \zr^{j_0} \beta_0) = v_q(\alpha_0^r - \beta_0^r)$ and $v_q(\alpha_0 - \zr^k \beta_0) = 0$ for $k \neq j_0$. 
Since $q \neq r$, Lemma~\ref{cor: Dif 2 Roots ppr} implies that, for any $0 \leq j, k \leq r-1$ such that $j \neq k$, we have $v_q(\gamma_k - \gamma_j) > 0$ if and only if $j \equiv -k - j_0 \mod r$. It follows that
\begin{equation}\label{eq:twins2}
        v_q (\gamma_{k} - \gamma_{-k-j_0}) = v_q(\alpha_0 - \zr^{j_0} \beta_0)=v_q(\alpha_0^r - \beta_0^r) = \tfrac{p}{2} \, v_q(ab),
    \end{equation}
    where the last equality follows by Lemma~\ref{lem: useful eqs} (4). In particular, 
    %letting $i_0$ be such that $-2 i_0 \equiv j_0 \mod r$, 
    any choice of $k \neq i_0$ satisfies $- k - i_0 \not \equiv j_0 \mod r$, implying that $v_q(\gamma_k - \gamma_{i_0}) = 0$. We deduce that $\gamma_{i_0}$ is an isolated root in the cluster picture. Moreover, using the same argument as in the proof of Theorem~\ref{thm: Clusters rrp Q} (1) and \eqref{eq:twins2},
     %Moreover, for any $1 \leq k \leq \frac{r-1}{2}$, we have $k - (2i_0 - k) \equiv j_0 \mod r$.
     %so $v_q(\gamma_{k} - \gamma_{2i_0 - k}) > 0$. 
     %We can explicitly compute this valuation:
    %The conventions made in Remark~\ref{rmk:conventionrrp} yield $ v_q(\gamma_0 - \gamma_j) = 0$ for any $1 \leq j \leq r-1$. Thus $\gamma_0$ is an isolated root and does not belong to any cluster other than $\Rroots$. 
    %Similarly, for any $1 \leq k \leq r-1$, Lemma \ref{lem:generalproduct} combined with the convention made in Remark~\ref{rmk:conventionrrp} asserts that $v_q( a + \zr^{-k-j} b) > 0$ if and only if $- k  -j \equiv 0 \bmod \, r$. 
    % \begin{equation*}
    %     v_q (\gamma_{k} - \gamma_{2i_0 - k}) = v_q \left( \prod_{\substack{0 \leq j \leq r-1 \\ j \not \equiv i_0 + k \mod r}} (\gamma_{i_0 + k} - \gamma_j) \right) = v_q \left( \frac{r (a^r + b^r)}{\zr^k(a + \zr^{-2k} b)} \right) = v_q(a^r + b^r).
    % \end{equation*}
    there are $\frac{r-1}{2}$ twins of depth $\frac{p}{2}v_q(ab)$, each of them consisting of the roots $\lbrace \gamma_{k}, \gamma_{2i_0 - k} \rbrace$ where $k\neq i_0$.
  
	$(2)$ Suppose $q =r$ and $r \nmid ab$. Then $v_r \left(\alpha_0^r - \beta_0^r \right) = v_r ( 4 i \sqrt{a}^p \sqrt{b}^p ) = 0$ by Lemma~\ref{lem: useful eqs}, so for any $0 \leq j, k \leq r-1$ with $j \neq k$, we have $v_r (\alpha_0 - \zr^{-k-j} \beta_0) = 0$, and thus $v_r(\gamma_k - \gamma_j) = v_r( 1 - \zr^{j-k}) = \frac{1}{r-1}$ from Lemma~\ref{cor: Dif 2 Roots ppr}.
	
	$(3)$ Suppose $q = r$, and $r \mid ab$. We use Remark~\ref{rmk: i0 j0 ppr} to compute 
    %then $v_r( \alpha_0^r - \beta_0^r) = \frac{p}{2} \, v_r(ab) > \frac{3}{2}$ and $r v_r(1 - \zr) = \frac{r}{r-1} \leq \frac{3}{2}$. Note that $v_{r}(\alpha_0)v_{r}(\beta_0)=0$ by Lemma~\ref{lem: alpha0 beta0 coprime}. Thus we can apply Lemma \ref{lem:generalproduct} to $\alpha_0^r - \beta_0^r$. 
    %For any $1 \leq j \leq r-1$, we have $v_r(\alpha_0 - \zr^{j} \beta_0) = v_r(1 - \zr) = \frac{1}{r-1}$. We also deduce that 
    \begin{equation*}
		v_{r}(\alpha_0 - \beta_0) = v_{r}(\alpha_0^r - \beta_0^r) - \sum_{j=1}^{r-1} v_r(\alpha_0 - \zr^{j} \beta_0) = \tfrac{p}{2} \, v_r(ab) - 1.
	\end{equation*}
    Therefore, when $k = 0$, we get $v_r(\gamma_0 - \gamma_j) = v_r((1 - \zr^{j})(\alpha_0 - \zr^{-j} \beta_0)) = 2 v_r(1- \zr) = \frac{2}{r-1}$ for $j \neq k$. If $1 \leq k \leq r-1$ then, for any $j \neq k$, Lemma~\ref{cor: Dif 2 Roots ppr} yields
	\begin{equation*}
		v_r(\gamma_k - \gamma_j) = v_r((1 - \zr^{j-k}) (\alpha_0 - \zr^{-j-k} \beta_0)) = \begin{cases}
			\frac{2}{r-1} & \text{ if } j \not \equiv -k \bmod\:r, \\
			\frac{1}{r-1} + \frac{p}{2} \, v_r(ab)  -1 & \text{ if } j \equiv -k \bmod\:r.
		\end{cases}
	\end{equation*}
	Therefore all the roots belong to a cluster having depth $\frac{2}{r-1}$, the root $\gamma_0$ is isolated, and there are $\frac{r-1}{2}$ twins, each with relative depth $\frac{p}{2} \, v_r(ab) - \frac{r}{r-1}$.
\end{proof}

\begin{cor}\label{cor: ClusterCminusQpl}
Let $\mathfrakq$ be an odd place of $\Qzrplus$ lying above the rational prime $q$. The cluster picture of $\Cminus / \Qzrplus_{\mathfrakq}$ is obtained from the cluster picture of $\Cminus / \Q_q$ as in Theorem~\ref{thm: ClusterCminusQ}, by multiplying all depths by the index of ramification $e_{\Qzrplus_{\mathfrakq} / \Q_q} = \begin{cases}
1 & \text{ if } q \neq r, \\
\frac{r-1}{2} & \text{ if } q = r.
\end{cases}$
\end{cor}

\begin{proof}
This is an immediate consequence of Theorem~\ref{thm: ClusterCminusQ} and \eqref{eq:valsrelation}.  
\end{proof}

\begin{remark} \label{rmk: el lema de la noche}
    When $q \mid a$, Lemma~\ref{lem: useful eqs} \textit{(5)} implies that $v_q(\alpha_0^r + c^r) > 0$. Note that if $v_q(a)>0$ then $v_q(c)=0$, otherwise this would contradict the coprimality of $a$, $b$ and $c$. In this case, Lemma~\ref{lem: useful eqs} \textit{(5)} also implies that $v_q(\alpha_0)=0$ by the same reasoning. If $q\neq r$, then Lemma~\ref{lem:generalproduct} applied with $\alpha = \alpha_0$, $\beta = c$ states that there is some $0 \leq k_0 \leq r-1$ such that $v_q(\alpha_0 + \zr^{k_0} c) > 0$ and $v_q(\alpha_0 + \zr^j c) = 0$ for $j \neq k_0$. Recall from Remark~\ref{rmk: i0 j0 ppr} that $v_q(\alpha_{0} - \zr^{j_{0}} \beta_{0}) >0$. Now equality
    \[
    \left(\alpha_0 + \zr^{k_0} c \right) - \left(\alpha_0 - \zr^{j_0} \beta_0 \right) = \frac{\zr^{k_0} c}{\alpha_0} \left( \alpha_0 + \zr^{j_0 - k_0} c \right)
    \]
    implies $v_q(\alpha_0 + \zr^{j_0 - k_0} c)>0$. The definition of $i_0$ implies $k_0 \equiv -i_0 \mod r$, and so $v_q(\alpha_0 + \zr^{-i_0} c) > 0$. 
    We conclude that $v_q(\alpha_0 + \zr^{k_0} c) = v_q(\alpha_0 + \zr^{-i_0} c) = pv_q(a)/2$ by Lemma~\ref{lem: useful eqs} \textit{(5)}.
    Similarly when $q=r$, by Lemma~\ref{lem:generalproduct} applied with $\alpha=\alpha_0$ and $\beta=c$, we have $v_r(\alpha_0 + c) > \frac{1}{r-1}$ and $v_r(\alpha_0 + \zr^j c) = \frac{1}{r-1}$ for $1 \leq j \leq r-1$.
    
\end{remark}

We now compute the cluster pictures of the hyperelliptic curve $\Cplus / \Q$ and $\Cplus / \Qplus$. 

\begin{thm}\label{thm: ClusterCplusQ}
	 Let $q$ be an odd prime number of bad reduction for $\Cplus/\Q$. 
 
	\begin{enumerate}
		\item If $q \neq r$ and $q \mid a$, then the cluster picture of $\Cplus / \Q_q$ is
		\begin{center}
			\clusterpicture      
			\Root[D] {2} {first} {r2};
			\Root[D] {4} {r2} {r3};
			\ClusterLDName c1[][n][] = (r2)(r3);
			\Root[D] {2} {c1} {r4};
			\Root[D] {4} {r4} {r5};
			\ClusterLDName c2[][n][] = (r4)(r5);
			\Root[E] {6} {c2} {r6};
			\Root[E] {2} {r6} {r7};
			\Root[E] {2} {r7} {r8};
			\Root[D] {7} {r8} {r9};
			\Root[D] {4} {r9} {r10};
			\ClusterLDName c3[][n][] = (r9)(r10);
			\Root[D] {4} {c3} {r11};
			\Root[D] {4} {r11} {r12};
			\ClusterLDName c4[][2n][\lastroot \, , \gamma_{i_0}] = (r11)(r12);
			\ClusterLD c4[][\, 0] = (c1)(c1n)(c2)(c2n)(r6)(r7)(r8)(c3)(c3n)(c4)(c4n);
			\endclusterpicture
		\end{center}
  where $n = \frac{p}{2} \, v_q(a)$. Using the notation $i_0$ introduced in Definition~\ref{def: i0 j0 ppr}, each twin of depth $n$  consists of the pair of roots $\{\gamma_{k}, \gamma_{2i_0-k}\}$, for $0\leq k \leq \frac{r-1}{2}$, $k \not \equiv i_0 \mod r$.
		\vspace{1em}
		
		\item If $q \neq r$ and $q \mid b$, then the cluster picture of $\Cplus / \Q_q$ is
		\begin{center}
                \clusterpicture 
			\Root[D] {2} {first} {r2};
			\Root[D] {4} {r2} {r3};
			\ClusterLDName c1[][n][] = (r2)(r3);
			\Root[D] {2} {c1} {r4};
			\Root[D] {4} {r4} {r5};
			\ClusterLDName c2[][n][] = (r4)(r5);
			\Root[E] {6} {c2} {r6};
			\Root[E] {2} {r6} {r7};
			\Root[E] {2} {r7} {r8};
			\Root[D] {7} {r8} {r9};
			\Root[D] {4} {r9} {r10};
			\ClusterLDName c3[][n][] = (r9)(r10);
			\Root[D] {4} {c3} {r11};
			\Root[D] {4} {r11} {r12};
            \ClustercLDName c5[][][\lastroot] = (r11);
            \ClustercLDName c6[][][\gamma_{i_0}] = (r12);
			\ClusterLD c4[][\, 0] = (c1)(c1n)(c2)(c2n)(r6)(r7)(r8)(c3)(c3n)(c5)(c5n)(c6);
			\endclusterpicture
		\end{center}
   where $n = \frac{p}{2} \, v_q(b)$. Denoting by $\gamma_{i_0}$ the rightmost root, each twin consists of the pair of roots $\{\gamma_{k}, \gamma_{2i_0-k}\}$, for $0\leq k \leq \frac{r-1}{2}$, $k \not \equiv i_0 \mod r$.
		\vspace{1em}
		
		\item If $q = r$ and $r \nmid ab$, then the cluster picture of $\Cplus / \Q_r$ is 
		\begin{center}
			\clusterpicture            
			\Root[D] {2} {first} {r2};
			\Root[D] {4} {r2} {r3};
			\Root[D] {4} {r3} {r4};
			\Root[E] {4} {r4} {r5};
			\Root[E] {1} {r5} {r6};
			\Root[E] {1} {r6} {r7};
			\Root[D] {4} {r7} {r8};
			\Root[D] {4} {r8} {r9};
			\Root[D] {4} {r9} {r10};
			\ClusterLD c1[][\frac{1}{r-1}] = (r2)(r3)(r4)(r5)(r6)(r7)(r8)(r9)(r10);
			\Root[D] {2} {c1} {r11};
             \ClustercLDName c5[][][\lastroot] = (r11);
			\ClusterLD c4[][\, 0 \ ] = (c1)(c1n)(c2)(c2n)(r6)(r7)(r8)(c3)(c3n)(c5)(c5n);
			\endclusterpicture.
		\end{center}
 
		\vspace{1em}
		
		\item If $q = r$ and $r \mid a$, then the cluster picture of $\Cplus / \Q_r$ is
		\begin{center}
			\clusterpicture            
			\Root[D] {2} {first} {r2};
			\Root[D] {4} {r2} {r3};
			\ClusterLDName c1[][m][\gamma_1, \gamma_{-1}] = (r2)(r3);
			\Root[D] {2} {c1} {r4};
			\Root[D] {4} {r4} {r5};
			\ClusterLDName c2[][m][\gamma_2, \gamma_{-2}] = (r4)(r5);
			\Root[E] {6} {c2} {r6};
			\Root[E] {2} {r6} {r7};
			\Root[E] {2} {r7} {r8};
			\Root[D] {7} {r8} {r9};
			\Root[D] {4} {r9} {r10};
			\ClusterLDName c3[][m][\gamma_{\frac{r-1}{2}}, \gamma_{\frac{r+1}{2}}] = (r9)(r10);
			\Root[D] {4} {c3} {r11};
			\Root[D] {4} {r11} {r12};
			\ClusterLDName c4[][2m][\lastroot \, , \gamma_0] = (r11)(r12);
			\ClusterLD c4[][\, \dfrac{2}{r-1}] = (c1)(c1n)(c2)(c2n)(r6)(r7)(r8)(c3)(c3n)(c4)(c4n);
			\endclusterpicture
		\end{center}
  where $m = \frac{p}{2} \, v_r(a) - \frac{r}{r-1}$.
		\vspace{1em}
	
		\item If $q = r$ and $r \mid b$, then the cluster picture of $\Cplus / \Q_r$ is
		\begin{center}
			\clusterpicture            
			\Root[D] {2} {first} {r2};
			\Root[D] {4} {r2} {r3};
			\ClusterLDName c1[][m][\gamma_1, \gamma_{-1}] = (r2)(r3);
			\Root[D] {2} {c1} {r4};
			\Root[D] {4} {r4} {r5};
			\ClusterLDName c2[][m][\gamma_2, \gamma_{-2}] = (r4)(r5);
			\Root[E] {6} {c2} {r6};
			\Root[E] {2} {r6} {r7};
			\Root[E] {2} {r7} {r8};
			\Root[D] {7} {r8} {r9};
			\Root[D] {4} {r9} {r10};
			\ClusterLDName c3[][m][\gamma_{\frac{r-1}{2}}, \gamma_{\frac{r+1}{2}}] = (r9)(r10);
			\Root[D] {4} {c3} {r11};
    \ClustercLDName c7[][][\gamma_0] = (r11);
			\ClusterLD c4[][\frac{2}{r-1}] = (c1)(c1n)(c2)(c2n)(r6)(r7)(r8)(c3)(c3n)(c7)(c7n);
			\Root[D] {2} {c4} {r12};
    \ClustercLDName c6[][][\lastroot] = (r12);
			\ClusterLD c5[][\, 0] = (c1)(c1n)(c2)(c2n)(r6)(r7)(r8)(c3)(c3n)(c4)(c4n)(c6)(c6n);
			\endclusterpicture
		\end{center}
    where $m = \frac{p}{2} \, v_r(b) - \frac{r}{r-1}$.
	\end{enumerate}
\end{thm}

\begin{proof}
	By Theorem~\ref{thm: ClusterCminusQ}, it remains to compute the $q$-adic valuation of $\lastroot - \gamma_j$ for $0 \leq j \leq r-1$. 
    %We know from Lemma~\ref{lem: useful eqs} that $v_q(\alpha_0) v_q(\beta_0) = 0$. Thus $v_q(\sqrt{\alpha_0}) v_q(\sqrt{\beta_0}) = 0$.
 
    $(1)$ 
    %Note that $v_{q}(\alpha_0)v_{q}(\beta_0)=0$ by Lemma~\ref{lem: alpha0 beta0 coprime}. 
    If $q \mid a$ then Remark~\ref{rmk: el lema de la noche} with Proposition~\ref{prop extra root} give that $v_q(\gamma_{i_0}-\lastroot)=p \, v_q(a)>0$ and $v_q(\gamma_{k}-\lastroot)=0$ for all $0 \leq k \leq r-1$, $k \neq i_0$.
    We conclude that $\gamma_{i_0}$ and $\lastroot$ lie in a twin of depth $p \, v_q(a)$ and that the outer depth remains unchanged.
	
    $(2)$ If $q \mid b$ then $v_q \left( \prod_{j= 0}^{r-1} (\gamma_j - \lastroot)\right) = v_q(4a^p) = 0$ by Proposition~\ref{prop extra root} and the assumption $\gcd(a,b)=1$. 
    Thus $v_q(\gamma_j - \lastroot) = 0$ for every $0 \leq j \leq r-1$, and moreover $\lastroot$ does not belong to any cluster other than $\Rroots$.
	
    $(3)$ If $r \nmid ab$ then $v_r \left( \prod_{j= 0}^{r-1} (\gamma_j - \lastroot)\right) = 0$ by Proposition~\ref{prop extra root}. 
    Thus $\lastroot$ lies outside of the cluster containing the other roots.
	
    $(4)$ %Note that $v_{r}(\alpha_0)v_{r}(\beta_0)=0$ by Lemma~\ref{lem: useful eqs}. 
    If $r \mid a$ then Remark~\ref{rmk: el lema de la noche} combined with Lemma~\ref{prop extra root} gives 
    $ v_r(\gamma_{j} - \lastroot) = \frac{2}{r-1}$
    for all $1 \leq j \leq r-1$. Again, from Lemma~\ref{prop extra root} we get \vspace{-0.5em}
    \begin{equation*}
		v_r(\gamma_0-\lastroot)= 2(v_r \left(\alpha_0 + c \right))  = 2(v_r \left( \alpha_0^r + c^r \right) - \sum_{j= 1}^{r-1} v_r( \alpha_0 + \zr^{j} c)) = p \, v_r(a) - 2. \vspace{-0.25em}
	\end{equation*}
    Thus $\gamma_0$ and $\lastroot$ lie in a twin of relative depth $pv_r(a) -\frac{2r}{r-1}$.
    
    $(5)$ If $r \mid b$ then $v_r \left( \prod_{j= 0}^{r-1} (\gamma_j - \lastroot)\right) = 0$ by Proposition~\ref{prop extra root}, and moreover $\lastroot$ is an isolated root.
\end{proof}

\begin{cor}\label{cor: ClusterCplusQpl}
Let $\mathfrakq$ be an odd place of $\Qzrplus$ lying above the rational prime $q$. The cluster picture of $\Cplus / \Qzrplus_{\mathfrakq}$ is obtained from the cluster picture of $\Cplus / \Q_q$ as in Theorem~\ref{thm: ClusterCplusQ}, by multiplying all depths by the index of ramification $e_{\Qzrplus_{\mathfrakq} / \Q_q} = \begin{cases}
1 & \text{ if } q \neq r, \\
\frac{r-1}{2} & \text{ if } q = r. \end{cases}$
\end{cor}

\begin{proof}
 This is an immediate consequence of Theorem~\ref{thm: ClusterCplusQ} and \eqref{eq:valsrelation}. 
\end{proof}

%%%%%%%%%%%%%%%%%%%%%%%%%%%%%%%%%%%%%%%%%%%%%%%%%%%%%%%%%%%%%%%%%%%%%%
\subsection{Ramification index of the splitting field} \label{sec: Ramification indices}
%%%%%%%%%%%%%%%%%%%%%%%%%%%%%%%%%%%%%%%%%%%%%%%%%%%%%%%%%%%%%%%%%%%%%%%%%%%%%
Recall that, to use Theorem \ref{thm: allcond} to compute conductor exponents of $\Cminus$ and $\Cplus$, we need their cluster pictures together with the action of Galois on the roots. In this section, we give the necessary tools for the latter. 
%In order to better understand this, 
We begin by discussing the reducibility of $\gminus$, and then describe its splitting field, and its ramification index over $\Q_q$. Moreover, we study the discriminants of certain totally ramified extensions of local fields; this will be essential to compute the wild part of the conductor exponents in Section \ref{sec: Conductors ppr}.  
\medskip

%Recall the definition of $i_0$ given in Definition~\ref{def: i0 j0 ppr}. Recall that when $q = r$, we defined $i_0 = 0$.

%\begin{prop}\label{prop: elephant i0}
%    If $q \mid ab$, then $\gminus$ is reducible and $\gamma_{i_0} \in \Q_q$. 
%\end{prop}

%\begin{proof}
 %   The absolute Galois group of $\Q_q$ acts on the clusters, preserving depths and containments (\cite[Remark 3.2]{hyperusersguide}). When $q \mid ab$, Theorem~\ref{thm: ClusterCminusQ} implies that $\gamma_{i_0}$ is an isolated root in the cluster of $\Cminus / \Q_q$, so it is fixed by all the elements of $G_{\Q_q}$, hence the result.
%\end{proof}

\begin{defn} We let $\sfL \coloneqq \Q_q(\Rroots)$ be the splitting field of the polynomial $\gminus\in \Q_{q}[x]$. 
\end{defn}

\begin{remark}
    Recall that the set roots of $\gplus$ is $\Rroots^{+} = \Rroots \cup \lbrace \lastroot \rbrace$, with $\lastroot = -2c \in \Q_q$. Therefore, $L$ is also the splitting field of $\gplus$.
\end{remark}

\begin{lemma}\label{lem: alpha beta 2k}
	For any integer $n \geq 0$, we have $\alpha_{0}^{2n}+\beta_{0}^{2n} \in \Q_q(\gamma_0)$. 
\end{lemma}

\begin{proof}
    The polynomial $X^{2n} + Y^{2n}$ is symmetric, so there exists some $P \in \Q_q[X, Y]$ such that $X^{2n} + Y^{2n} = P(X+Y, \, XY)$. Thus from Lemma~\ref{lem: useful eqs}, $\alpha_0^{2n} + \beta_0^{2n} = P(\alpha_0 + \beta_0, \alpha_0 \beta_0) = P(\gamma_0, c^2) \in \Q_q(\gamma_0)$.
    %We prove the result by induction on $n$. The case $n = 0$ holds trivially. Therefore, we suppose that $n \geq 1$. Assuming that the result holds for any $0 \leq j \leq n-1$, we show that it also holds for $n$. By assumption, $\alpha_0^2 + \beta_0^2 \in \Q_r(\gamma_0)$, so taking its $n$-th power and reordering terms gives
	%\begin{align*}
		%(\alpha_{0}^{2}+\beta_{0}^{2})^{n} & = \sum_{j=0}^{n} \binom{n}{j} \,  \alpha_{0}^{2j}\beta_{0}^{2(n-j)} \\
	%	& = \sum_{k = 0}^{\left \lfloor (n-1) / 2 \right \rfloor} \binom{n}{k} \, \alpha_{0}^{2k} \beta_{0}^{2k} \left( \alpha_{0}^{2(n-2k)} + \beta_{0}^{2(n-2k)} \right) + \begin{cases}
	%		\binom{n}{n/2} \, \alpha_{0}^{n} \beta_{0}^{n} \text{ if } n \text{ is even,} \\
	%		0 \text{ otherwise.}
	%	\end{cases}
	%\end{align*}
	%As stated above, the left-hand side belongs to $\Q_q(\gamma_0)$. Combining the fact that $\alpha_0 \beta_0 = c^2 \in \Q_q$ with the induction hypothesis shows that $\alpha_0^{2n} + \beta_0^{2n}$ can be written as the difference of two elements of $\Q_q(\gamma_0)$; this completes the result. 
\end{proof}

We now describe the splitting field $L$. For any $1 \leq j \leq r-1$, recall that we defined $\omega_j = \zr^j + \zr^{-j}$ and $\omega= \omega_1$. Moreover, we let $\tau_j \coloneqq \zr^j - \zr^{-j} \in \overline{\Q_q}$ and $\tau \coloneqq \tau_1$. 

\begin{theorem}\label{thm:splittingfield} 
	The splitting field $L$ is given by $L=\Q_{q}(\omega, \tau\sqrt{-ab}, \gamma_{0})$.
\end{theorem}

\begin{proof}
	Let $L' \coloneqq \Q_{q}(\omega, \tau\sqrt{-ab}, \gamma_{0})$. We begin by proving that $L' \subseteq \sfL$. First of all, it is clear that $\gamma_0 \in \sfL$, and the equality $\omega=(\gamma_1+\gamma_{-1})/ \gamma_0$ shows that $\omega\in \sfL$. In order to prove that $\tau \sqrt{-ab} \in \sfL$, note that $\tau_j (\alpha_{0}-\beta_0)=\gamma_{j} -\gamma_{-j} \in L$. On the other hand, we have $(\alpha_0^r - \beta_0^r)/ \sqrt{-ab} = 4 (-ab)^{(p-1)/2}\in \Q_{q}$. Finally by the definition of $\phi_r$ (Definition~\ref{def: phir}) and the equality $\alpha_0\beta_0=c^2$, we have
	\begin{align*}
		\phi_r(\alpha_0, -\beta_0) = \sum_{j = 0}^{r-1} \alpha_0^j \beta_0^{r-1-j} & = \alpha_0^{(r-1)/2} \beta_0^{(r-1)/2} + \sum_{j=0}^{(r-3)/2} (\alpha_0^j\beta_0^{r-1-j}+\alpha_0^{r-1-j}\beta_0^{j}) \\
		& = c^{r-1} + \sum_{j=0}^{(r-3)/2} c^{2j}(\beta_0^{r-1-2j}+\alpha_0^{r-1-2j}).
	\end{align*}
	By Lemma~\ref{lem: alpha beta 2k}, every term in the sum belongs to $\Q_q(\gamma_0)$, so $\phi_r(\alpha_0, - \beta_0) \in \Q_q(\gamma_0) \subset \sfL$. 
 Then the equality
	\begin{equation}\label{eq: tau sqrt -ab}
		\tau_j \sqrt{-ab} = \tau_j (\alpha_0 - \beta_0) \cdot \frac{\sqrt{-ab}}{\alpha_0^r - \beta_0^r} \cdot \phi_r(\alpha_0, - \beta_0)
	\end{equation}
    implies that $\tau_j \sqrt{-ab} \in \sfL$, for any $1 \leq j \leq r-1$.
	
	We now prove the reverse inclusion $L \subseteq L'$. Fix $1 \leq j \leq r-1$. It is easy to check that
	\begin{equation*}
    \gamma_j + \gamma_{-j} = \omega_j \gamma_0 \quad \text{ and } \quad \gamma_j - \gamma_{-j} = \tau_j (\alpha_0 - \beta_0), \quad \text{ so } \quad  \gamma_j = \frac{1}{2} \left(\omega_j \gamma_0 + \tau_j (\alpha_0 - \beta_0) \right).
	\end{equation*}
    To prove that $\gamma_j \in L'$, it suffices to check that $\tau_j (\alpha_0 - \beta_0) \in L'$, since $\omega_j \in \Q_q(\omega) \subset L'$ and $\gamma_0 \in L'$. By Galois theory, $\Q_q(\omega)$ is the subfield of $\Q_q(\zr)$ fixed by the automorphism $\sigma \in \Gal(\Q_q(\zr) / \Q_q)$ described by $\sigma : \zr \mapsto \zr^{-1}$. Note that $\tau / \tau_j \in \Q_{q}(\omega)$ since this element is fixed by $\sigma$. Therefore $\tau_j \sqrt{-ab} = \frac{\tau_j}{\tau} \cdot \tau \sqrt{-ab}\in \Q_q(\omega, \tau \sqrt{-ab}) \subset L'$.
    We conclude that $\tau_j (\alpha_0 - \beta_0) \in L'$ by \eqref{eq: tau sqrt -ab}. 
    %\diana{17/09/2025: Do we need to consider separately the case where $\zeta_r \in \Q_q$ and the Galois argument doesn't work?}
\end{proof}

\begin{lemma}\label{lem: gamma0 quadratic}
    If $\gminus$ is reducible, then $\gamma_0 \in \Q_q(\tau\sqrt{-ab})$. Moreover, if $\zr \notin \Q_q$, then $\gamma_0 \in \Q_q$.
\end{lemma}

\begin{proof}
    Consider the following towers of extensions
    \begin{equation} \label{eq:towerestensions}
	\Q_{q}\overset{d_0\mid 2}{\subseteq} \Q_{q}(\sqrt{-ab})\overset{n_{0} \mid r}{\subseteq} \Q_{q}(\alpha_{0})\qquad \text{ and } \qquad 
	\Q_{q}\overset{n_{0}'\leq r}{\subseteq} \Q_{q}(\gamma_{0})\overset{d_{0}'\mid 2}{\subseteq} \Q_{q}( \alpha_0),
	\end{equation}
	where the superscripts denote the degree of the respective extension. If $\gminus$ is reducible, then $n_0' < r$, and then $n_0' d_0' = n_0 d_0 < 2r$. By Remark~\ref{rmk: alpha0 in Qq}, $n_0$ is either equal to $1$ or $r$. We show that $n_0=1$. If not then $n_0=r$, and hence $d_0=1$. Since $d_0'\mid 2$ and $n_0'd_0'=n_0d_0$, we obtain $n_0'=r$. This contradicts the assumption that $g_{r}^{-}$ is reducible.
    
    Therefore, $n_0 = 1$. In this case, $\alpha_0\in \Q_q(\sqrt{-ab})$. Recall that $\beta_0 = c^2 / \alpha_0$. If $\zr \in \Q_q$, then $\Q_q(\tau \sqrt{-ab}) = \Q_q(\sqrt{-ab})$, and so $\gamma_0 = \alpha_0 + c^2 / \alpha_0$ belongs to $\Q_q(\tau \sqrt{-ab})$. 
    Now suppose $\zr \notin \Q_q$. If $d_0=1$ then $\alpha_0\in \Q_q$ from which it immediately follows that $\gamma_0=\alpha_0+c^2/\alpha_0\in\Q_q$. If $d_0=2$, let $\sigma$ (resp. $\Norm$, $\Trace$) be the non-trivial automorphism (resp. the norm and trace maps) of the extension $\Q_q(\sqrt{-ab}) / \Q_q$. The definition of $\alpha_0$ combined with Lemma~\ref{lem: useful eqs} gives $\sigma(\alpha_0)^r = \beta_0^r = c^{2r} / \alpha_0^r$, so $\Norm(\alpha_0)^r = c^{2r}$. But $\zr \notin \Q_q$, so we deduce that $\Norm(\alpha_0) = c^2$, so $\sigma(\alpha_0) = \beta_0$, and $\gamma_0 = \Trace(\alpha_0) \in \Q_q$. 
    \begin{comment}
    If $\zr \notin \Q_q$, let $\sigma$ (resp. $\Norm$, $\Trace$) be the non-trivial automorphism (resp. the norm and trace maps) of the extension $\Q_q(\sqrt{-ab}) / \Q_q$. The definition of $\alpha_0$ combined with Lemma~\ref{lem: useful eqs} gives $\sigma(\alpha_0)^r = \beta_0^r = c^{2r} / \alpha_0^r$, so $\Norm(\alpha_0)^r = c^{2r}$. But $\zr \notin \Q_q$, so we deduce that $\Norm(\alpha_0) = c^2$, so $\sigma(\alpha_0) = \beta_0$, and $\gamma_0 = \Trace(\alpha_0) \in \Q_q$. 
    \end{comment}
\end{proof}

\begin{prop}\label{prop: elephant i0}
    If $q \mid ab$, then $\gminus$ is reducible and $\gamma_{i_0} \in \Q_q$. 
\end{prop}

\begin{proof}
    The absolute Galois group of $\Q_q$ acts on the clusters, preserving depths and containments (\cite[Remark 3.2]{hyperusersguide}). When $q \mid ab$, Theorem~\ref{thm: ClusterCminusQ} implies that $\gamma_{i_0}$ is an isolated root in the cluster of $\Cminus / \Q_q$, so it is fixed by all the elements of $G_{\Q_q}$, hence the result.
\end{proof}

The purpose of the following lemma is to assist in computing the ramification index of $L/\Q_{q}$ when $q=r$, through the use of certain subextensions.

\begin{lemma}\label{lem:fielddiagram}
The diagram in Figure~\ref{fig:RamDiagram} below presents the ramification indices of the displayed extensions of $\Q_r$.
        
\end{lemma}

\begin{proof}
Let $K=\Q_{r}(\omega, \tau\sqrt{-ab})$. Recall that $L=\Q_{r}(\mathcal{R})=K(\gamma_0)=\Q_{r}(\omega, \tau\sqrt{-ab}, \gamma_0)$ by Theorem~\ref{thm:splittingfield}, and that $K(\gamma_0) = K(\alpha_0)$. Thus the inclusion of fields are shown in Figure~\ref{fig:RamDiagram} is clear from the definition of $\alpha_0, \beta_0$ and $\gamma_0$ (Definition~\ref{defn: Rootsppr}). The ramification indices within the extensions given by
\[
\Q_{r}\overset{e\, = \, r-1}{\subseteq} \Q_{r}(\zeta_r) \overset{e\, \mid\, 2}{\subseteq} \Q_{r}(\zeta_r, \alpha_0^{r})\overset{e \,\mid\, r}{\subseteq} \Q_{r}(\zeta_r, \alpha_0)
\]
and 
\[
\Q_{r}\overset{e\, \mid \, 2}{\subseteq} \Q_{r}(\sqrt{-ab})=\Q_{r}(\alpha_0^r) \overset{e\, \mid\, r-1}{\subseteq} \Q_{r}(\zeta_r, \alpha_0^{r})\overset{e \,\mid\, r}{\subseteq} \Q_{r}(\zeta_r, \alpha_0)
\]
are an immediate consequence of the definition of $\alpha_0$ (Definition~\ref{defn: Rootsppr}). 
Note that the ramification index of the extension $\Q_{r}(\zeta_r, \alpha_0)/L$ divides $2$ since $\Q_{r}(\zeta_r)/\Q_{r}(\omega)$ is a quadratic extension. 
We now compute the ramification index of the extension $K / \Q_r$. Recall that $\Q_r(\omega) / \Q_r$ is totally ramified of degree $\frac{r-1}{2}$; it thus remains to compute the ramification index of the extension $K / \Q_r(\omega)$. The defining polynomial $x^2 + ab(\omega_2 - 2)$ has discriminant $\Delta:=4ab(2 - \omega_2)$. Recall that $v_{\rQzrplus}(\omega_2 - 2) = 1$ and $v_{\rQzrplus}(ab) = \frac{r-1}{2} v_r(ab)$. 
If $K=\Q_r(\omega)$, then $v_\mathfrak{r}(\Delta)$ is even because $\Delta$ defines a trivial extension.
        Suppose that $K / \Q_r(\omega)$ is quadratic, then 
         it is ramified if and only if $v_{\rQzrplus}(\Delta)$ is odd (e.g.~\cite{Cohen}[Theorem 10.2.9]). Since $r$ is odd, we have $v_{\rQzrplus}(\Delta) = 1 + v_{\rQzrplus}(ab) = 1 + \frac{r-1}{2} v_r(ab)$, and we conclude that
	\begin{equation*}
		e_{K / \Q_r(\omega)} = \begin{cases}
			2 & \text{ if } \frac{r-1}{2} v_r(ab) \text{ is even,} \\
			1 & \text{ if } \frac{r-1}{2} v_r(ab) \text{ is odd,}
		\end{cases} \quad \text{ so } \quad e_{K / \Q_r} = \begin{cases}
			r-1 & \text{ if } \frac{r-1}{2} v_r(ab) \text{ is even,} \\
			\frac{r-1}{2} & \text{ if } \frac{r-1}{2} v_r(ab) \text{ is odd.}
		\end{cases}
    \end{equation*}

The remainder of the ramification indices indicated in Figure~\ref{fig:RamDiagram} now follow from the multiplicativity of ramification indices in field extensions.

\end{proof}
   \begin{figure}[h!]

        \adjustbox{scale=0.8,center}{
            \begin{tikzcd}
                {} & {\mathbb{Q}_{r}(\zeta_r, \alpha_0) } &  & \\
                L = K(\gamma_0) \arrow[ru, "e\; \mid\; 2", no head] & {\mathbb{Q}_{r}(\zeta_r, \alpha_0^r)} \arrow[u, "e\; \mid\; r", no head]  & {} & \\
                & & \mathbb{Q}_{r}(\alpha_0)=\mathbb{Q}_{r}(\beta_0) \arrow[swap, luu, "e=r-1", no head] &  \\
                {K\coloneqq \mathbb{Q}_{r}(\omega, \tau\sqrt{-ab})} \arrow[d, swap, "e\; \mid\; 2", no head] \arrow[ruu, "e\; \mid\; 2", no head] \arrow[uu, "e\; \mid\; r", no head] & \mathbb{Q}_{r}(\zeta_r) \arrow[uu, "e\; \mid\; 2", no head] &  &      \\
                \mathbb{Q}_{r}(\omega) \arrow[ru, "e=2", no head] & & \mathbb{Q}_{r}(\sqrt{-ab})=\mathbb{Q}_{r}(\alpha_0^r) \arrow[uu, "e\; \mid\; r", no head] \arrow[luuu, "e\; \mid\; r-1", no head] & \mathbb{Q}_{r}(\gamma_0) \arrow[luu, "e\; \mid \; 2", no head] \\
                & \mathbb{Q}_{r} \arrow[uu, "e=r-1", no head] \arrow[ru, "e\; \mid\; 2", no head] \arrow[swap, rru, "e\; \mid\; r", no head] \arrow[lu, "e=\frac{r-1}{2}", no head] & &       
        \end{tikzcd}}
        \caption{\centering Diagram of field inclusions when $q= r$. The symbol $e$ denotes the ramification index of the respective extension.}
        \label{fig:RamDiagram}
        \end{figure}
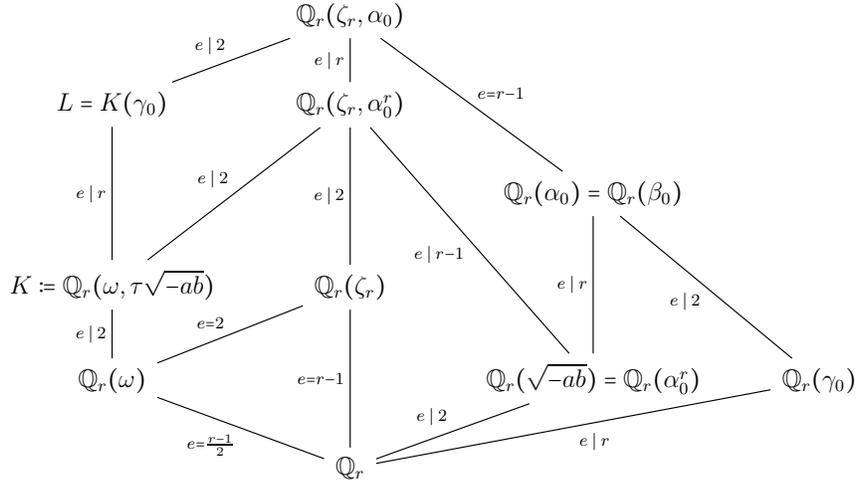
We now compute the ramification index of the extension $L / \Q_q$. This is related to the reducibility of $\gminus$.

\begin{thm} \label{thm: Ramif index split field} 
	The ramification index $e_{L / \Q_q}$ of $L / \Q_q$ is described as follows.
	\begin{enumerate}
            \item If $q \neq r$ and $q \mid ab$, then $\gminus$ is reducible and $e_{L / \Q_q} \leq 2$.

            \item If $q = r$ and $r\nmid ab$ we have 
		\begin{equation*}
			e_{L/\Q_{r}} = \begin{cases}
				(r-1) & \text{if} \ g_{r}^{-} \text{ is reducible}, \\
				r(r-1) & \text{if} \ g_{r}^{-} \text{ is irreducible}.
		\end{cases}
            \end{equation*}
            
            \item If $q = r$ and $r\mid ab$, then $g_{r}^{-}$ is reducible and  
		\begin{equation*}
			e_{L/\Q_{r}}=
			\begin{cases}
				(r-1) & \text{if } \frac{r-1}{2} v_{r}(ab) \text{ is even,} \\
				\frac{(r-1)}{2} & \text{if }\frac{r-1}{2} v_{r}(ab)\text{ is odd.} \end{cases}
		\end{equation*}
	\end{enumerate}
\end{thm}

\begin{proof}
	
        $(1)$ Assume first $q \neq r$ with $q \mid ab$. 
        %We know from Proposition~\ref{prop: elephant i0} that $\gminus$ is reducible and $\gamma_{i_0} \in \Q_q$. We claim that $\gamma_0 \in \Q_q(\sqrt{-ab})$. If $\zr \notin \Q_q$, then $i_0 = 0$, and the claim follows. If $\zr \in \Q_q$, then consider the following towers of field extensions
	%\begin{equation} \label{eq:towerestensions}
	%\Q_{q}\overset{d_0\mid 2}{\subseteq} \Q_{q}(\sqrt{-ab})\overset{n_{0} \leq r}{\subseteq} \Q_{q}(\zr^{i_0}\alpha_{0})\qquad \text{ and } \qquad 
	%\Q_{q}\overset{n_{0}'\leq r}{\subseteq} \Q_{q}(\gamma_{i_0})\overset{d_{0}'\mid 2}{\subseteq} \Q_{q}(\zr^{i_0} \alpha_0),
	%\end{equation}
	%where the superscripts denote the degree of the respective extension. Notice that since $\zr \in \Q_q$, we have $\Q_{q}(\zr^{i_0}\alpha_{0})=\Q_{q}(\alpha_{0})$. 
    %Hence $\Q_{q}(\sqrt{-ab})\overset{n_{0} \leq r}{\subseteq} \Q_{q}(\alpha_{0})$ and since it is a Kummer extension and we assume $\zr \in \Q_q$, we have that $n_{0} \in \{1, r\}$.
    %From \eqref{eq:towerestensions}, we have the equality of degree estensions $d_0 \, n_0 = d_0' \, n_0'$. 
    %But since $\gamma_{i_{0}} \in \Q_{q}$, we have $n_{0}'=1$ and the equality of degrees is $d_0 \, n_0 = d_0'$ with $d_{0}' \mid 2$. This implies that $n_{0} = 1$; as otherwise we would have $n_{0} = r \geq 5$ and $d_{0} \, r > d_{0}'$, giving a contradiction. In particular $\zr^{i_0} \alpha_0 \in \Q_{q}(\sqrt{-ab})$ and thus also $ \alpha_0 \in \Q_{q}(\sqrt{-ab})$, which implies that $\gamma_{0} \in \Q_{q}(\sqrt{-ab})$. 
    Since $\gminus$ is reducible by Proposition~\ref{prop: elephant i0}, then Theorem~\ref{thm:splittingfield} and Lemma~\ref{lem: gamma0 quadratic} show that $L$ is equal to $K \coloneqq \Q_q(\omega, \tau \sqrt{-ab})$. But the extension $\Q_q(\omega) / \Q_q$ is unramified as $q \neq r$, and so $e_{K / \Q_q} = e_{K / \Q_q(\omega)} \leq 2$ as $K / \Q_q(\omega)$ is quadratic.
        From now on, we assume that $q = r$. By the proof of Lemma~\ref{lem:fielddiagram}, 
	\begin{equation*}
		e_{K / \Q_r(\omega)} = \begin{cases}
			2 & \text{ if } \frac{r-1}{2} v_r(ab) \text{ is even,} \\
			1 & \text{ if } \frac{r-1}{2} v_r(ab) \text{ is odd,}
		\end{cases} \quad \text{ so } \quad e_{K / \Q_r} = \begin{cases}
			r-1 & \text{ if } \frac{r-1}{2} v_r(ab) \text{ is even,} \\
			\frac{r-1}{2} & \text{ if } \frac{r-1}{2} v_r(ab) \text{ is odd.}
		\end{cases}
	\end{equation*}

        $(2)$ Assume now that $r \nmid ab$. Since $\frac{r-1}{2} v_r(ab) = 0$ is even, we see from above that $e_{K / \Q_r} = r-1$. If $\gminus$ is reducible then Theorem~\ref{thm:splittingfield} and Lemma~\ref{lem: gamma0 quadratic} show that the splitting field $L$ equals $K$, so we obtain  $e_{L / \Q_r} = e_{K / \Q_r}$. We assume from now on that $\gminus$ is irreducible, implying $\left[\Q_r(\gamma_0) : \Q_r \right] = r$. Since $r$ is coprime with 2 and $r-1$, by Lemma~\ref{lem:fielddiagram}, it follows that $$[L : K]= \left[\Q_r(\alpha_0,\zr) : \Q_r(\alpha_0^r,\zr) \right] =  \left[\Q_r(\alpha_0) : \Q_r(\alpha_0^r) \right]  = r.$$
        Since $r$ is prime, we have $e_{L/K},\; e_{\Q_r(\zr, \alpha_0) / \Q_r(\zr, \, \alpha_0^r)},\; e_{\Q_r(\alpha_0)/\Q_r(\alpha_0^{r})} \in \{1, r\}$. But $r$ is coprime with 2 and $r-1$, hence, by Lemma~\ref{lem:fielddiagram}, $e_{L/K} =r$ precisely when $e_{\Q_r(\zr, \alpha_0) / \Q_r(\zr, \, \alpha_0^r)}=r$; this is true if and only if $e_{\Q_r(\alpha_0)/\Q_r(\alpha_0^{r})}=r$. We deduce the following equalities of ramification indices:
        \begin{equation*}
            e_{L \, / \, K} \, = \,   e_{\Q_r(\zr, \, \alpha_0) \, / \, \Q_r(\zr, \, \alpha_0^r)} \, = \, e_{\Q_r(\alpha_0) \, / \, \Q_r(\alpha_0^r)}.
        \end{equation*}
        Recall from Lemma~\ref{lem: useful eqs} (6) that $v_r(\alpha_0) v_r(\beta_0)= 0$. Therefore, at least one of $\alpha_0$ and $\beta_0$ is a unit in the ring of integers of $\Q_r(\alpha_0) = \Q_r(\beta_0)$. Hence, we conclude the proof using Theorem \ref{thm: rth power extensions} (we can apply such theorem because $\Q_r(\sqrt{-ab}) / \Q_r$ is unramified, as $r \nmid ab$). This implies that $\Q_r(\alpha_0) / \Q_r(\alpha_0^r)$ is totally ramified. So, $e_{L / K} = r$, and $e_{L / \Q_r} = r(r-1)$. 

         $(3)$ Assume that $r \mid ab$. Then, by Proposition~\ref{prop: elephant i0}, the polynomial $\gminus$ is reducible and $\gamma_{i_0} \in \Q_q$, and since $q=r$ we have $i_{0}=0$. Therefore, $L = K$, and we obtain  $e_{L / \Q_r} = e_{K / \Q_r}$.
\end{proof}

\medskip

\begin{thm}\label{thm: rth power extensions}
	Let $F / \Q_r$ be a finite unramified extension of degree $f$, let $u \in \O_F^{*} \setminus \left( \O_F^{*} \right)^{r}$ and consider $M \coloneqq F(u^{1/r})$. Then $M / F$ is totally ramified of degree $r$. Moreover the polynomial $E(x) \coloneqq \left( x + u^{r^{f-1}} \right)^r - u$ is Eisenstein, and the discriminant $\Delta(M / F)$ has $r$-adic valuation $v_{r} (\Delta (M / F)) = r$.
\end{thm}

\begin{proof}
	Since $u^{1/r}$ does not belong to $F$, \cite[Section VI, Theorem 9.1]{Lang} implies that the polynomial $P(x) \coloneqq x^r - u$ is irreducible over $F$. We claim that $M / F$ is not Galois. Then by \cite[Chapter III, Theorem 2]{SerreLF}, $M/F$ cannot be unramified, and since it has a prime degree it has to be totally ramified. Indeed, the splitting field of $P$ is $F(\zr, u^{1/r})$. Note that $\zr \notin F$ as $F / \Q_r$ is unramified.
    Hence, $\zr \notin M$, implying that $M$ is strictly contained in the splitting field of $P$.
    %This implies that $[F(\zr) : F] \mid r-1$. However, $P$ is irreducible, so $[M : F] = r$ and $M \cap F(\zr) = F$. Therefore, $M$ is strictly contained in $F(\zr, u^{1/r})$, hence proves the claim.
	
	Let us prove now that $E(x) = P(x + u^{r^{f-1}})$ is Eisenstein. First of all, notice that $E$ is irreducible as $P$ is. Using the binomial formula, we have 
	\begin{equation*}
		E(x) = x^r + \sum_{j = 1}^{r-1} \binom{r}{j} \left( u^{r^{f-1}}\right)^{j} x^{r-j} + u^{r^f} - u.
	\end{equation*}
	For $0 \leq j \leq r$, denote by $E_{j}$ the $j$-th coefficient of $E(x)$, so that $E(x) = \sum_{j = 0}^{r} E_j x^j$. We denote by $\mathfrak{p}_r$ the unique prime of $\cO_F$. The extension $F / \Q_r$ is unramified, so if $j \neq 0, r$, we have $v_{\mathfrak{p}_r} (E_j) = v_r \left( \binom{r}{j}\right)$ as $u$ is a unit in $\O_F$, and $v_r \left( \binom{r}{j}\right) = 1$ as $r \| \binom{r}{j}$. In particular $v_{\mathfrak{p}_r} (E_{1}) = 1$. 
    
    On the other hand, $v_{\mathfrak{p}_r} (E_0) = 1$.  Indeed, $F / \Q_r$ is unramified, so the residue field $\O_F / \mathfrak{p}_r$ has order $r^f$. Its group of units is cyclic of order $r^f - 1$, hence, $E_0 = u \left( u^{r^{f} - 1} - 1 \right) \equiv 0 \mod \, \mathfrak{p}_r$. Now, if $E_0 \equiv 0 \bmod \, \mathfrak{p}_r^2$ then the Newton polygon of $E$ would contain two slopes: one joining the points $(0, v_{\mathfrak{p}_r}(E_0)), (1, 1)$, and another one joining $(1, 1), (r, 0)$. This would contradict the irreducibility of $E$. We conclude that $E$ is indeed Eisenstein, implying that $u^{1/r}-u^{r^{(f-1)}}$ is a uniformizer of $M$ (see e.g. \cite[pg. 23]{Frohlich}.) It follows that $\mathcal{O}_{M}=\mathcal{O}_{F}[u^{1/r}-u^{r^{(f-1)}}]=\mathcal{O}_{F}[u^{1/r}]$ ,and so $\Delta(M/F)= \disc(x^r-u)$ up to units. Finally, 
	\begin{equation*}
		\disc( x^r - u) = (-1)^{(r(r-1)/2)}r^ru^{(r-1)},  
	\end{equation*}
	and so $v_{r}(M/F)=r$.
\end{proof}

\begin{cor}\label{cor:discriminant}
    The discriminant of $\Q_r(\sqrt{-ab}) / \Q_r(\alpha_0^r)$ has $r$-adic valuation
    \begin{equation*}
        v_{r} \left(\Delta(\Q_r(\alpha_0) / \Q_r(\alpha_0^r) \right) = \begin{cases}
            r & \text{ if } r \nmid ab \text{ and } \gminus \text{ is irreducible,} \\
            0 & \text{ otherwise.}
        \end{cases}
    \end{equation*}   
\end{cor}

\begin{proof}
    If $r \mid ab$ or $\gminus$ is reducible, then $\Q_r(\alpha_0) = \Q_r(\alpha_0^r)$ (see proof of Lemma~\ref{lem: gamma0 quadratic} and Proposition~\ref{prop: elephant i0}). If $r \nmid ab$ and $\gminus$ is irreducible, then $\Q_r(\alpha_0^r) / \Q_r$ is unramified, and one of $\alpha_0^r$ and $\beta_0^r$ is a unit in $\O_{\Q_r(\alpha_0^r)} = \O_{\Q_r(\beta_0^r)}$ by Lemma~\ref{lem: useful eqs}. The result follows from Theorem~\ref{thm: rth power extensions}.
\end{proof}

%%%%%%%%%%%%%%%%%%%%%%%%%%%%%%%%%%%%%%%%%%%%%%%%%%%%%%%%%%%%%%%%%%
\subsection{Conductor exponents for \texorpdfstring{$\Cminus$}{C-} and \texorpdfstring{$\Cplus$}{C+} at odd places}\label{sec: Conductors ppr}~
%%%%%%%%%%%%%%%%%%%%%%%%%%%%%%%%%%%%%%%%%%%%%%%%%%%%%%%%%%%%%%%%%%%%%%%%%%%%%
In this section, we compute the conductor exponent of $\Cminus / \Q_q$, $\Cplus / \Q_q$ and $\Cminus / \Qplus_{\qQzrplus}$, $\Cplus / \Qplus_{\qQzrplus}$ at odd places of bad reduction from their cluster pictures. Recall that we denote by $\cR = \{\gamma_0, \dots, \gamma_{r-1}\}$ the set of roots of $g^-(x)$, and $\mathcal{R}^{+} = \mathcal{R}\cup \{\lastroot\}$ the set of roots of $\gplus(x)$, where $\lastroot = -2c$.

\begin{thm}\label{thm: conductorCminusQ}
The conductor exponent of $\Cminus / \Q$ at an odd prime of bad reduction $q$ is 
\begin{equation*}
    \condexp{\Cminus / \Q_q} = \begin{cases}
        \frac{r-1}{2} & \text{ if } q \neq r \text{ and } q \mid ab, \\ 
        r-1 & \text{ if } q = r, \, \text{and } \gminus \text{ is reducible over } \Q_r \, ,\\
        r & \text{ if } q = r, \, \text{and } \gminus \text{ is irreducible over } \Q_r \, . \\
        %r-1 & \text{ if } q = r, \text{ and } r \mid ab.        
    \end{cases}
\end{equation*}
\end{thm}

\begin{proof}
    \textbf{Case 1:} Assume first $q \neq r$, with $q \mid ab$. We check that the cluster picture of $\Cminus / \Q_q$, shown in Theorem~\ref{thm: ClusterCminusQ}(1), satisfies the semistability criterion (Definition \ref{def:ss}). By Theorem~\ref{thm: Ramif index split field}, the extension $\Q_{q}(\mathcal{R})/\Q_{q}$
    has ramification degree at most $2$. The proper clusters are $\mathcal{R}$ and the twins.
    %$\mathfrak{t}_{j} \coloneqq \lbrace \gamma_{j}, \gamma_{-j} \rbrace$ for $1 \leq j \leq \frac{r-1}{2}$. 
    The former is clearly invariant under the action of inertia. The twins $\mathfrak{t}_{j}$ are also invariant under inertia: by Theorem~\ref{thm:tameindex} we have $d_{\Rroots}^{\ast} =0$ or $1$. Hence $d_{\Rroots}^{\ast} \in \Z$, and $[I_{\Q_{r}}:I_{\mathfrak{t_{j}}}]=1$ since the only cluster strictly containing $\mathfrak{t}_{j}$ is $\Rroots$.
    Finally, the only principal cluster is $\mathcal{R}$, which has depth $d_{\Rroots} = 0$. Using the definition of $\nu_{\Rroots}$ in Definition \ref{defclusterinv}, one checks that $\nu_{\Rroots} = 0$. We deduce that $\Cminus / \Q_q$ is semistable by Theorem~\ref{thm:ss}, and Theorem~\ref{thm:condss} now implies that $\condexp{\Cminus / \Q_q}$ equals the number of twins appearing in its cluster picture, \ie \; $\frac{r-1}{2}$.
    
In the next three cases, we use Theorem \ref{thm: allcond} to compute conductor exponents. 

    \textbf{Case 2:} Assume now $q = r$, with $r \nmid ab$ and $\gminus$ reducible over $\Q_r$. By Theorem~\ref{thm: Ramif index split field}, $\Q_r(\Rroots) / \Q_r$ is a tame extension. Then, Proposition~\ref{cor:nwild0} implies that $\condwild{\Cminus / \Q_r} = 0$. We now compute $\condtame{\Cminus / \Q_r}$. The cluster picture depicted in Theorem~\ref{thm: ClusterCminusQ} shows that the candidates for $U$, \ie \; odd clusters $\neq \Rroots$, are the singletons $\lbrace \gamma_j \rbrace$. Fix $0\leq j \leq r-1$, then $\mathcal{P}( \lbrace \gamma_{j} \rbrace )=\mathcal{R}$, which has depth $d_{\Rroots} = \frac{1}{r-1}$. Since the inertia group stabilizes $\Rroots$, we have $\left[ I_{\Q_r} : I_{\Rroots} \right] = 1$. One can compute $\tilde{\lambda}_{\mathcal{R}}=\frac{r}{2(r-1)}$, and then we have \begin{equation*}
        \xi_{\mathcal{R}}(\tilde{\lambda}_{\mathcal{R}}) = v_2(2(r-1)) > v_2(r-1) = \xi_{\mathcal{R}}(d_{\mathcal{R}}).
    \end{equation*}
    Therefore $\lbrace \gamma_j \rbrace \notin U$, for any $0 \leq j \leq r-1$, so $U = \emptyset$. To compute $V$, observe that the only proper non-\"ubereven cluster is $\Rroots$. Recall that $\xi_{\Rroots}(\tilde{\lambda}_{\Rroots})
    \neq 0$, so $\Rroots \notin V$, and thus $V = \emptyset$. Since $\Cminus$ has genus $\frac{r-1}{2}$ and $\gminus$ is monic, we conclude that $\condtame{\Cminus / \Q_r} = r-1$.

    \textbf{Case 3:} Assume that $q = r$, with $r \nmid ab$ and $\gminus$ irreducible over $\Q_r$. The computation of $\condtame{\Cminus / \Q_r}$ is the same as in the previous case, since the clusters are the same and inertia still stabilizes $\Rroots$. We now compute $\condwild{\Cminus / \Q_r}$. Since $\gminus$ is irreducible, $G_{\Q_{r}}$ acts transitively on $\mathcal{R}$ and there is a single orbit in $\mathcal{R}/G_{\Q_{r}}$. Theorem~\ref{thm: allcond} then implies that 
    \begin{equation*}
    \condwild{\Cminus/\Q_{r}} = 
    v_r(\Delta_{\Q_{r}(\gamma_0)/\Q_{r}}) - \left[ \Q_r(\gamma_0) : \Q_r \right] + f_{\Q_{r}(\gamma_0)/\Q_{r}}.
    \end{equation*}
    Since $\gminus$ is irreducible, $\left[ \Q_r(\gamma_0) : \Q_r \right] = r$. Recall that $\Q_r(\sqrt{-ab}) = \Q_r(\alpha_0^r)$. By Theorem~\ref{thm: rth power extensions}, $\Q_r(\alpha_0) / \Q_r(\alpha_0^r)$ is totally ramified. By Lemma~\ref{lem:fielddiagram}, we see that $\Q_r(\gamma_0) / \Q_r$ is totally ramified which implies that $f_{\Q_{r}(\gamma_0)/\Q_{r}} = 1$. By Corollary~\ref{cor:discriminant}, we have $v_{r}(\Delta_{\Q_{r}(\alpha_{0}) / \Q_{r}(\sqrt{-ab})}) = r$.
    
    \begin{itemize}
        \item If $\sqrt{-ab} \in \Q_r$ then $\Q_r(\gamma_0)=\Q_r(\alpha_0)$. In this case, $\mathfrak{p}_{r} = (r)$ is the unique prime of the ring of integers of $\Q_{r}(\sqrt{-ab})=\Q_{r}$ and $v_{r}(\Delta_{\Q_{r}(\gamma_{0})/\Q_{r}}) =r$. 
        
        \item If $\sqrt{-ab} \notin \Q_r$, since $r \nmid ab$ then $\Q_r(\sqrt{-ab}) / \Q_r$ is unramified, so $v_r(\Delta_{\Q_r(\sqrt{-ab}) / \Q_r}) = 0$. On the other hand, by Corollary~\ref{cor:discriminant}, $\Delta_{\Q_r(\alpha_0) / \Q_r(\sqrt{-ab})}$ equals $r^r$ up to a unit.   Since $\Norm_{\Q_r(\sqrt{-ab}) / \Q_r}(r) = r^2$ then the growth of the discriminant in tower extensions (see \cite[Corollary 2.10, p. 202]{Neukirch}) gives $v_r \left( \Delta_{\Q_r(\alpha_0) / \Q_r} \right) = 2r$. Applying the same formula to the tower of extensions $\Q_r \subset \Q_r(\gamma_0) \subset \Q_r(\alpha_0)$ gives $v_{r}(\Delta_{\Q_{r}(\gamma_{0}) / \Q_r})=r$ as $\Q_r(\alpha_0)/ \Q_r(\gamma_0)$ is unramified. 
    \end{itemize}
    As a consequence we obtain in both cases $\condwild{\Cminus/\Q_{r}} =1$, Thus $\condexp{\Cminus / \Q_r} = r$.

    \textbf{Case 4:} Assume now that $q = r$ and $r \mid ab$. By Theorem~\ref{thm: Ramif index split field}, the extension $\Q_r(\Rroots) / \Q_r$ is tame. Therefore, Proposition~\ref{cor:nwild0} yields $\condwild{\Cminus / \Q_r} = 0$. We now compute $\condtame{\Cminus / \Q_r}$. For $ 1 \leq j \leq \frac{r-1}{2}$, we denote by $\mathfrakt_j$ the twin $\mathfrakt_j \coloneqq \left \lbrace \gamma_j, \gamma_{-j} \right \rbrace$. As in Cases 2 and 3, the candidates for $U$, \ie \; odd clusters $\neq \Rroots$, are the singletons $\lbrace \gamma_j \rbrace$. On the other hand, candidates for $V$ are $\Rroots$ and $\mathfrakt_j$ for $1 \leq j \leq \frac{r-1}{2}$.
    
    We first show that $\lbrace \gamma_0 \rbrace \notin U$. Its parent is $\Rroots$ which has depth $d_{\Rroots}=\frac{2}{r-1}$, and is stabilized by inertia i.e. $[I_{\Q_{r}}:I_{\Rroots}]=1$. One can compute $\tilde{\lambda}_{\Rroots} = \frac{1}{r-1}$ and then $\xi_{\mathcal{R}}(\tilde{\lambda}_{\mathcal{R}}) = v_2(r-1) > v_2(\frac{r-1}{2}) = \xi_{\mathcal{R}}(d_{\mathcal{R}})$. We deduce that $\gamma_0 \notin U$. Fix now $1 \leq j \leq r-1$. The parent of $\lbrace \gamma_j \rbrace$ is $\mathfrakt_j$ if $1 \leq j \leq \frac{r-1}{2}$, and $\mathfrakt_{r-j}$ if $\frac{r+1}{2} \leq j \leq r-1$. For simplicity, assume that $1 \leq j \leq \frac{r-1}{2}$. Observe first that $\mathfrakt_j$ is not an orphan of $\Rroots$. Indeed, $I_{\mathcal{R}}=I_{L}$ where $L = \Q_r(\Rroots)$. By the computation of the ramification indices in Theorem~\ref{thm: Ramif index split field} we have either $I_{L}/I_{\Q_{r}}$ has order $r-1$ or $\frac{r-1}{2}$. Hence, in the case $e_{L/\mathbb{Q}_r} = r-1$, the elements $\sigma_j:\zeta_r \to \zeta_r^j$ with $1 \leq j \leq r-1 $ generate $I_{L}/I_{\Q_{r}}$; and in the case 
    where $e_{L/\Q_r} = \frac{r-1}{2}$, only the elements $\sigma_j:\zeta_r \to \zeta_r^j$ with $\sigma_j$ of order smaller or equal than $\frac{r-1}{2}$ give $I_{L}/I_{\Q_{r}}$. As a result, the orbits of the roots of $\Rroots \setminus \{\gamma_{0}\}$ either have length $r-1$ or $\frac{r-1}{2}$. 
    Note that $\mathfrak{t}_{j}$ is not fixed under the action of $I_{\Rroots}$ in either case, and is therefore not an orphan of $\Rroots$. 
    %In either case, no twin $\mathfrak{t}_{j}$ is fixed under the action of $I_{\Rroots}$ and so no $\mathfrak{t}_{j}$ is orphan of $\Rroots$. 
    Since the only cluster strictly containing $\mathfrakt_j$ is $\Rroots$, Theorem~\ref{thm:tameindex} yields $[I_{\mathbb{Q}_r} : I_{\mathfrakt_j}] = \mathrm{denom } \ d_{\Rroots} = \mathrm{denom }\  \frac{2}{r-1} = \frac{r-1}{2}$. Recall from Theorem~\ref{thm: ClusterCminusQ} that $d_{\mathfrak{t}_{j}}=\frac{p}{2} v_{r}(ab)-\frac{r-2}{r-1}$. One can compute $\tilde{\lambda}_{\mathfrak{t}_{j}}=\frac{p}{2}v_{r}(ab)$. We deduce that
    \begin{equation*}
        \xi_{\mathfrak{t}_{j}}(\tilde{\lambda}_{\mathfrak{t}_{j}})=
        \begin{cases}
        1& \text{ if } v_{r}(ab)\frac{r-1}{2} \text{ is odd}\\
        0& \text{ if } v_{r}(ab)\frac{r-1}{2} \text{ is even}
        \end{cases} \quad \text{ and } \quad \xi_{\mathfrak{t}_{j}}(d_{\mathfrak{t}_{j}})=
        \begin{cases}
        0& \text{ if } v_{r}(ab)\frac{r-1}{2} \text{ is odd}\\
        1& \text{ if } v_{r}(ab)\frac{r-1}{2} \text{ is even.}
        \end{cases}
    \end{equation*}
    If $v_{r}(ab)\frac{r-1}{2}$ is odd then $U=V=\emptyset$, and we deduce that $\condtame{\Cminus / \Q_r} =r-1$. On the other hand, if $v_{r}(ab)\frac{r-1}{2}$ is even then 
    \begin{equation*}
        U=\left \lbrace \{\gamma_{j}\}\; :\; 1\leq j\leq r-1 \right \rbrace \quad \text{ and } \quad V= \left \lbrace \mathfrakt_{j} \; :\; 1\leq j\leq \frac{r-1}{2} \right \rbrace.
    \end{equation*}
    But in this case, Theorem~\ref{thm: Ramif index split field} yields $e_{L/\Q_r}=r-1$ giving that the elements $\sigma_j:\zeta_r \to \zeta_r^j$ belong to $I_{\Q_r}$. Therefore $\gamma_{1}, \ldots, \gamma_{r-1}$ lie in a unique orbit under the action of $I_{\Q_r}$. Therefore $\# U/I_{\Q_{r}}=\# V/I_{\Q_{r}}=1$ and, we obtain $\condtame{\Cminus/\Q_{r}}=r-1$.
\end{proof}

\begin{thm}\label{thm: conductorCminusQpl}
	Let $\qQzrplus$ be an odd prime of $\Qplus$ of bad reduction for $\Cminus/\Qzrplus$. 
The conductor exponent of $\Cminus / \Qplus$ at $\qQzrplus$ is 
\begin{equation*}
    \condexp{\Cminus / \Qplus_{\qQzrplus}} = \begin{cases}
        \frac{r-1}{2} & \text{ if } \qQzrplus \nmid r \text{ and } q \mid ab, \\ 
        r-1 & \text{ if } \qQzrplus \mid r, \text{ and } \gminus \text{ is reducible over } \Q_r \, ,\\
        \frac{3(r-1)}{2} & \text{ if } \qQzrplus \mid r, \text{ and } \gminus \text{ is irreducible over } \Q_r \, . \\
        %r-1 & \text{ if } \qQzrplus = \rQzrplus, \text{ and } \, r \mid ab.
    \end{cases}
\end{equation*}
\end{thm}

\begin{proof}
We know from \eqref{eq: Isom compl+ext} that $\Qplus_{\qQzrplus} \simeq \Q_q(\omega)$. Moreover, by Theorem~\ref{thm:splittingfield}, $\omega$ belongs to the splitting field $\Q_q(\Rroots)$, so we can identify $\Qplus_{\qQzrplus}(\Rroots)$ as an extension of $\Q_q(\omega)$. Recall that $\mathfrak r$ is the unique prime above $r$ in $\Qplus$. 

\textbf{Case 1:} First assume $\qQzrplus \nmid r$, with $\qQzrplus \mid ab$. $\Cminus/\Qzrplus$ is semistable by the proof of Theorem~\ref{thm: conductorCminusQ}, case 1. The result follows immediately from \cite[Proposition 3.15]{Liu} and Theorem~\ref{thm:condss}. 
%The curve $\Cminus/\Q_q$ is semistable by the proof of Theorem~\ref{thm: conductorCminusQ}.
%Hence $\Cminus/\Qplus_\mathfrakq$ is also semistable by \cite[Proposition 3.15]{Liu}, and Theorem~\ref{thm:condss} implies that $\condexp{\Cminus / \Qplus_{\mathfrakq}}$ equals the number of twins appearing in its cluster picture, \ie \; $\frac{r-1}{2}$ as seen in Corollary~\ref{cor: ClusterCminusQpl}.

In the next cases, as $\Qplus_{\rQzrplus}/\Q_r$ is tamely ramified, the computation of $\condwild{\Crrp / \Qplus_{\rQzrplus}}$ follows immediately from the proof of Theorem~\ref{thm: conductorCminusQ} and Lemma~\ref{lem:wildthoughtame}, therefore we focus on the computation of $\condtame{\Crrp / \Qplus_{\rQzrplus}}$.

\textbf{Case 2:} Assume now $\qQzrplus\mid r$ (so that $\qQzrplus=\mathfrak r$), with $r \nmid ab$ and $\gminus$ reducible over $\Q_r$. 
%Theorem~\ref{thm: Ramif index split field} states that $\Qplus_{\rQzrplus}(\Rroots) / \Qplus_{\rQzrplus}$ is tamely ramified. Then from Proposition~\ref{cor:nwild0}, we deduce that $\condwild{\Cminus / \Qplus_{\rQzrplus}} = 0$. We now compute $\condtame{\Cminus / \Qplus_{\rQzrplus}}$. 
Corollary~\ref{cor: ClusterCminusQpl} shows that the candidates for $U$, \ie \; odd clusters $\neq \Rroots$, are the singletons $\lbrace \gamma_j \rbrace$. Fix $0\leq j \leq r-1$. Then $\mathcal{P}(\{\gamma_{j}\})=\mathcal{R}$ which has depth $d_{\Rroots} = \frac{1}{2}$. Since the inertia group stabilizes $\Rroots$, we have $[ I_{\Qplus_{\rQzrplus}} : I_{\Rroots} ] = 1$. One can compute $\tilde{\lambda}_{\mathcal{R}}=\frac{r}{4}$, and moreover 
\begin{equation*}
    \xi_{\mathcal{R}}(\tilde{\lambda}_{\mathcal{R}}) = v_2(4) > v_2(2) = \xi_{\mathcal{R}}(d_{\mathcal{R}}).
\end{equation*}
Therefore $\lbrace \gamma_j \rbrace \notin U$ for any $0 \leq j \leq r-1$,  and $U = \emptyset$. To compute $V$, observe that the only proper non-\"ubereven cluster is $\Rroots$. Recall that $\xi_{\Rroots}(\tilde{\lambda}_{\Rroots}) = 2 \neq 0$. Thus $\Rroots \notin V$, and $V = \emptyset$. Recall that $\Cminus$ has genus $\frac{r-1}{2}$ and $\gminus$ is monic. We conclude that $\condtame{\Cminus / \Qplus_{\rQzrplus}} = r-1$.

\textbf{Case 3:} Assume now $\qQzrplus\mid r$ (so that $\qQzrplus=\mathfrak r$), with $r \nmid ab$ and $\gminus$ irreducible over $\Q_r$. The description of the cluster picture is the same as in the previous case, and inertia still stabilizes $\Rroots$. Therefore $\condtame{\Cminus / \Qplus_{\rQzrplus}} = r-1$. 
%Now the extension $\Qplus_{\rQzrplus} / \Q_q$ is tamely ramified, and we know from the proof of Theorem~\ref{thm: conductorCminusQ} that $\condwild{\Cminus / \Q_q} = 1$. We apply Lemma~\ref{lem:wildthoughtame} to get $\condwild{\Cminus / \Qplus_{\rQzrplus}} = \frac{r-1}{2}$, and we obtain $\condexp{\Cminus / \Qplus_{\rQzrplus}} = r - 1 + \frac{r-1}{2} = \frac{3(r-1)}{2}$.

\textbf{Case 4:} Assume now that $\qQzrplus\mid r$ (so that $\qQzrplus=\mathfrak r$) and $r \mid ab$. 
%We know from Theorem~\ref{thm: Ramif index split field} that the extension $\Qplus_{\rQzrplus}(\Rroots) / \Qplus_{\rQzrplus}$ is tamely ramified. Therefore, Proposition~\ref{cor:nwild0} implies that $\condwild{\Cminus / \Qplus_{\rQzrplus}} = 0$. We now compute $\condtame{\Cminus / \Qplus_{\rQzrplus}}$. 
For $ 1 \leq j \leq \frac{r-1}{2}$, we denote by $\mathfrakt_j$ the twin $\mathfrakt_j \coloneqq \left \lbrace \gamma_j, \gamma_{-j} \right \rbrace$. As in Cases 2 and 3 of this proof, the candidates for $U$, \ie \; odd clusters $\neq \Rroots$, are the singletons $\lbrace \gamma_j \rbrace$. On the other hand, the candidates for $V$ are $\Rroots$ and $\mathfrakt_j$ for $1 \leq j \leq \frac{r-1}{2}$. 

Let us first show that $\lbrace \gamma_0 \rbrace \notin U$. Recall from Corollary~\ref{cor: ClusterCminusQpl} that the parent parent of $\gamma_0$ is $\Rroots$ which has depth $d_{\Rroots}= 1$, and is stabilized by inertia i.e. $[ I_{\Qplus_{\rQzrplus}} : I_{\Rroots}] = 1$. One can check that $\tilde{\lambda}_{\Rroots} = \frac{r}{2}$ and then $\xi_{\mathcal{R}}(\tilde{\lambda}_{\mathcal{R}}) = v_2(2) > v_2(1) = \xi_{\mathcal{R}}(d_{\mathcal{R}})$. We deduce that $\gamma_0 \notin U$. Fix now $1 \leq j \leq r-1$. The parent of $\lbrace \gamma_j \rbrace$ is $\mathfrakt_j$ if $1 \leq j \leq \frac{r-1}{2}$, and $\mathfrakt_{r-j}$ if $\frac{r+1}{2} \leq j \leq r-1$. For simplicity, assume that $1 \leq j \leq \frac{r-1}{2}$. The only cluster strictly containing $\mathfrakt_j$ is $\Rroots$, and since $d_{\Rroots} = 1$, we have $d^{\ast}_{\Rroots} = 1$ whether $\mathfrakt_j$ is an orphan or not. We deduce from Theorem~\ref{thm:tameindex} that $[I_{\Qplus_{\rQzrplus}} : I_{\mathfrakt_j}]=1$. Recall from Corollary~\ref{cor: ClusterCminusQpl} that $d_{\mathfrak{t}_{j}}=\frac{p}{2}\cdot \frac{r-1}{2}v_{r}(ab)-\frac{r-2}{2}$. We compute $\tilde{\lambda}_{\mathfrak{t}_{j}} = \frac{p(r-1)}{4} \, v_{r}(ab)$. We deduce that 
\begin{equation*}
    \xi_{\mathfrak{t}_{j}}(\tilde{\lambda}_{\mathfrak{t}_{j}})=
    \begin{cases}
    1& \text{ if } v_{r}(ab)\frac{r-1}{2} \text{ is odd},\\
    0& \text{ if } v_{r}(ab)\frac{r-1}{2} \text{ is even},
    \end{cases} \quad \text{ and } \quad \xi_{\mathfrak{t}_{j}}(d_{\mathfrak{t}_{j}})=
    \begin{cases}
    0& \text{ if } v_{r}(ab)\frac{r-1}{2} \text{ is odd},\\
    1& \text{ if } v_{r}(ab)\frac{r-1}{2} \text{ is even.}
    \end{cases}
\end{equation*}
If $v_{r}(ab)\frac{r-1}{2}$ is odd then $U=V=\emptyset$, and we deduce that $\condtame{\Cminus / \Qplus_{\rQzrplus}} =r-1$. On the other hand, if $v_{r}(ab)\frac{r-1}{2}$ is even then
\begin{equation*}
    U=\left \lbrace \lbrace \gamma_{j} \rbrace \; :\; 1\leq j\leq r-1 \right \rbrace \quad \text{ and } \quad V= \left \lbrace \mathfrakt_{j} \; :\; 1\leq j\leq \frac{r-1}{2} \right \rbrace.
\end{equation*}
In this case, Theorem~\ref{thm: Ramif index split field} gives $e_{L/\Qplus_{\rQzrplus}}=2$, so $I_{\Qplus_{\mathfrakr}}$ stabilizes each twin $\mathfrakt_j$, and each of these lies in a different orbit under the action of $I_{\Qplus_{\rQzrplus}}$. Hence $\# V / I_{\Qplus_{\rQzrplus}} = \frac{r-1}{2}$. On the other hand, the equality $[I_{\Qplus_{\rQzrplus}} : I_{\mathfrakt_j}]=1$ implies that singletons are regrouped by pairs under inertia. Therefore, $\# U / I_{\Qplus_{\rQzrplus}} = \frac{r-1}{2}$. We conclude that $\condtame{\Cminus / \Qplus_{\rQzrplus}} = r-1$ by Theorem~\ref{thm: allcond}.
\end{proof}

Just as in Corollary~\ref{cor:twistCr}, one can find a curve with lower conductor exponent.

\begin{cor} \label{cor:twistCmius}
    When $r \mid ab$, the twist of $\Cminus$ by a uniformizer of $\Qplus_{\rQzrplus}$ is semistable and has conductor exponent $\frac{r-1}{2}$.
\end{cor}

\begin{proof}
    Just as in the proof of Corollary~\ref{cor:twistCr}, the cluster picture of the twist is the same one as that of $\Cminus / \Qzrplus_{\rQzrplus}$, but the leading term of the defining polynomial has valuation $1$.
    Theorem~\ref{thm: Ramif index split field},states that, when $r \mid ab$, we have $e_{\Qplus_{\rQzrplus}(\Rroots) / \Qplus_{\rQzrplus}} \leq 2$. The description of the cluster picture and the proof of Theorem~\ref{thm: ClusterCminusQ} shows that $I_{\Qplus_{\rQzrplus}}$ does not permute the twins. One checks that $\nu_{\Rroots}$ for the twist is even, so the semistability criterion (Definition \ref{def:ss}) is satisfied for the twist, and the computation of the conductor exponent follows from Theorem~\ref{thm:condss}.
\end{proof}

\begin{theorem}\label{thm: conductorCplusQ}
The conductor exponent of $\Cplus/\Q$ at an odd prime of bad reduction $q$ is
\begin{equation*}
 \condexp{\Cplus/\Q_{q}} = \begin{cases}
        \frac{r-1}{2} & \text{if } q \neq r \text{ and } q \mid ab, \\
        r-1 & \text{if } q=r, \ r \nmid ab  \text{ and } \gminus \text{ is reducible over } \Q_r \, , \\
        r & \text{if } q=r, \ r\nmid ab \text{ and } \gminus \text{ is irreducible over } \Q_r \, , \\
        r-2 & \text{if } q=r, \text{ and } r \mid a,\\
        r-1 & \text{if } q=r, \text{ and } r \mid b.\\
 \end{cases}
\end{equation*}
\end{theorem}

\begin{proof} By Remark \ref{rem:wildplus}, $\condwild{\Cplus/\Q_q}=\condwild{\Cminus/\Q_q}$. The latter is given by the proof of Theorem~\ref{thm: conductorCminusQ}, and we, therefore, focus on the computations of $\condtame{\Cplus/\Q_q}$. 

\textbf{Case 1:} Assume first $q \neq r$, with $q \mid ab$. Suppose $q\mid a$. Let $\mathfrakt \coloneqq \lbrace \gamma_{i_0}, \lastroot \rbrace$; this is the new twin appearing in the cluster picture of $\Cplus / \Q_q$ given by Theorem~\ref{thm: ClusterCplusQ}, and is not present in the cluster picture of $\Cminus / \Q_q$. Since both $\gamma_{i_0}, \lastroot \in \Q_q$ (see Proposition~\ref{prop: elephant i0} for justification of the inclusion of the former root), the twin $\mathfrakt$ is clearly invariant under the action of $I_{\Q_q}$. Moreover, $d_{\mathfrakt} = \nu_{\mathfrakt}= 2p v_q(ab) \in \Z$. We deduce that $\Cplus / \Q_q$ satisfies the semistability criterion (Definition\ref{def:ss}) so that it is semistable by Theorem \ref{thm:ss}. There are $\frac{r+1}{2}$ twins in the cluster picture, and since $\Rroots$ is \"ubereven, the conductor exponent is $\condexp{\Cplus / \Q_q} = \frac{r-1}{2}$.

Now suppose $q\mid b$. As in Theorem~\ref{thm: ClusterCplusQ}, the cluster picture of $\Cplus/ \Q_{q}$ differs by the cluster picture of $\Cminus/ \Q_{q}$ by the singleton $\{\lastroot\}$. In this case, $\condexp{\Cplus/ \Q_{q}}=\frac{r-1}{2}$ by case 1 of the proof of Theorem~\ref{thm: conductorCminusQ} in combination with Theorem~\ref{thm:ss}.

\textbf{Case 2:}
Assume that $q = r$, with $r \nmid ab$. The tame part of the conductor does not depend on the reducibility of $\gminus$.
The cluster picture depicted in Theorem~\ref{thm: ClusterCplusQ} shows that the candidates for $U$, \ie \; odd clusters $\neq \Rroots^+$, are the singletons 
$\{\gamma_{j}\}$ for $0\leq j\leq r$ and $\mathcal{R}$.  Observe that $\mathcal{P}(\{\lastroot\})=\mathcal{P}(\mathcal{R})=\mathcal{R^+}$ which has depth $d_{\mathcal{R}^{+}}=0$. 
Since the inertia group stabilizes $\Rroots^+$ we have $\left[ I_{\Q_r} : I_{\Rroots^+} \right] = 1$. One can compute $\tilde{\lambda}_{\mathcal{R}^+}=0$ and then we have \begin{equation*}
        \xi_{\mathcal{R}^+}(\tilde{\lambda}_{\mathcal{R}^+}) =\xi_{\mathcal{R}^+}(d_{\mathcal{R}^+})=0.
    \end{equation*}
Thus $\{\lastroot\}, \mathcal{R}\in U$. For any $0\leq j\leq r-1$, then $\mathcal{P}(\{\gamma_{j}\})=\mathcal{R}$. 
    From the proof of Theorem~\ref{thm: conductorCminusQ} cases 2 and 3, we see that $\{\gamma_{j} \} \not\in U$. 
    Thus $U=\{\mathcal{R},\; \{\lastroot\}\}$. 
    
    To compute $V$, observe that the only proper non-\"ubereven clusters are $\Rroots$ and $\cR^+$. Recall that $\xi_{\Rroots^+}(\tilde{\lambda}_{\Rroots^+})
    = 0$, and as computed in the proof of Theorem~\ref{thm: conductorCminusQ} cases 2 and 3, 
    $\xi_{\mathcal{R}}(\tilde{\lambda}_{\mathcal{R}})\neq 0$. Thus $V = \{ \cR^+ \}$. Trivially, $\Rroots, \Rroots^+$ and $\{\lastroot\}$ are fixed by the action of $I_{\Q_r}$. Therefore $\#(U/I_{\Q_r})=2$ and $\#(V/I_{\Q_r})=1$. Since $\Cplus$ has genus $\frac{r-1}{2}$, $\gplus$ is monic and $\mid \Rroots^{+} \mid$ is even, we conclude that $\condtame{\Cplus / \Q_r} = r-1$ by Theorem~\ref{thm: allcond}. 
 
\textbf{Case 3:}
Assume now $q = r$ with $r \mid a$. For $ 1 \leq j \leq \frac{r-1}{2}$, we denote by $\mathfrakt_j$ the twin $\mathfrakt_j \coloneqq \left \lbrace \gamma_j, \gamma_{-j} \right \rbrace$, and we let $\mathfrakt \coloneqq \left \lbrace \gamma_0, \lastroot \right \rbrace$. The candidates for $U$ are the singletons $\lbrace \gamma_j \rbrace$. We have $P( \lbrace \gamma_0 \rbrace) = P( \lbrace \lastroot \rbrace) = \mathfrak t$, which has depth $d_{\mathfrakt} = pv_r(a) - 2$. We compute $\tilde{\lambda}_{\mathfrakt} = p v_r(a) - 1$, and $[I_{\Q_r} : I_{\mathfrakt}] = 1$ since $\gamma_0, \lastroot \in \Q_r$. We deduce that $\xi_{\mathfrakt}(d_{\mathfrakt})=0$ and $\xi_{\mathfrakt}(\tilde{\lambda}_{\mathfrakt})=0$, which implies that $\gamma_{0}, \lastroot\in U$. Fix now $1 \leq j \leq r-1$. Recall from the proof of Theorem~\ref{thm: conductorCminusQ} that $P(\lbrace \gamma_j \rbrace) = \mathfrakt_j$, $d_{\mathfrakt_j} = \frac{p}{2} v_r(a) + \frac{2-r}{r-1}$, and $[ I_{\Q_r} : I_{\mathfrakt_j} ] = \frac{r-1}{2}$. We compute $\tilde{\lambda}_{\mathfrakt_j} = \frac{p}{2} v_r(a) + \frac{1}{r-1}$, which gives 
\begin{align*}
    \xi_{\mathfrakt_j}(d_{\mathfrakt_j}) & = \max \left \lbrace 1 -v_{2}\left( \frac{p(r-1)}{2} \, v_{r}(a)+2-r \right), \; 0 \right \rbrace, \\
    \xi_{\mathfrakt_j}(\tilde{\lambda}_{\mathfrakt_j}) & = \max \left \lbrace 1 -v_{2}\left( \frac{p(r-1)}{2} \, v_{r}(a)+ 1 \right), \; 0 \right \rbrace.
\end{align*}
Since $2-r$ and $1$ are both odd, we deduce that $\xi_{\mathfrakt_j}(d_{\mathfrakt_j}) = \xi_{\mathfrakt_j}(\tilde{\lambda}_{\mathfrakt_j})$. Thus, $\lbrace \gamma_j \rbrace \in U$ and we deduce that $U = \left \lbrace \lbrace \gamma_0 \rbrace , \ldots, \lbrace \lastroot \rbrace \right \rbrace$. On the other hand, the candidates for $V$ are $\mathfrakt$ and $\mathfrakt_j$ for $1 \leq j \leq \frac{r-1}{2}$. Recall from above that $\xi_{\mathfrakt}(\tilde{\lambda}_{\mathfrakt})=0$ i.e. $\mathfrakt \in V$. Moreover, one can check that $\xi_{\mathfrakt_j}(\tilde{\lambda}_{\mathfrakt_j}) = 0$ if and only if $\frac{r-1}{2} v_r(a)$ is odd. This yields 
\begin{equation*}
    V = \begin{cases}
        \{\mathfrak{t}\} & \text{if }  v_{r}(a)\frac{r-1}{2} \text{ is even},\\
        \left \lbrace \mathfrak{t}, \mathfrak{t}_{1},\dots, \mathfrak{t}_{\frac{r-1}{2}} \right \rbrace & \text{if } v_{r}(a)\frac{r-1}{2} \text{ is odd}.
\end{cases}
\end{equation*}
Since $\lastroot \in \Q_r$, we have $\Q_r(\Rroots^{+}) = \Q_r(\Rroots)$, and Theorem~\ref{thm: Ramif index split field} gives $e_{\Q_r(\Rroots^{+}) / \Q_r} = r-1$ (resp. $\frac{r-1}{2}$) if $v_r(a) \frac{r-1}{2}$ is even (resp. odd). Hence, $\{\gamma_{1}, \ldots, \gamma_{r-1}\}$ forms one orbit under the action of $I_{\Q_r}$ (resp. breaks into two orbits) if $v_{r}(a)\frac{r-1}{2}$ is even (resp. if $v_{r}(a)\frac{r-1}{2}$ is odd). We deduce
\begin{equation*}
    \# U/I_{\Q_{r}} = \begin{cases}
        3 & \text{if } v_{r}(a)\frac{r-1}{2}\text{ is even},\\
        4 & \text{if } v_{r}(a)\frac{r-1}{2}\text{ is odd},
\end{cases} \qquad \# V/I_{\Q_{r}} = \begin{cases}
        1 & \text{if } v_{r}(a)\frac{r-1}{2}\text{ is even},\\
        2 & \text{if } v_{r}(a)\frac{r-1}{2}\text{ is odd}.
\end{cases}
\end{equation*}
In all cases we obtain $\condtame{\Cplus / \Q_r} = r-2$. 

\textbf{Case 4:}
Assume now $q = r$ with $r \mid b$. The cluster picture depicted in Theorem~\ref{thm: ClusterCplusQ} shows that the candidates for $U$, \ie \; odd clusters $\neq \Rroots^+$, are the singletons 
$\lbrace \gamma_{j} \rbrace $ for $0\leq j\leq r$ and $\mathcal{R}$. Observe that $\mathcal{P}(\mathcal{R})=\mathcal{P}(\{\lastroot\})=\mathcal{R}^{+}$. 
Moreover, $d_{\mathcal{R}^{+}}=\tilde{\lambda}_{\mathcal{R}^{+}}=0$. 
Thus $\mathcal{R}, \{\lastroot\}\in U$. 
From the proof of Theorem~\ref{thm: conductorCminusQ} case 4, we see that $\{\gamma_{0}\}\not\in U$. 
Now let $1\leq j\leq r-1$. 
By the proof of Theorem~\ref{thm: conductorCminusQ} case 4, $\{\gamma_{j}\}\in U$ if and only if $v_{r}(ab)\frac{r-1}{2}$ is even.
Thus
\[
U=
\begin{cases}

\{\mathcal{R}, \{\gamma_j\}: 1\leq j \leq r \} & \text{if } v_{r}(ab)\frac{r-1}{2}\text{ is even},\\
\{\mathcal{R}, \{\lastroot\}\} & \text{if } v_{r}(ab)\frac{r-1}{2}\text{ is odd}.
\end{cases}
\]
 To compute $V$, observe that the only proper non-\"ubereven clusters are the twins $\mathfrak{t}_{j}$ for $1\leq j\leq\frac{r-1}{2}$, 
$\mathcal{R^{+}}$ and $\mathcal{R}$. 
 It is straightforward to check that $\xi_{\cR^+}(\tilde{\lambda}_{\mathcal{R}^{+}})=0$, so $\mathcal{R}^{+}\in V$. 
By the proof of Theorem~\ref{thm: conductorCminusQ} case 4, $\mathfrak{t}_{j}\in V$ if and only if $v_{r}(ab)\frac{r-1}{2}$ is even, and furthermore the same proof implies that $\mathcal{R}\not\in V$.
Thus
\[
V=
\begin{cases}

\{\mathcal{R^{+}} \} \cup \{ \mathfrak{t}_{j}: 1\leq j \leq \frac{r-1}{2}\} & 
\text{if } v_{r}(ab)\frac{r-1}{2}\text{ is even},\\
\{\mathcal{R^{+}}\} & \text{if } v_{r}(ab)\frac{r-1}{2}\text{ is odd}.
\end{cases}
\]
Suppose $v_{r}(ab)\frac{r-1}{2}$ is odd. Then $\Rroots^{+}$ and $\lbrace \lastroot \rbrace$ lie in distinct orbits under the action of inertia so that $\#(U/I_{\Q_{r}})=2$. We deduce that  $\condtame{\Cplus/\Q_{r}}=r-1$. 
Suppose $v_{r}(ab)\frac{r-1}{2}$ is even. 
Then by Theorem~\ref{thm: Ramif index split field}, we have $e_{\Q_{r}(\Rroots^{+})/\Q_{r}}=r-1$, so the roots $\{\gamma_{1},\dots, \gamma_{r-1}\}$ form one orbit under the action of inertia. 
Thus $\#(U/I_{\Q_{r}})=3$. By the same argument, the twins $\{\mathfrak t_1, \dots, \mathfrak t_{(r-1)/2}\}$ form one orbit under the action of inertia, so $\#(V/I_{\Q_{r}})=2$. Since $\mid \Rroots^{+} \mid $ is even, we conclude $\condtame{\Cplus/\Q_{r}}=r-1$.

\end{proof}

\begin{theorem}\label{thm: conductorCplusQpl}
	Let $\qQzrplus$ be an odd prime of $\Qplus$ of bad reduction for $\Cplus/\Qzrplus$. 
The conductor exponent at $\mathfrak{q}$ is
\begin{equation*}
 \condexp{\Cplus/\Qzrplus_{\qQzrplus}} = \begin{cases}
        \frac{r-1}{2} & \text{if } \mathfrakq \nmid r \text{ and } q \mid ab, \\
        r-1 & \text{if } \mathfrakq\mid r, \ r \nmid ab  \text{ and } \gminus \text{ is reducible over } \Q_r \, , \\
        \frac{3(r-1)}{2} & \text{if } \mathfrakq\mid r, \ r\nmid ab \text{ and } \gminus \text{ is irreducible over } \Q_r \, , \\
        \frac{r-1}{2} & \text{if } \mathfrakq\mid r, \text{ and } r \mid a,\\
        r-1 & \text{if } \mathfrakq\mid r, \text{ and } r \mid b.\\
 \end{cases}
\end{equation*}

\end{theorem}

\begin{proof} 
By Remark \ref{rem:wildplus} again, $\condwild{\Cplus/\Qplus_{\qQzrplus}}=\condwild{\Cminus/\Qplus_{\qQzrplus}}$. The latter is computed in the proof of Theorem~\ref{thm: conductorCminusQpl}, and we therefore focus on the computations of $\condtame{\Cplus/\Qplus_{\qQzrplus}}$. 
Recall that $\mathfrak r$ is the unique prime above $r$ in $\Qplus$.

\textbf{Case 1:} Assume first $\qQzrplus \nmid r$, with $\qQzrplus \mid ab$. 
 By the proof of Case 1 in Theorem~\ref{thm: conductorCplusQ}, \cite[Proposition 3.15]{Liu} and Theorem~\ref{thm:condss} the result follows immediately.

%By the proof of Theorem~\ref{thm: conductorCplusQ}, $\Cplus/\Q_q$ is semistable. Therefore $\Cplus/\Qplus_\mathfrakq$ is semistable by \cite[Proposition 3.15]{Liu}. 
%Suppose $r\mid a$. Since $\Rroots^+$ is \"ubereven by Corollary~\ref{cor: ClusterCplusQpl}, Theorem~\ref{thm:condss} implies that $\condexp{\Cplus / \Qplus_{\mathfrakq}}$ equals the number of twins appearing in its cluster picture minus one, \ie \;  $\condexp{\Cplus/\Qplus_{\qQzrplus}} = \frac{r-1}{2}$. Similarly, if $r\mid b$ then $\condexp{\Cplus/\Qplus_{\qQzrplus}} = \frac{r-1}{2}$ by Theorem~\ref{thm:condss}. 

In the next three cases, we use Theorem \ref{thm: allcond} to compute conductor exponents. 

\textbf{Case 2:}
Assume now $\qQzrplus \mid r$ (so that $\qQzrplus=\rQzrplus$) with $r \nmid ab$. The tame part of the conductor does not depend on the reducibility of $\gminus$. The cluster picture depicted in Corollary~\ref{cor: ClusterCplusQpl} shows that the candidates for $U$, \ie \;  odd clusters $\neq \Rroots^+$, are the singletons $\lbrace \gamma_j \rbrace$, and $\cR$. Observe that $\mathcal{P}(\{\lastroot\})=\mathcal{P}(\mathcal{R})=\mathcal{R^+}$ whose depth is 
$d_{\mathcal{R}^{+}}=0$. 
Since the inertia group stabilizes $\Rroots^+$ we have $\left[ I_{\Qplus_{\rQzrplus}} : I_{\Rroots^+} \right] = 1$. One can compute $\tilde{\lambda}_{\mathcal{R}^+}=0$ and then we have \begin{equation*}
        \xi_{\mathcal{R}^+}(\tilde{\lambda}_{\mathcal{R}}) =\xi_{\mathcal{R}^+}(d_{\mathcal{R}})=0.
    \end{equation*}
 Thus $\{\lastroot\}, \mathcal{R}\in U$. Now let $0\leq j\leq r-1$. 
    Then $\mathcal{P}(\{\gamma_{j}\})=\mathcal{R}.$ 
From the proof of Theorem~\ref{thm: conductorCminusQpl} Cases 2 and 3, we see that $\{\gamma_{j}\}\not\in U$, so $U=\{\mathcal{R},\; \{\lastroot\}\}$.

To compute $V$, observe that the only proper non-\"ubereven clusters are $\Rroots$ and $\cR^+$. Recall from above that $\xi_{\Rroots^+}(\tilde{\lambda}_{\Rroots^+})= 0$, and as computed in the proof of Theorem~\ref{thm: conductorCminusQpl} case 2, $\xi_{\mathcal{R}}(\tilde{\lambda}_{\mathcal{R}})\neq 0$. Thus $V = \{ \cR^+ \}$. Trivially,  both $U$ and $V$ are fixed by the action of $I_{\Qplus_{\rQzrplus}}$. Therefore $\#(U/I_{\Qplus_{\rQzrplus}})=2$ and $\#(V/I_{\Qplus_{\rQzrplus}})=1$. 
 Since $\Cplus$ has genus $\frac{r-1}{2}$ and $\gplus$ is monic, and $\mid \Rroots^+ \mid$ is even, we conclude that $\condtame{\Cplus / {\Qplus_{\rQzrplus}}} = r-1$ by Theorem~\ref{thm: allcond}.

\textbf{Case 3:} Assume now $\qQzrplus \mid r$ (so that $\qQzrplus=\rQzrplus$), with $r \mid a$. We check that the cluster picture of $\Cplus / \Qplus_{\rQzrplus}$, shown in Corollary~\ref{cor: ClusterCplusQpl} (4), satisfies the semistability criterion (Definition~\ref{def:ss}).
% given in the proof of .
By Theorem~\ref{thm: Ramif index split field}, the extension $\Qplus_{\rQzrplus}(\mathcal{R}^+)/\Qplus_{\rQzrplus}$ has ramification degree at most $2$. The proper clusters are $\mathcal{R}^+$ and the twins $\mathfrak{t}\coloneqq\{\gamma_0, \lastroot\}$, $\mathfrak{t}_{j} \coloneqq \lbrace \gamma_{j}, \gamma_{-j} \rbrace$ for $1 \leq j \leq \frac{r-1}{2}$. The former is invariant under the action of inertia, as its elements belong to $\Q_r$. Moreover, the twins $\mathfrak{t}_{j}$ for $1 \leq j \leq \frac{r-1}{2}$ are fixed by inertia by the argument in the proof Theorem~\ref{thm: conductorCminusQpl} case 4. The only principal cluster is the whole 
cluster $\mathcal{R}^{+}$ with $d_{\mathcal{R}^{+}}=1\in \Z$ and $v_{\mathcal{R}^{+}}=r+1\in 2\Z$. 
We deduce that $\Cplus / \Qplus_{\rQzrplus} $ is semistable, and Theorem~\ref{thm:condss} implies that $\condexp{\Cplus / \Qplus_{\rQzrplus}}$ equals the number of twins appearing in its cluster picture minus one, \ie \; $\frac{r-1}{2}$.

\textbf{Case 4:} Assume now $\qQzrplus \mid r$ (so that $\qQzrplus=\rQzrplus$), with $r \mid b$. The cluster picture depicted in Corollary~\ref{cor: ClusterCplusQpl} shows that the candidates for $U$, \textit{i.e.} odd clusters $\neq \Rroots^+$, are the singletons 
$\lbrace \gamma_{j} \rbrace $ for $0\leq j\leq r$ and $\mathcal{R}$. Observe that $\mathcal{P}(\mathcal{R})=\mathcal{P}(\{\lastroot\})=\mathcal{R}^{+}$. Moreover, $d_{\mathcal{R}^{+}}=\tilde{\lambda}_{\mathcal{R}^{+}}=0$, so $\mathcal{R}, \{\lastroot\}\in U$. 
From the proof of Theorem~\ref{thm: conductorCminusQpl} case 4, we see that $\{\gamma_{0}\}\not\in U$. 
Now let $1\leq j\leq r-1$. 
By the proof of Theorem~\ref{thm: conductorCminusQpl} case 4, $\{\gamma_{j}\}\in U$ if and only if $v_{r}(ab)\frac{r-1}{2}$ is even.
Thus 
\[U=
\begin{cases}
\{\mathcal{R}, \{\lastroot\}\} & \text{if } v_{r}(ab)\frac{r-1}{2}\text{ is odd},\\
\{\mathcal{R}, \{\gamma_j\}: 1\leq j \leq r \} & \text{if } v_{r}(ab)\frac{r-1}{2}\text{ is even}.
\end{cases}
\]
To compute $V$, observe that the only proper non-\"ubereven clusters are the twins $\mathfrak{t}_{j}$ for $1\leq j\leq\frac{r-1}{2}$, 
$\mathcal{R^{+}}$ and $\mathcal{R}$. 
One can check that $\xi_{\cR^+}(\tilde{\lambda}_{\mathcal{R}^{+}})=0$, so $\mathcal{R}^{+}\in V$. 
By the proof of Theorem~\ref{thm: conductorCminusQpl} case 4, $\mathfrak{t}_{j}\in V$ if and only if $v_{r}(ab)\frac{r-1}{2}$ is even. Furthermore the same proof implies that $\mathcal{R}\not\in V$.
Thus
\[
V=
\begin{cases}
\{\mathcal{R^{+}}\} & \text{if } v_{r}(ab)\frac{r-1}{2}\text{ is odd},\\
\{\mathcal{R^{+}}, \mathfrak{t}_{j}: 1\leq j \leq \frac{r-1}{2}\} & 
\text{if } v_{r}(ab)\frac{r-1}{2}\text{ is even}.
\end{cases}
\]
Suppose $v_{r}(ab)\frac{r-1}{2}$ is odd. Again, $\Rroots^{+}$ and $\{ \lastroot \}$ lie in different orbits under the action of $I_{\Qplus_{\rQzrplus}}$, so $\condtame{\Cplus/\Qplus_{\rQzrplus}}=r-1$. 
Suppose $v_{r}(ab)\frac{r-1}{2}$ is even. 
Then by Theorem~\ref{thm: Ramif index split field}, we have that $e_{\Qzrplus_{\rQzrplus}(\cR^+)/\Qzrplus_{\rQzrplus}} = 2$, giving that the elements $\sigma_j:\zeta_r \to \zeta_r^{-1}$ generates $I_{\Qzrplus}$. Thus the roots $\{\gamma_{1},\dots, \gamma_{r-1}\}$ are paired under the action of inertia, forming $\frac{r-1}{2}$ orbits, and we get $\#(U/I_{\Qplus_{\rQzrplus}})=\frac{r+3}{2}$. By the same argument, 
the twins $\{\mathfrak t_1, \dots, \mathfrak t_{(r-1)/2}\}$ are fixed under the action of inertia, forming $\frac{r-1}{2}$ orbits, so $\#(V/I_{\Qplus_{\rQzrplus}})=\frac{r+1}{2}$. We conclude that $\condtame{\Cplus/\Qplus_{\rQzrplus}}=r-1$. 
\end{proof}

\begin{cor} \label{cor:twistCplus}
    When $r \mid b$, the twist of $\Cplus / \Qplus_{\rQzrplus}$ by a uniformizer of $\Qplus_{\rQzrplus}$ is semistable and has conductor exponent $\frac{r-1}{2}$.
\end{cor}

\begin{proof}
    Just as in the proofs of Corollaries~\ref{cor:twistCr} and \ref{cor:twistCmius}, the cluster picture of the twist is the same as that of $\Cplus / \Qplus_{\rQzrplus}$, but the valuation of the leading term equals $1$. Let us prove that the semistability criterion (Definition \ref{def:ss}) is now satisfied. Since $r \mid b$, we know from Theorem~\ref{thm: Ramif index split field} that $e_{\Qplus_{\rQzrplus}(\Rroots) / \Qplus_{\rQzrplus}} \leq 2$. We also know from the proof of Theorem~\ref{thm: ClusterCminusQ} that $\Rroots$ and all the twins are invariant under the action of $I_{\Qplus_{\rQzrplus}}$. Note that $\Rroots^{+}$ is not principal, so the only principal cluster is $\Rroots$. Since the valuation of the leading term is one, we check that $\nu_{\Rroots}$ for the twist is even. The computation of the conductor exponent follows from Theorem~\ref{thm:condss}.
\end{proof}

\appendix
%%%%%%%%%%%%%%%%%%%%%%%%%%%%%%%%%%%%%%%%%%%%%%%%%%%%%%%%%%%%%%%%%%%%%%%%%%%%%%%%%%%%%%%%%%%%%%%%%%%%%%%%%%%%%%%%%%%%%%%%%%%%%%%%%%%%%%%%%%%%%%%%%%%%%%%%%%
\section{\large Compatibility of conductors, by Martin Azon}\label{sec:compatibility}
%%%%%%%%%%%%%%%%%%%%%%%%%%%%%%%%%%%%%%%%%%%%%%%%%%%%%%%%%%%%%%%%%%%%%%%%%%%%%%%%%%%%%%%%%%%%%%%%%%%%%%%%%%%%%%%%%%%%%%%%%%%%%%%%%%%%%%%%%%%%%%%%%%%%%%%%%%
In this section, we show that the conductor computations obtained in Theorem \ref{thm:conductorQplus} are coherent with those found in \cite[§5]{BCDF} and \cite[§7]{chen2022modular}. This compatibility depends on the Jacobian of the curve studied in \cite{BCDF, chen2022modular} being of $\Gltwo$-type.

%%%%%%%%%%%%%%%%%%%%%%%%%%%%%%%%%%%%%%%%%%%%%%%%%%%%%%%%%%%%%%%%%%%%%%%%%%%%%
\subsection{Abelian varieties of \texorpdfstring{$\Gltwo$}{Gltwo}-type}\label{Gltwo ab var}
%%%%%%%%%%%%%%%%%%%%%%%%%%%%%%%%%%%%%%%%%%%%%%%%%%%%%%%%%%%%%%%%%%%%%%%%%%%%%
We begin by recalling the definition of an abelian variety of $\Gltwo$-type.

\begin{defn}\label{Def Gltwo}
    Let $A$ be an abelian variety defined over a number field $L$. We say that $A / L$ is of $\Gltwo$-type if there exists an embedding 
    \begin{equation*}
        F \hookrightarrow \End_{L}(A) \otimes_{\Z} \Q,
    \end{equation*}
    where $F$ is a number field with $[F : \Q] = \dim(A).$    
\end{defn}

Consider now an abelian variety $A$ over a number field $L$, and let $g$ be its dimension. Assume $A$ is of $\Gltwo$-type, and let $F$ be a number field acting on it by endomorphisms as in Definition~\ref{Def Gltwo}. Fix, once and for all, a prime number $\ell$. Denote by $\TlA$ the $\ell$-adic Tate module of $A$, and $\Vl\coloneqq \TlA \otimes_{\Zl} \Ql$. The absolute Galois group of $L$, denoted by $G_L$, acts on $\Vl$, giving rise to a representation 
\begin{equation*}
    \rhol : G_L \longrightarrow \Aut_{\Ql}(\Vl) \simeq \GL_{2g}(\Ql).
\end{equation*}

On the other hand, $A$ being of $\Gltwo$-type implies that $F$ acts on $\Vl$ linearly, so setting $\Fl \coloneqq F \otimes_{\Q} \Ql$, we see that $\Vl$ carries a $\Fl$-module structure. Moreover, the action of $G_L$ on $\Vl$ is $\Fl$-linear, so $\Vl$ is a $\Fl[G_L]$-module.

For any prime ideal $\lambda$ in $F$ such that $\lambda \mid \ell$, we define $\Flam$ to be the $\lambda$-adic completion of $F$. The local field $\Flam$ carries a structure of $\Fl$-module via the projection induced by the isomorphism of $\Ql$-vector spaces $\Fl \simeq \prod_{\lambda \mid \ell} \Flam$. We define $\Vlam \coloneqq \Vl \otimes_{\Fl} \Flam$, which also carries a structure of $\Fl$-module.
%\newpage

\begin{prop}
    For any prime ideal $\lambda \mid \ell$, $\Vlam$ is a $2$-dimensional $\Flam$-vector space.
\end{prop}

\begin{proof}
    This is a particular case of \cite{Ribet76}*{Theorem 2.1.1}.
\end{proof}

The Galois group $G_L$ also acts $\Fl$-linearly on $\Vlam$. Since $\Flam$ is embedded in $\Fl$, then $G_L$ also acts $\Flam$-linearly on $\Vlam$, giving rise to a representation 
\begin{equation}\label{Rep rholamb}
    \rholam : G_L \longrightarrow \Aut_{\Flam}(\Vlam) \simeq \Gltwo(\Flam).
\end{equation}
In \cite{BCDF} (resp. \cite{chen2022modular}), the authors compute the conductor of the representations $\rholam$ above, in the particular case where $A$ is the Jacobian of the curve $\Crrp$ (resp. $\Cminus, \Cplus$) given in Definition \ref{def:cr} (resp. \ref{def: Cminus & Cplus}), and $F = L = \Q(\zr+\zr^{-1})$. In Theorem \ref{thm:conductorQplus} (resp. \ref{thm: conductorCminusQpl}, \ref{thm: conductorCplusQpl}), we computed the conductor exponents of $\rhol$ at odd primes in the same setting. 
The goal of this section is to show how the conductors of $\rho_{\ell}$ and $\rho_{\lambda}$ relate to each other, and check that the respective computations agree.

\begin{prop}\label{Decomp Vell prop}
    We have an isomorphism of $\Fl[G_L]$-modules 
    \begin{equation}\label{Decomp Vell}
        \Vl \simeq \bigoplus_{\lambda \mid \ell} \Vlam.
    \end{equation}
\end{prop}

\begin{proof}
    Since there are finitely many $\lambda$'s dividing $\ell$, we have an isomorphism of $\Fl$-modules $\Fl \simeq \bigoplus_{\lambda \mid \ell} \Flam$.
    Tensoring with $\Vl$, we obtain an isomorphism of $\Fl$-modules
    \begin{equation*}
        \Vl \simeq \Vl \otimes_{\Fl} \Fl \simeq \Vl \otimes_{\Fl} \left( \bigoplus_{\lambda \mid \ell} \Flam \right) \simeq \bigoplus_{\lambda \mid \ell} \Vl \otimes_{\Fl} \Flam \simeq \bigoplus_{\lambda \mid \ell} \Vlam.
    \end{equation*}
    The isomorphism is explicitly given by $\theta : P \mapsto \left( P \otimes \iota_{\lambda}(1) \right)_{\lambda \mid \ell \, }$, where $\iota_{\lambda} : F \hookrightarrow \Flam$ is the canonical embedding. On any $\Vlam = \Vl \otimes_{\Fl} \Flam$, the Galois group $G_L$ acts on the first component of the tensor product. From this it is clear that $\theta$ is $G_L$-equivariant and therefore establishes an isomorphism of $\Fl[G_L]$-modules.
\end{proof}

The link between the representation $\rhol$ and the $\rholam$'s is as follows. Fix $\lambda \mid \ell$: if one forgets the structure of $\Fl$-module on $\Vlam$ and sees this as a $\Ql$-vector space, then the isomorphism of $\Fl$-modules \eqref{Decomp Vell} becomes an isomorphism of $\Ql$-vector spaces. The difference now is that the dimension of $\Vlam$ as a $\Ql$-vector space is $2 \, [\Flam : \Ql]$. Hence, by forgetting the structure of $\Fl$-module on $\Vlam$, $\rholam$ induces a representation 
\begin{equation}
    \rholamtilde : G_L \longrightarrow \Aut_{\Ql}(\Vlam) \simeq \GL_{2 \, [\Flam : \Ql]}(\Ql),
\end{equation}
and therefore isomorphism \eqref{Decomp Vell} becomes an isomorphism of $\Ql$-representations 
\begin{equation}\label{Decomp rhol}
    \rhol \simeq \bigoplus_{\lambda \mid \ell} \rholamtilde.
\end{equation}
In the next section, we study the relationship between the different $\rholam$'s as $\lambda$ varies.
%%%%%%%%%%%%%%%%%%%%%%%%%%%%%%%%%%%%%%%%%%%%%%%%%%%%%%%%%%%%%%%%%%%%%%%%%%%%%
\subsection{Compatible systems of representations}
%%%%%%%%%%%%%%%%%%%%%%%%%%%%%%%%%%%%%%%%%%%%%%%%%%%%%%%%%%%%%%%%%%%%%%%%%%%%%

Let $K$ be a local field. Recall from \eqref{eq: Ramif groups} the definition of the ramification groups of $G_K$ in upper numbering. We begin with the definition of the conductor of an $\ell$-adic representation. We refer the reader to \cite{Ulmer} for further details.

\begin{defn}
    Let $E$ be a finite extension of $\Ql$ and $V$ be a $E$-vector space. Consider a continuous representation $\rho : G_K \rightarrow \GL(V)$. The conductor $\condexp{\rho}$ of $\rho$ is the integer defined by
    \begin{equation*}
        \condexp{\rho} \coloneqq \int_{-1}^{\infty} \codim V^{G_L^{s}} \,  d s. 
    \end{equation*}
    To ease notation, whenever $L$ is a number field, $v$ is a finite place of $L$ and $\rho : G_L \rightarrow \GL(V)$ is a continuous representation, we let $\condv{\rho} \coloneqq \condexp{\rho |_{G_{L_v}}}$.
\end{defn}

If $C / K$ is a curve and $\rhol : G_K \rightarrow \Aut_{\Ql} (\Vl)$ is the representation attached to the torsion points of the Jacobian of $C$ (\textit{c.f.} Section~\ref{subsec:CJC}), we have $\condexp{C/K} = \condexp{\rhol}$.

\begin{remark}
    When $\rho$ has finite image, it factors through a finite extension of $K$. In this setting, $\condexp{\rho}$ is simply the \textit{Artin conductor} of $\rho$ (see \cite[Chapter VII, §11]{Neukirch}).
\end{remark}

The following definition makes use of Weil-Deligne representations (abbreviated as $\WD$ from now on). We refer the reader to \cite{Ulmer} for a quick survey on the definitions and basic facts about $\WD$-representations, and their link to $\ell$-adic representations. We adopt here the terminology used in \cite{Gee} rather than the one used in \cite{Ribet76}.

\begin{defn}\label{Comp syst def}
    Let $F$ be a number field and $\Pplaces$ its set of finite places. Let $L$ be a number field and fix $S$ a finite set of finite places of $L$. For any $\lambda \in \Pplaces$, let $\Slam$ be the set of finite places $v$ of $L$ such that $v$ and $\lambda$ lie over the same rational prime $\ell$. For any $\lambda \in \Pplaces$, let $\Vlam$ be a $n$-dimensional $\Flam$-vector space where $G_L$ acts continuously. We say that $(\Vlam)_{\lambda \in \Pplaces}$  forms a $F$-rational weakly compatible system of representations of $G_L$ (with exceptional set $S$) if:
    \begin{enumerate}
        \item\label{Syst unram} for all $\lambda \in \Pplaces$, the representation $\Vlam$ is unramified outside $S \cup \Slam$,
        \item for any finite place $v$ of $L$ not in $S$ there exists a polynomial $p_{v}(T) \in F[T]$ such that
        \begin{equation*}
            p_{v}(T) = \det \left( 1 - T \Frob_{v} \mid \Vlam \right) \in F_{\lambda}[T] \quad \forall \lambda \text{ such that } v \notin \Slam.
        \end{equation*}
    \end{enumerate}
    If in addition, the following condition is fulfilled, we say that $(\Vlam)_{\lambda \in \Pplaces}$ forms a strictly compatible system:
    \begin{enumerate}[resume]
        \item for any $v \in S$, there exists a Frobenius semisimple $\WD$-representation $(r_{v}, N_{v})$ of $L_{v}$ such that 
        \begin{equation}\label{WD associated}
            \WD \left(\left. \Vlam\right |_{G_{L_v}} \right)^{\mathrm{ss}} = (r_v, N_v) \quad \forall \lambda \text{ such that } v \notin \Slam .
        \end{equation}
    \end{enumerate}
\end{defn}

Note that if $(\Vlam)_{\lambda \in \Pplaces}$ is a compatible system of representations of $G_L$,  Definition \ref{Comp syst def} $(1)$ implies that for any $\lambda \in \Pplaces$ and any $v \notin S \cup S_{\lambda}$ we have $\condv{\Vlam} = 0$. The next Lemma uses the notion of the conductor of a WD-representation. We refer to \cite[§8]{Ulmer} for details.

\begin{lemma}\label{Cond constant}
    Consider a strictly compatible system of $F$-rational representations $(\Vlam)_{\lambda \in \Pplaces}$ of $G_L$. Then for any $v \in S$, $\condv{\Vlam}$ is independent of $\lambda$.
\end{lemma}

\begin{proof}
    Fix $v \in S$, $\lambda \in \Pplaces$ and consider the $\WD$-representation $(r_v, N_v)$ associated to $\Vlam |_{G_{L_v}}$ as in \eqref{WD associated}. In \cite[Section 8]{Ulmer}, the author shows that $\condv{\Vlam}$ matches the conductor of the $\WD$-representation associated to $\Vlam$. By \eqref{WD associated}, this is the conductor of $(r_v, N_v)$, which is independent of $\lambda$, hence the result.
\end{proof}

We can now state the link between abelian varieties of $\Gltwo$-type and strictly compatible systems.

\begin{theorem}\label{Strict comp ab var}
    Let $L$ be a number field, $A/ L$ be an abelian variety of $\Gltwo$-type and let $F$ be a number field acting on $A$ by endomorphisms. Then the system of representations $(\rholam)_{\lambda \in \Pplaces}$ associated with $A$ as in Section \ref{Gltwo ab var} forms a strictly compatible system with an exceptional set consisting of the places of $L$ where $A$ has bad reduction. 
\end{theorem}

\begin{proof}
    This is a direct consequence of \cite[Proposition 2.8.1]{BCGP}.
    %The $\rholam$'s forming a compatible system is proved in \cite[Section 11.10]{Shimura}. Strict compatibility is established in \cite[Section 2.4]{Fontaine} (see \cite[Section 2.30 - 2.39]{Gee} for a quick survey).
\end{proof}

\begin{remark}
    Strict compatibility of a system $(\rholam)_{\lambda \in \Pplaces}$ is known to be true if the representations $\rholam$ are all modular. In the particular case of \cite{BCDF} (resp. \cite{chen2022modular}), where $A$ is the Jacobian of $\Crrp$ (resp. $\Cminus, \Cplus$) and $L = F = \Q(\zr + \zr^{-1})$ (denoted by $\Qplus$ above), the authors prove that the representations $\rholam$ are modular. By \cite{Carayol}, this gives another proof of the strict compatibility of the system $(\rholam)_{\lambda \in \Pplaces}$. 
\end{remark}

\begin{cor}\label{Cond indep lam}
    With hypothesis as in Theorem \ref{Strict comp ab var}, $\condv{\rholam}$ is independent of $\lambda$.
\end{cor}
\vspace{0.25em}
%%%%%%%%%%%%%%%%%%%%%%%%%%%%%%%%%%%%%%%%%%%%%%%%%%%%%%%%%%%%%%%%%%%%%%%%%%%%%
\subsection{Artin conductors of \texorpdfstring{$\rhol$}{rhol} and \texorpdfstring{$\rholam$}{rholam}}
%%%%%%%%%%%%%%%%%%%%%%%%%%%%%%%%%%%%%%%%%%%%%%%%%%%%%%%%%%%%%%%%%%%%%%%%%%%%%
From now on, we come back to the setting of Section \ref{Gltwo ab var} and we relate here the Artin conductors of the representations $\rhol$ and $\rholam$. Let $S$ be the set of places of bad reduction of $A$. When $v$ is a place of good reduction, the Néron-Ogg-Shafarevich criterion (\cite{SerreTate}) yields $\condv{\rhol} = \condv{\rholam} = 0$.

\begin{proposition}\label{Relation cond}
    Let $\rhol$ and $\rholam$ be as in subsection \ref{Gltwo ab var}, then for any $v \in S$ we have
    \vspace{0.25em}
    \begin{equation}\label{rel conductors}
        \condv{\rhol} = [F : \Q] \, \condv{\rholam} .
    \end{equation} 
\end{proposition}

\begin{proof}
We begin by relating $\condv{\rholam}$ to $\condv{\rholamtilde}$. For any $u \geq -1$, denote by $G^{u}_{L_v} \subseteq G_{L_v}$ the $u$-th ramification group in upper numbering of $G_{L_v}$. By definition of the conductor, we have
\begin{align*}
    \condv{\rholam} & =\int_{-1}^{\infty} \codim_{\Flam} \Vlam^{G^{u}_{L_v}}\, du, \quad \\
    \condv{\rholamtilde} & =\int_{-1}^{\infty} \codim_{\Ql} \Vlam^{G^{u}_{L_v}}\, du. \quad
\end{align*}
For any vector space $V$ over $\Flam$, we have $\codim_{\Ql}(V) = [\Flam : \Ql] \, \codim_{\Flam}(V)$. Applying this to $\Vlam^{G^{u}_{L_v}}$ for $u \geq -1$ yields
    \begin{equation}\label{cond intermed}
        \condv{\rholamtilde} = [\Flam : \Ql] \, \condv{\rholam}.
    \end{equation}
    Since the conductor of a direct sum of representations is the sum of the conductors, isomorphism \eqref{Decomp rhol} combined with \eqref{cond intermed} give
    \begin{equation*}
        \condv{\rhol} = \sum_{\lambda \mid \ell} \condv{\rholamtilde} = \sum_{\lambda \mid \ell} \, [\Flam : \Ql] \, \condv{\rholam}.    
    \end{equation*}
    Corollary \ref{Cond indep lam} shows that $\condv{\rholam}$ is independent of $\lambda$, so using $\sum_{\lambda \mid \ell} \, [\Flam : \Ql] = [F : \Q]$, we obtain the desired result.
    
\end{proof}

We know from \cite[§3]{BCDF} that $\Jac(\Crrp) / \Qplus$ is an abelian variety of $\Gltwo$-type. In the same way, Theorem 1 (resp. Remark p.1058) in \cite{TTV} states that $\Jac(\Cminus) / \Qplus$ (resp. $\Jac(\Cplus) / \Qplus$) is of $\Gltwo$-type too. We deduce now the link between the conductor exponent of the curves above and those of the $2$-dimensional representations arising from the $\Gltwo$-type Jacobians.

\begin{prop}\label{cor: Results are compatible1}
    Let $C \in \lbrace \Crrp, \Cminus, \Cplus \rbrace$, and let $\qQzrplus$ be an odd place of $\Qplus$. Consider $\rho_{\lambda, \, C}$, the representation attached to the Jacobian of $C / \Qplus_{\qQzrplus}$ as in \eqref{Rep rholamb}. Then we have
    \begin{equation*}
        \condexp{C /\Qplus_{\qQzrplus}} = \frac{r-1}{2} \, \mathfrak{n}_{\qQzrplus}(\rho_{\lambda, \, C}).
    \end{equation*}
\end{prop}

\begin{proof}
    This follows from Proposition \ref{Relation cond}, combined with $[\Qplus:\Q]= \text{genus}(C) = \frac{r-1}{2}$.
\end{proof}

\begin{cor}
    Theorem \ref{thm:conductorQplus} is compatible with the results in \cite[§5]{BCDF}. When specializing to $r = 5$, Theorems~\ref{thm: conductorCminusQpl} and \ref{thm: conductorCplusQpl} are compatible with the results in \cite[§7]{chen2022modular}.
\end{cor}

\bibliographystyle{amsplain}
\bibliography{bib}
\end{document}